\documentclass[12pt]{article}

\usepackage[latin1]{inputenc}

   \usepackage[francais,english]{babel}

   \usepackage{graphicx} 

\usepackage{amsmath}    
\usepackage{amssymb}

\input diagrams.sty
\usepackage[all]{xy}

\def\to{\mathchoice{\longrightarrow}{\rightarrow}{\rightarrow}{\rightarrow}}

\def\To#1{\mathchoice{\xrightarrow{\textstyle\kern4pt#1\kern3pt}}{\stackrel{#1}{\longrightarrow}}{}{}}

\newcommand{\CC}{\mathbb{C}}
\newcommand{\QQ}{\mathbb{Q}}
\newcommand{\cH}{\mathcal{H}}

\DeclareMathOperator{\can}{can}
\DeclareMathOperator{\var}{var}
\DeclareMathOperator{\Id}{Id}

\begin{document}
 
\pagestyle{myheadings}

\title{Interaction de strates cons\'ecutives pour les cycles \'evanescents  III : \\ Le cas de la valeur propre 1.}

\author{par Daniel Barlet\footnote{Universit\'e Henri Poincar\'e et Institut Universitaire de France\newline
Institut Elie Cartan (Nancy) UMR 7502 UHP-CNRS-INRIA\newline
BP 239 - F - 54506 Vandoeuvre-l\`es-Nancy Cedex.France.\newline
e-mail : barlet@iecn.u-nancy.fr}, \\ avec un appendice de Claude Sabbah\footnote{Centre de Math\'ematiques Laurent Schwartz,\newline 
Ecole Polytechnique,  UMR 7640 CNRS\newline
e-mail sabbah@math.polytechnique.fr}}

\date{06/09/05}

\maketitle

\markright{Interaction de strates  III.}

\section*{Abstract.}

\noindent This text is a study of the missing case in our article [B.91], that is to say the eigenvalue 1 case. Of course this is a more involved situation because the existence of the smooth stratum for the hypersurface \ $\lbrace f = 0 \rbrace $ \ forces to consider three strata for the nearby cycles. And we already know that the smooth stratum is always "tangled" if it is not alone (see [B.84b] and the introduction of [B.03]). 

\noindent The new phenomenon is the role played here by a "new" cohomology group, denote by \ $ H^n_{c\,\cap\,S}(F)_{=1} $, of the Milnor's fiber of \ $f$ \ at the origin. It has the same dimension as \ $H^n(F)_{=1}$ \  and \ $H^n_c(F)_{=1}$\footnote{The case \ $n = 2$, that is to say the case of an hypersurface in \ $\mathbb{C}^3$,  is analoguous but special.} and it leads to a non trivial factorization of the canonical map \ $ can : H^n_c(F)_{=1} \rightarrow H^n(F)_{=1}$, and to a monodromic isomorphism of variation \ $ var :  H^n_{c\,\cap\,S}(F)_{=1} \rightarrow H^n_c(F)_{=1}$. It gives a canonical hermitian form
 $$ \mathcal{H} : H^n_{c\,\cap\,S}(F)_{=1}\times H^n(F)_{=1} \rightarrow \mathbb{C}$$ 
 which is non degenerate\footnote{ Here again the \ $n = 2$ \ is analoguous but special.}. This generalizes the case of an isolated singularity for the eigenvalue 1 (see [B.90] and [B.97] ).
 
 \noindent The "overtangling" phenomenon for strata associated to the eigenvalue 1 implies the existence of triple poles at negative integers (with big enough absolute value) for the meromorphic continuation of the distribution \ $ \int_X \vert f \vert^{2\lambda} \Box $ \  for  functions \ $f$ \ having semi-simple local monodromies at each singular point of \ $ \lbrace f = 0 \rbrace .$

\bigskip

\noindent {\bf AMS  Classification  2000} : 32-S-25, 32-S-40, 32-S-50.

\bigskip

\noindent {\bf Key Words} : Hypersurface, Non isolated Singularity, Vanishing Cycles,Tangling of Strata.

\bigskip

\section* { R\'esum\'e.}
 
 \bigskip
 
 \noindent Ce texte \'etudie le cas manquant de notre article [B.91], \`a savoir le cas de la valeur propre 1. C'est \'evidemment un cas plus compliqu\'e que celui d'une valeur propre  $\not= 1$ \ puisque la pr\'esence de la strate des points lisses de l'hypersurface \ $\lbrace f = 0 \rbrace $ \ oblige \`a consid\'erer trois strates pour les cycles proches. Et l'on sait d\'ej\`a que cette strate lisse est toujours "emm\^el\'ee" en pr\'esence d'une autre strate (voir [B.84b] et l'introduction de [B.03]).
 
 \noindent Le ph\'enom\`ene nouveau est le r\^ole jou\'e  ici par un "nouveau" groupe de cohomologie, not\'e \ $ H^n_{c\,\cap\,S}(F)_{=1} $, de la fibre de Milnor de \ $f$ \ \`a l'origine, qui est de m\^eme dimension que \ $H^n(F)_{=1}$ \ et \ $H^n_c(F)_{=1}$\footnote{Le cas \ $n = 2 $, c'est \`a dire le cas d'une hypersurface de \ $\mathbb{C}^3$,  est analogue mais sp\'ecial.} et qui donne lieu \`a une factorisation non triviale de l'application canonique \ $ can : H^n_c(F)_{=1} \rightarrow H^n(F)_{=1}$, et \`a un isomorphisme monodromique de variation \ $ var :  H^n_{c\,\cap\,S}(F)_{=1} \rightarrow H^n_c(F)_{=1}$.  On en d\'eduit une forme hermitienne canonique 
 $$ \mathcal{H} : H^n_{c\,\cap\,S}(F)_{=1}\times H^n(F)_{=1} \rightarrow \mathbb{C}$$ 
 qui est non d\'eg\'en\'er\'ee\footnote{ La encore le cas \ $n = 2$ est analogue mais sp\'ecial.}. Ceci g\'en\'eralise le cas d'une singularit\'e isol\'ee pour la valeur propre 1 (voir [B.90] et [B.97] ).
 
 \noindent  Le ph\'enom\`ene de "suremm\^element" de strates pour la valeur 1 donne l'existence de p\^oles triples aux entiers n\'egatifs (assez grands en valeur absolue) pour le prolongement m\'eromorphe de la distribution \ $ \int_X \vert f \vert^{2\lambda} \Box $ \ en pr\'esence de monodromies locales semi-simples pour \ $f$ \ en chaque point du lieu singulier de \ $ \lbrace f = 0 \rbrace .$
 
 \bigskip

\noindent {\bf AMS  Classification  2000}: 32-S-25, 32-S-40, 32-S-50.
 
 \bigskip
 
 \noindent {\bf Mots Clefs} : Hypersurface, Singularit\'e non isol\'ee, Cycles \'Evanescents, Emm\^element de Strates.

\section*{Introduction.}


\subsection{}
Le but du pr\'esent article est de compl\'eter l'\'etude du ph\'enom\`ene \\ d'emm\^element de strates  (cons\'ecutives) pour les cycles \'evanescents, faite dans [B.91]\footnote{Voir \'egalement [B.03].} pour une valeur propre de la monodromie diff\'erente de  1 par l'\'etude du  cas de la valeur propre 1 de la monodromie. La situation que nous consid\'erons ici est le cas d'un germe non constant de fonction holomorphe \ $ \tilde{f} : (\mathbb{C}^{n+1}, 0) \rightarrow (\mathbb{C}, 0)$ \ v\'erifiant la propri\'et\'e suivante :

\noindent Il existe un germe de courbe \ $(S, 0)$ \ \`a l'origine de \ $\mathbb{C}^{n+1}$ \ tel qu'en chaque point \ $x$ \ en dehors de cette courbe la monodromie locale de \ $f$ \ agissant sur la cohomologie r\'eduite de la fibre de Milnor de \ $f$ \ en \ $x$ \ ne pr\'esente pas la valeur propre 1. 

\noindent Bien sur, le germe de  courbe  \ $(S, 0)$ \ est  "\`a priori" contenu dans le lieu singulier de l'hypersuface \ $f^{-1}(0)$ \ mais ce lieu singulier sera en g\'en\'eral  de dimension plus grande que    1.

\bigskip


\subsection{}
 De m\^eme que pour le cas d'une valeur propre \ $\not= 1$ \ on rencontre cette situation de la fa{\c c}on suivante : 

\noindent On consid\`ere un germe non constant de fonction holomorphe 
$$ \tilde{F} :  (\mathbb{C}^{N}, 0) \rightarrow (\mathbb{C}, 0)$$
et on denote par \ $\Sigma_0 \subset \Sigma_1 $ \ les deux plus grosses strates du lieu singulier de \ $ \tilde{F}^{-1}(0) $ \ le long desquelles la valeur propre  1  de la monodromie apparait (sur la cohomologie \underline{r\'eduite}\footnote{Ce sont donc les supports des deux premiers faisceaux de cohomologie non nuls du complexe \`a cohomologie constructible  des cycles \'evanescents de \ $\tilde{F}$ \ relatif \`a la valeur propre  1  de la monodromie, que l'on notera \ $ \psi_{\tilde{F}, 1}.$}). Supposons que la codimension de \ $\Sigma_0$ \ soit \'egale \`a \ $n+1$ \ et que celle de \ $\Sigma_1$ \ soit \'egale \`a \ $n$\footnote{C'est dans cette hypoth\`ese, qui est "g\'en\'erique", que l'adjectif "cons\'ecutif" prend son sens dans l'expression "strates cons\'ecutives": \ $\Sigma_0$ \ est de codimension  1  dans \ $\Sigma_1$.}. Effectuons une section plane transverse en un point g\'en\'erique de \ $\Sigma_0 $ \ par un plan  \ $P$ \ de dimension \ $n+1$. Un tel plan "g\'en\'erique" est non caract\'eristique pour le \ $D-$module correspondant au faisceau pervers \ $\psi_{\tilde{F, 1}}$ \ et la restriction \ $ \tilde{f} : = \tilde{F}\vert_P $ \ v\'erifie nos hypoth\`eses et d\'ecrit la situation locale pour \ $\psi_{\tilde{F, 1}}$ \ dans un voisinage du point g\'en\'erique de \ $\Sigma_0 $ \ que l'on a consid\'er\'e\footnote{Grace au th\'eor\`eme de restriction non carat\'eristique (voir [B.M.89]) et \`a la correspondance de Riemann Hilbert (voir [K.84]).}.

\bigskip

\subsection{Les hypoth\`eses standards pour la valeur propre 1.}

\bigskip

\noindent Soit \ $X$ \ un voisinage ouvert de Stein connexe de l'origine dans \ $\mathbb{C}^{n+1}$ \ et supposons que \ $ n \geq 2 $\footnote{Dans le cas \ $n=1$ \ le ph\'enom\`ene d'emm\^element de strates pour la valeur propre 1  que nous allons \'etudier n'appara\^it pas, m\^eme pour une fonction non r\'eduite. Le cas o\`u  \ $n=1$ \ sera consid\'er\'e dans la suite  sous l'hypoth\`ese plus forte d'une singularit\'e isol\'ee pour la valeur propre  1 (notion d\'ej\`a \'etudi\'ee dans [B.90]  et [B.97] ). Ce cas sera essentiel quand on \'etudiera la situation d'une section hyperplane transverse en un point g\'en\'erique de la courbe \ $S$.}. Soit \ $ f : X \rightarrow \mathbb{C}$ \ une fonction holomorphe non constante. Notons pour \ $p \geq 0$ \ par \ $H^p(0)$ \ le faisceau constructible sur \ $ Y : = f^{-1}(0)$ \ donn\'e par le sous-faisceau spectral de la monodromie de \ $f$ \ pour la valeur propre  1, agissant sur le \ $p-$i\`eme faisceau de cohomologie du complexe \ $\psi_f$ \ des cycles \'evanescents. Donc le germe en chaque point \ $y \in Y$ \ de \ $H^p(0)$ \ est le \ sous-espace spectral pour la valeur propre 1 de la monodromie agissant sur le $p-$i\`eme groupe de la cohomologie \underline{r\'eduite} de la fibre de Milnor de \ $f$ \ en \ $y$.

\noindent Nous ferons les hypoth\`eses suivantes :
\begin{itemize}
\item{(i)}  Le faisceau \ $H^n(0)$ \ est concentr\'e \`a l'origine.
\item{(ii)}  Le faisceau \ $H^{n-1}(0) $ \ est concentr\'e sur une courbe \ $S$ \ admettant l'origine pour unique point singulier. De plus, chaque composante irr\'eductible de \ $S$ \ contient l'origine et est un disque topologique.
\item{(iii)} La restriction de \ $H^{n-1}(0) $ \ \`a \ $S^* : = S \setminus \lbrace 0 \rbrace$ \ est un syst\`eme local.
\item{(iv)} Pour tout \ $ i \geq 2$ \ on suppose que \ $ H^{n-i}(0) \equiv 0 .$
\end{itemize}

\bigskip


\subsubsection{ Remarques.}
\begin{itemize}
\item{} Seule la condition (iv)  est r\'eellement restrictive dans notre situation locale d'un germe \`a l'origine de \ $\mathbb{C}^{n+1}$ \ en raison de la perversit\'e de \ $ \psi_{f,1}.$
\item{} Bien sur, ces hypoth\`eses sont satisfaites (quitte \`a se restreindre \`a un voisinage ouvert de l'origine assez petit) d\`es que le lieu singulier est de dimension inf\'erieur ou \'egal \`a  1. Pour \ $n=2$ \  c'est le cas d\`es que l'on consid\`ere un germe non constant de fonction holomorphe qui est r\'eduit.
\item{} Pour \ $n = 3$ \ ces hypoth\`eses seront  satisfaites (quitte \`a nouveau \`a  se restreindre ...) pour un germe r\'eduit   d\`es que les singularit\'es transverses aux composantes irr\'eductibles de dimension 2 du lieu singulier  sont des courbes r\'eduites et  irr\'eductibles\footnote{En effet, pour une courbe plane r\'eduite et   irr\'eductible, la valeur propre 1  de la monodromie n'apparait pas.}.
\end{itemize}

\bigskip


\subsection{Les r\'esultats.}
 D\'ecrivons les principaux r\'esultats que nous obtenons. La nouveaut\'e principale par rapport au cas d'une valeur propre diff\'erente de \ $1$, trait\'e dans [B.91], est l'apparition d'une "cohomologie interm\'ediaire" pour la partie spectrale relative \`a la valeur propre \ $1$ \ de la monodromie, pour la fibre de Milnor de \ $f$ \ \`a l'origine. Nous la noterons par \ $H^n_{c\,\cap\,S}(0)$ \ et elle est introduite au paragraphe 4. Elle correspond \`a des cycles \'evanescents dont la limite \`a l'origine  rencontre  le sous-ensemble analytique \ $S$ \ suivant un compact (mais la limite peut rencontrer le bord de \ $f^{-1}(0)$ \ loin de \ $S$). Ceci est pr\'ecis\'e dans l'appendice I.

\noindent On a alors deux applications naturelles d'\'elargissement de support
$$ H^n_c(0) \overset{can_{c}}{\rightarrow} H^n_{c\,\cap\,S}(0) \overset{can_{c\,\cap\,S}} {\rightarrow} H^n(0) $$
dont la compos\'ee est l'application canonique usuelle
$$ can : H^n_c(0)\rightarrow H^n(0) .$$
Dans le cas d'une singularit\'e isol\'ee pour la valeur propre 1, l'application \ $ can_{c\,\cap\,S} : H^n_{c\,\cap\,S}(0) \rightarrow H^n(0) $ \ est un isomorphisme. Ce n'est plus le cas sous les hypoth\`eses standards pr\'ecis\'ees plus haut et nous dirons que le ph\'enom\`ene de {\bf suremm\^element pour la valeur propre 1} se produit dans cette situation pr\'ecis\'ement quand cette application n'est pas un isomorphisme.

En suivant la m\^eme ligne de preuve que dans [B.91] nous commencerons par \'etablir le 

\bigskip

\subsubsection{Th\'eor\`eme 1.}

\smallskip

\noindent {\it Sous les hypoth\`eses standards pour la valeur propre  1, on a, pour \ $n\geq 3$, la suite exacte : }

$$ 0 \rightarrow H^1_{\lbrace 0 \rbrace}(S, H^{n-1}(0)) \overset{i}{\rightarrow} H^n_{c\,\cap\, S}(0) \overset{can_{c\,\cap\,S}}{\longrightarrow} H^n(0) \overset{\theta}{\rightarrow} H^1(S^*, H^{n-1}(0)) .$$

\bigskip

\noindent Ce th\'eor\`eme est prouv\'e au paragraphe 4, ainsi qu'un r\'esultat analogue pour \ $n = 2$.

\noindent Maintenant l'existence d'une forme hermitienne localement constante et non d\'eg\'en\'er\'ee sur le syst\`eme local \ $H^{n-1}(0)$ \ d\'eduite de la forme hermitienne canonique des sections hyperplanes transverses \`a \ $S^*$ \ qui pr\'esentent une singularit\'e isol\'ee pour la valeur propre \ $1$ \ (voir [B.90]), donne l'in\'egalit\'e\footnote{Puisque \ $H^{n-1}(0)$ \ n'a pas de section non nulle \`a support l'origine.}
$$ \dim H^1_{\lbrace 0 \rbrace}(S, H^{n-1}(0)) \leq \dim H^1(S^*, H^{n-1}(0)) .$$
On en d\'eduit  l'in\'egalit\'e
$$ \dim H^n_{c\,\cap\,S}(0) \geq \dim H^n(0) .$$
Le r\'esultat le plus difficile sera de d\'emontrer que l'on a, en fait, toujours \'egalit\'e de ces deux dimensions sous les hypoth\`eses standards, pour  \ $n \geq 3 $. 

\subsubsection{Th\'eor\`eme 2.}

\textit{Sous les hypoth\`eses standards pour la valeur propre 1 et  pour \ $n \geq 3 $ \ on a les propri\'et\'es suivantes
\begin{itemize}
\item{1)} $\dim H^n_{c\,\cap\,S}(0) = \dim H^n(0) $.
\item{2)} Il existe une application naturelle lin\'eaire et monodromique \\
$ \tilde{var} : H^n_{c\,\cap\,S}(0) \rightarrow H^n_c(0) $ \ qui est un isomorphisme d'espace vectoriels monodromiques.
\item{3)} Il existe une  forme hermitienne "canonique"
$$ \mathcal{H} :  H^n_{c\,\cap\,S}(0)\times H^n(0) \rightarrow \mathbb{C} $$
 et elle est non d\'eg\'en\'er\'ee.
\end{itemize}
Pour \ $n = 2 $ \ tout ceci reste vrai \`a condition de remplacer \ $H^2_{c\,\cap\,S}(0)$ \ par son quotient par \ $j(H^1(S^*, \mathbb{C}))$, o\`u \ $j$ \ est une injection naturelle qui sera d\'efinie au 4.12 .}

\bigskip

\noindent Notre preuve passe par la construction d'une application de variation, compatible aux monodromies
$$ var : H^n_{c\,\cap\,S}(0) \rightarrow H^n_c(0) \ ,$$
construction qui sera men\'ee \`a bien au paragraphe 6. Nous prouverons ensuite au paragraphe 8 que cette variation est injective. Ceci impliquera l'\'egalit\'e cherch\'ee puisque la dualit\'e de Poincar\'e sur la fibre de Milnor \`a l'origine donne l'\'egalit\'e des dimensions de \ $H^n_c(0)$ \ et \ $H^n(0)$.

\noindent L'injectivit\'e de la variation (qui est donc finalement un isomorphisme gr\^ace \`a l'\'egalit\'e des dimensions), donne \'egalement la non d\'eg\'en\'erescence de la forme hermitienne\footnote{Elle est sesquilin\'eaire, et \ $( id\times can_{c\,\cap\,S})^*(\mathcal{H})$ \ est hermitienne sur \ $ H^n_{c\,\cap\,S}(0) \times  H^n_{c\,\cap\,S}(0) $.} canonique g\'en\'eralis\'ee :
$$ \mathcal{H} :  H^n_{c\,\cap\,S}(0) \times H^n(0) \rightarrow \mathbb{C} .$$
qui est compatible aux monodromies et reli\'ee, comme dans le cas d'une singularit\'e isol\'ee (voir [Loe. 86]) et dans le cas d'une singularit\'e isol\'ee pour la valeur propre \ $1$ \ (voir [B.97]), \`a la dualit\'e de Poincar\'e (hermitienne), not\'ee \ $ < , > $, par la formule
$$ \mathcal{H}(e,e') = < \tilde{var}(e) , e' > \quad \forall (e,e') \in  H^n_{c\,\cap\,S}(0) \times H^n(0) $$
o\`u \ $ \tilde{var} = var\circ \Theta $ , l'automorphisme \ $\Theta$ \ de \ $H^n_{c\,\cap\,S}(0)$ \ \'etant d\'efini par 
$$\Theta : = \sum_{k=0}^{\infty} (-1)^k \frac{(T-1)^k}{k+1} = \frac{Log\, T}{T} $$
o\`u \ $T$ \ d\'enote la monodromie. Ceci est montr\'e au paragraphe 7.

\noindent Nous concluons en \'etendant les th\'eor\`emes 13 et 14 de [B.91] au cas de la valeur propre 1. On obtient ainsi les deux th\'eor\`emes suivants qui, je l'esp\`ere, seront suffisants pour donner une approche "filtr\'ee" (donc non localis\'ee en \ $f$) du ph\'enom\`ene d'emm\^element de strates  dans le cas d'une fonction holomorphe \`a lieu singulier de dimension 1, au moins dans le cas consid\'er\'e dans  [B.04].

\bigskip

\subsubsection{Th\'eor\`eme 3.}

\textit{Sous les hypoth\`eses standards pour la valeur propre 1, consid\'erons \ $e \in H^n(0) $ \ v\'erifiant \ $\mathcal{N}^k(e) = 0 $ \ avec \ $k \geq k_0$\footnote{Rappelons que, par d\'efinition, \ $k_0$ \ est l'ordre de nilpotence de la monodromie agissant sur le syst\`eme local \ $H^{n-1}(0)$ \ sur \ $S^*$.}. Alors les conditions suivantes sont \'equivalentes :}
\begin{itemize}
\item{1)} \ $\tilde{ob}_k(e) = 0$ .
\item{2)}\textit{ Si \ $w \in \Gamma(Y, Ker\,\delta^n(k))$ \ v\'erifie \ $ r^n(k)(w) = e $, pour chaque \ $j \in \mathbb{Z} $, les fonctionnelles analytiques
$$ P_{k+l} \big(\lambda = 0, \int_X \vert f \vert^{2\lambda} \bar{f}^{-j}\frac{df}{f}\wedge \bar{w}_k \wedge \square \big) $$
sont nulles pour \ $l \geq 2$.}
\end{itemize}

\noindent\textit{ De plus, pour v\'erifier 2) il suffit de le faire pour \ $l = 2$ \ et pour \ $j \in [0,n] $.}

\bigskip

\subsubsection{Th\'eor\`eme 4.}

\textit{Sous les hypoth\`eses standards pour la valeur propre 1, supposons que le prolongement m\'eromorphe de \ $ \int_X \vert f \vert^{2\lambda} \square $ \ admette un p\^ole d'ordre \ $ \geq k $ \ en un entier n\'egatif, avec \ $ k \geq \sup{(k_0, k_1)} + 2$ \ o\`u \ $k_0$ \ et \ $k_1$ \ sont respectivement les ordres de nilpotence de la monodromie agissant sur le syst\`eme local \ $ H^{n-1}(0)_{\vert S^*} $ \ et sur l'espace vectoriel \ $H^n(0)$.}

\noindent \textit{Alors l'ordre des p\^oles aux entiers n\'egatifs assez grand est exactement \'egal \`a \ $k$ \ et l'on a \ $ k_0 \leq k_1 = k - 2 $ \ avec  \ $\tilde{ob}_{k_1} \not\equiv 0 .$}

\bigskip

\noindent L'\'etude de "r\'eciproques" analogues aux th\'eor\`emes 15 et 16  de [B.91] n'est pas abord\'ee ici. Ces questions ne semblent pas du tout triviales \`a r\'esoudre pour la valeur propre 1. D'ailleurs les r\'esultats de ce type dans [B.91] sont probablement les plus d\'elicats de cet article.

\noindent Nous n'avons pas, non plus, abord\'e ici les raffinements de [B.91] qui sont obtenus dans [B.03].\\

\noindent Nous terminons cette \'etude par deux appendices.\\
Le premier donne une description topologique na{\"i}ve (\`a la Milnor) du groupe \ $H^n_{c \cap S}(0)$. Le second, \'ecrit par Claude Sabbah, donne une interpr\'etation en terme des complexes de cycles \'evanescents et de cycles proches des ph\'enom\`enes topologiques sous-jacents \`a cet article.

\bigskip

\newpage

\noindent Voici le plan de cet article

\begin{itemize}
\item{1.} D\'efinitions et rappels.
\item{2.} Le cas d'une singularit\'e isol\'ee pour la valeur propre 1.
\item{3.} La forme hermitienne canonique sur le syst\`eme local \ $ H^{n-1}(0)\vert_{S^*}.$
\item{4.} L'espace vectoriel monodromique \ $H^n_{c\,\cap\,S}(0).$
\item{5.} La suite exacte longue de J. Leray g\'en\'eralis\'ee.
\item{6.} Construction de \ $\tilde{var} : H^n_{c\,\cap\,S}(0) \rightarrow H^n_c(0).$
\item{7.} La forme hermitienne canonique \ $\mathcal{H} : H^n_{c\,\cap\,S}(0)\times H^n(0) \rightarrow \mathbb{C}.$
\item{8.} Injectivit\'e de la variation pour \ $n \geq 3$. Le cas \ $ n = 2$.
\item{9.} Applications et un exemple pour conclure.
\item{Appendice I :} Description topologique de \ $H^n_{c\,\cap\,S}(0).$
\item{Appendice II :} Interpr\'etation en termes de complexes de cycles proches et \'evanescents {\em par Claude Sabbah}.
\end{itemize}

\section{ D\'efinitions et rappels.}

\subsection{}

Sous les hypoth\`eses standards, introduisons les complexes de faisceaux \  $(\Omega^{\bullet}(k), \delta_0 )$ \ et \ $(\mathcal{E}^{\bullet}(k),\delta_0^{\bullet} )$.  Rappelons  que \ $\Omega^{\bullet}$ (resp. $\mathcal{E}^{\bullet} $ ) d\'esigne le faisceau des formes m\'eromorphes (resp. semi-m\'eromorphes) de degr\'e  $^{\bullet}$  sur \ $X$ \ \`a p\^oles dans \ $Y$ \ restreint (topologiquement) \`a \ $Y$, que 
 $$ \Omega^{\bullet}(k) : =  \Omega^{\bullet}  \otimes \mathbb{C}^k , \quad \mathcal{E}^{\bullet}(k) : = \mathcal{E}^{\bullet} \otimes \mathbb{C}^k $$ 
 et que 
$$\delta_0(w) : = dw - \frac{df}{f}\wedge\  _kN(w) $$
 o\`u \ $ _kN \in End(\mathbb{C}^k ) $ \ est d\'efini dans la base canonique \ $ e_1, \cdots , e_k $ \ par \\   $ _kN (e_j) = e_{j+1} ,  j \in [1,k] $ \ avec la convention \ $ e_{k+1} = 0 $ ; cet endomorphisme  agit sur \ $  \Omega^{\bullet}(k) $ \\ (resp. sur \ $\mathcal{E}^{\bullet}(k)$ ) \ par \ $_k\mathcal{N} : = Id \otimes\,  _kN$. 
 
 \noindent Ceci se "visualise" de la fa{\c c}on suivante :

$$ \delta_0(w) = \delta_0 \begin{pmatrix} {w_k} \\ {\vdots }\\ {w_2} \\ {w_1} \end{pmatrix} = \begin{pmatrix}{dw_k - \frac{df}{f}\wedge w_{k-1}} \\  {\vdots }\\ {dw_2 -\frac{df}{f}\wedge w_1} \\{dw_1} \end{pmatrix} .$$
De plus l'inclusion \'evidente de \ $(\Omega^{\bullet}(k), \delta_0 )$ \ dans \  $(\mathcal{E}^{\bullet}(k),\delta_0^{\bullet} )$ \ est un quasi-isomorphisme (voir [B.91], prop.1, p. 444.)

\bigskip

\subsection{}
On a les morphismes de complexes de degr\'e  0  suivants, pour tout  couple \ $(k, k')$ \ d'entiers :
\begin{align*} 
\pi_{k+k', k} :  (\mathcal{E}^{\bullet}(k+k'), \delta_0^{\bullet}) \longrightarrow (\mathcal{E}^{\bullet}(k), \delta_0^{\bullet}) \\
j_{k, k+k'}    :   (\mathcal{E}^{\bullet}(k), \delta_0^{\bullet}) \longrightarrow  (\mathcal{E}^{\bullet}(k+k'), \delta_0^{\bullet}) \\
\end{align*}
donn\'es par les projections et  injections \'evidentes :
\begin{align*}
 \mathbb{C}^{k+k'} \simeq \mathbb{C}^k \oplus \mathbb{C}^{k'} \rightarrow \mathbb{C}^k  \\
 \mathbb{C}^k  \hookrightarrow  \mathbb{C}^{k+k'} \simeq \mathbb{C}^k \oplus \mathbb{C}^{k'} .
 \end{align*}
 Ils se "visualisent" de la fa{\c c}on suivante :
$$ \pi_{k+k', k}\begin{pmatrix} {w_{k+k'}}\\ {\vdots}\\ {w_k}\\ {\vdots} \\ {w_1} \end{pmatrix} = \begin{pmatrix} {w_k} \\ {\vdots} \\ {w_1} \end{pmatrix} \quad \quad  j_{k, k+k'} \begin{pmatrix} {w_k} \\ {\vdots} \\ {w_1} \end{pmatrix} = \begin{pmatrix} {w_k} \\ {\vdots} \\ {w_1} \\ {0} \\ {\vdots} \\ {0}  \end{pmatrix} ,$$
o\`u le dernier vecteur colonne comprend \ $k'$ \ z\'eros.

\subsection{}
 On d\'efinit \'egalement, pour \ $k \geq 1 $, le morphisme  de complexe de degr\'e \ $0$  :
$$ _k \mathcal{N} :  (\mathcal{E}^{\bullet}(k), \delta_0^{\bullet}) \rightarrow  (\mathcal{E}^{\bullet}(k), \delta_0^{\bullet}) $$
en posant  \ $ _k \mathcal{N}  = j_{k-1,k}\circ \pi_{k,k-1} $\ ; ce  qui se "visualise" comme suit :
$$ _k \mathcal{N}\begin{pmatrix} {w_k}\\ {\vdots} \\ {w_2} \\{w_1} \end{pmatrix} = \begin{pmatrix} {w_{k-1}}\\ {\vdots} \\ {w_1} \\ {0} \end{pmatrix} .$$
On a, bien sur, \ $ (_k \mathcal{N})^k = 0 $ \ sur \ $\mathcal{E}^{\bullet}(k).$ La monodromie est alors d\'efinie sur ce complexe par \ $ T : = exp(-2i\pi. _k\mathcal{N}).$

\bigskip

\subsection{}
Enfin nous utiliserons \'egalement les morphismes  de complexes\footnote{ En fait on a \ $ \tau_{k,k'}\circ \delta_0 + \delta_0 \circ \tau_{k,k'} = 0.$}  de degr\'e +1
$$ \tau_{k,k'} :  (\mathcal{E}^{\bullet}(k), \delta_0^{\bullet}) \rightarrow  (\mathcal{E}^{\bullet}(k'), \delta_0^{\bullet}) $$
donn\'es par :
$$ \tau_{k,k'}\begin{pmatrix} {w_k}\\ {\vdots} \\ {w_2} \\ {w_1} \end{pmatrix} = \begin{pmatrix}{0} \\ {\vdots} \\ {0} \\ {\frac{df}{f}\wedge w_k} \end{pmatrix} \Bigg\rbrace k' .$$

\subsection{}
Notons par \ $h^i(k)$ \ le i-\`eme  faisceau de cohomologie du complexe \ $(\mathcal{E}^{\bullet}(k),\delta_0 )$. Le morphisme de complexe
$$ r^{\bullet} :  (\Omega^{\bullet}(k), \delta_0 ) \rightarrow (\Omega^{\bullet}_{X/f}[f^{-1}], d_{/f})\vert_{Y}  $$
d\'efini par \ $ \tilde{r}^{\bullet}(w) = w_k $ \ pour \ $ w \in  \Omega^{\bullet}(k) $, donne la suite exacte courte suivante, pour tout \ $p \geq 1 $  (voir la prop.1 p. 427-428 de [B.91] et la proposition (1.9) pour le lien entre \ $\tilde{r}^{\bullet}$ \ et \ $r^{\bullet}$) :
\begin{equation}
0 \rightarrow H^{p-1}(0)\big/ Im \mathcal{N}^k_{p-1} \rightarrow h^p(k) \overset{r^p(k)}{\rightarrow} Ker \mathcal{N}^k_p \rightarrow 0 
\end{equation}
o\'u \ $ -2i\pi\mathcal{N}_q $ \ d\'esigne le logarithme nilpotent de la monodromie de \ $f$ \ agissant\footnote{de fa{\c c}on unipotente.} sur le faisceau \ $ H^q(0) $ \ pour \ $ q \geq 0 $. On prendra garde que pour \ $q=0$ \ on consid\`ere la cohomologie r\'eduite (donc les "vrais" cycles \'evanescents et non les cycles proches.)

\bigskip

\subsection{}
 Sous les hypoth\`eses standards la suite exacte (1)   montre que les seuls faisceaux de cohomologie non nuls du complexe \ $(\mathcal{E}^{\bullet}(k),\delta_0 )$ \ sont les \ $h^i(k)$ \ pour \ $ i =0,1,n-1,n,n+1 .$ 

\noindent Supposons d'abord \ $ n \geq 3$. On a alors des isomorphismes monodromiques \ $\underline{\mathbb{C}}_Y \simeq h^0(k) $ \ et  \ $\underline{\mathbb{C}}_Y \simeq h^1(k) $ \ donn\'es par 
$$ \lambda \in \underline{\mathbb{C}}_Y \rightarrow \begin{pmatrix} {\lambda} \\ {0} \\ {.} \\ {.} \\ {.} \\ {0} \end{pmatrix} \quad {\rm et} \quad  \lambda \in \underline{\mathbb{C}}_Y \rightarrow \begin{pmatrix} {0} \\ {.} \\ {.} \\ {.} \\ {0} \\ { \lambda \frac{df}{f} } \end{pmatrix} .$$
Le morphisme \ $\tau_k $ \ induit l'isomorphisme \'evident .

\noindent De plus, grace \`a la nullit\'e de \ $H^{n-2}(0)$, on a un isomorphisme
$$ r^{n-1}(k) : h^{n-1}(k) \rightarrow Ker\,\mathcal{N}_{n-1}^k \subset H^{n-1}(0) \quad \forall k \geq 1 $$
ainsi que la suite exacte pour \ $k' \gg 1$
$$ 0 \rightarrow h^{n-1}(k') \overset{\tau_{k',k}}{\rightarrow} h^n(k) \rightarrow Ker\,\mathcal{N}_n^k \rightarrow 0 \quad \forall k \geq 1 .$$
Comme les faisceaux \ $ h^{n-1}(0) $ \ et \ $h^n(0)$ \ sont \`a  support dans \ $S$ \ et que \ $H^n(0)$ \ est \`a support l'origine, on en d\'eduit que sur \ $ S^*$ \ on a un isomorphisme pour \ $ k \gg 1 $
$$ \tau_k : h^{n-1}(k)\vert_{S^*} \simeq h^n(k)\vert_{S^*} .$$
Alors, la suite exacte pour \ $ k \gg 1$ 
$$ 0 \rightarrow H^{n-1}(0) \overset{i^{n-1}}{\rightarrow} h^n(k) \overset{\tau_k}{\rightarrow}h^{n+1}(k) \rightarrow 0 $$
montre que \ $h^{n+1}(k) $ \ est \`a support l'origine et isomorphe \`a \ $H^n(0)$ (pour \ $k \gg1$).

\bigskip

\subsection{}
Pour \ $ n = 2 $, on a la suite exacte
$$ 0 \rightarrow h^0(k) \overset{\tau_k}{\rightarrow} h^1(k) \overset{r^1(k)}{\rightarrow} Ker\, \mathcal{N}_1^k \rightarrow 0 \quad \forall k \geq 1 .$$
Comme \ $\tau_k $ \ induit un isomorphisme sur \ $Y \setminus S$ \ ( et les faisceaux \ $h^i(k)\vert_{Y\setminus S} $ \ sont isomorphes \`a \ $ \underline{\mathbb{C}}_{Y\setminus S} , \ i = 0,1$ ), on aura dans ce cas des isomorphismes monodromiques \ $ h^1(k)\big/ \tau_k h^0(k) \simeq Ker\,\mathcal{N}_1^k \subset H^1(0) $ \ ainsi que \ $ h^0(k) \simeq  \underline{\mathbb{C}}_{Y} .$

\noindent Pour \ $n = 2$ \ on remplacera l'isomorphisme \ $ r^{n-1}(k) : h^{n-1}(k) \rightarrow H^{n-1}(0) $ \ d\^u \`a la nullit\'e de l'application  \ $\tau_k : h^{n-2}(k) \rightarrow h^{n-1}(k) $ \ que l'on a pour \ $n \geq 3$ \ par l'isomorphisme  
$$ \underset{k \rightarrow \infty}{lim} \ h^1(k) \rightarrow H^1(0) $$ 
qui r\'esulte de l'isomorphisme monodromique (plus concret) pour \ $k \geq  k_0$
$$ r^{1}(k+1) : j_{k,k+1}(h^{1}(k)) \rightarrow H^1(0). $$
\noindent Comme l'application \ $\tau_k$ \ n'est pas compatible \`a la limite inductive sur \ $k$, il sera plus commode d'utiliser directement l'isomorphisme \ $ \underset{k \rightarrow \infty}{lim} \ h^1(k) \rightarrow H^1(0) $ \ que la suite exacte 
$$ 0 \rightarrow h^0(k) \overset{\tau_k}{\rightarrow} h^1(k) \overset{r^1(k)}{\rightarrow} Ker\, \mathcal{N}_1^k \rightarrow 0 \quad \forall k \geq 1 .$$

Ce point est important puisqu'il permettra de traiter le cas \ $ n = 2$ \ de fa{\c c}on assez analogue au cas "g\'en\'eral" \ $ n \geq 3 $.

\bigskip

\subsection{} 
On d\'eduit alors facilement de la constructibilit\'e des faisceaux \ $ h^i(k) $ \ et de la contractibilit\'e de \ $Y$ la propri\'et\'e suivante\footnote{Nous omettrons dor\'enavant l'indice de la diff\'erentielle \ $ \delta = \delta_0$.} :
\begin{equation*}
 \mathbb{H}^n(Y, (\mathcal{E}^{\bullet}(k), \delta^{\bullet})) \big/ \tau_k \big(\mathbb{H}^{n-1}(Y, (\mathcal{E}^{\bullet}(k), \delta^{\bullet}))\big)\quad  \simeq Ker\, \mathcal{N}_n^k \subset H^n(0) 
\end{equation*}

\bigskip

\noindent Nous allons \'etablir directement une assertion analogue pour les supports compacts ; cette approche directe s'applique \'egalement aux supports ferm\'es. Par contre la d\'eg\'en\'erescence de la suite spectrale \'evoqu\'ee dans ce cas n'est plus vraie pour les supports compacts. Par ailleurs ceci motivera notre d\'efinition de l'espace vectoriel monodromique \ $H^n_{c\,\cap\, S}(0)$ \ et rendra "\'evidente" les applications "naturelles" \ $H^n_c(0) \rightarrow H^n_{c\,\cap\, S}(0)$ \ et \ $ H^n_{c\,\cap\, S}(0) \rightarrow H^n(0)$.

\bigskip

\subsection{\bf Proposition.}

Pour chaque \ $k \geq 1$, on a une application surjective

$$ r^n_c(k) :  \mathbb{H}_c^n(Y, (\mathcal{E}(k)^{\bullet}, \delta^{\bullet})) \rightarrow  Ker\, \mathcal{N}_{n,c}^k \subset H_c^n(0) $$
  dont le noyau contient \ $ \tau_k \ \mathbb{H}_c^{n-1}(Y, (\mathcal{E}(k)^{\bullet}, \delta^{\bullet}))$. \\
Ces applications  v\'erifient \ $ r^n_c(k) \circ\, _k\mathcal{N} = \mathcal{N}_{n,c} \circ r^n_c(k) $ \ et elles sont compatibles aux applications \ $ j_{k,k+k'} $.
  
  \noindent De plus, ces applications induisent un isomorphisme compatible aux monodromies entre
  $$ \underset{k \rightarrow \infty}{lim}\ \mathbb{H}_c^n(Y, (\mathcal{E}(k)^{\bullet}, \delta^{\bullet})) $$
  et \ $H^n_c(0).$
  
  \bigskip
  
  \noindent \textit{Preuve.} En fait il est plus facile de d\'efinir l'application
  $$\tilde{r}^n_c(k) :  \mathbb{H}_c^n(Y, (\mathcal{E}(k)^{\bullet}, \delta^{\bullet})) \rightarrow  Ker\, \mathcal{N}_{n,c}^k \subset H_c^n(0) $$
  qui est reli\'ee \`a \ $r^n_c(k)$ \ par la formule \ $ r^n_c(k) = exp(-\mathcal{N}_{n,c}.Log (s_0) )\circ \tilde{r}^n_c(k)$, o\`u \ $ F : = f^{-1}(s_0) $ \ est la fibre de Milnor de \ $f$ \ \`a l'origine, et o\`u \ $s_0 \in D \cap \mathbb{R}^{+*} $ \ est un point base (et  on a fix\'e \ $ Log (s_0) \in \mathbb{R} .$)
  
  \noindent Comme \ $exp(-\mathcal{N}_{n,c}.Log (s_0) )$ \ est un automorphisme de \ $ H^n_c(0)$ \ compatible avec la monodromie \ $ T_c = exp(-2i\pi.\mathcal{N}_{n,c}) $ \ il suffit de prouver la proposition pour l'application \ $\tilde{r}^n_c(k)$\footnote{Pour une discussion plus d\'etaill\'ee de la relation entre \ $r^n(k)$ \ et \ $\tilde{r}^n(k)$ \  on pourra consulter la Remarque/Erratum qui suit le Corollaire 1 du Th\'eor\`eme 1 de [B.02].}.

  \noindent Pour \ $w \in \Gamma_{c/f}(X,\mathcal{E}^n(k)) \cap Ker\,\delta $ \ posons
  \begin{equation}
  \tilde{r}^n_c(k)[w] = [w_k\vert_F ] \in H^n_c(0) .
  \end{equation}
  En effet, la relation \ $ \delta w = 0 $ \ donne \ $ dw_k = \frac{df}{f}\wedge w_{k-1} $ \ ce qui montre que la restriction de \ $w$ \ \`a \ $F$ \ est bien d-ferm\'ee (et \`a support compact).
  
  \noindent Pour voir la surjectivit\'e, il suffit de constater que si \ $e \in H^n_c(0) $ \ v\'erifie \ $ \mathcal{N}^k_c(e)= 0 $ \ alors \ $ \varepsilon : = \sum_0^{k-1} \frac{(Log s )^j}{j!} \mathcal{N}^j_c(e) $ \ est une section uniforme du fibr\'e de Gauss-Manin \`a support propre de \ $f$ \ qui se calcule \`a l'aide du complexe \ $ (\mathcal{E}^{\bullet},\delta^{\bullet}) $ (voir [B.91] Prop.1 p.444) ce qui permet de trouver un ant\'ec\'edent \`a \ $e$ \ de la m\^eme fa{\c c}on que dans le cas \`a support ferm\'e.
  
  \noindent Pour d\'emontrer la relation  \ $ r^n_c(k) \circ\,_k\mathcal{N} = \mathcal{N}_{n,c} \circ r^n_c(k) $ \ nous devons montrer que, si \ $ \delta w = 0 $ \ on a \ $ [w_{k-1}\vert_F] = \mathcal{N}_c([w_k\vert_F]) .$
  
  \noindent Soit \ $ \tilde{X} : = X \underset{D}{\times} H $ \ o\`u \ $ H \overset{exp}{\rightarrow} D^* $ \ est le rev\^etement universel . D'apr\`es Milnor [Mi.68]  on a un isomorphisme \ $\mathcal{C}^{\infty}$ \ au dessus de \ $D^*$ \ $\Phi : \tilde{X} \rightarrow F \times H $. La n-forme  \ $\mathcal{C}^{\infty}$ 
  $$ \Omega : = \sum_0^{k-1} \frac{(-1)^j}{j!}.(\zeta - Log (s_0))^j .\tilde{w}_{k-j}  $$
  o\`u \ $ \tilde{w_h} = ((\Phi)^{-1})^*(w_h) $ \ v\'erifie \ $ d\Omega = 0 $ \ sur \ $ F \times H$. Elle v\'erifie 
  \'egalement \ $ \Omega\vert_{F\times \lbrace Log (s_0) \rbrace} = w_k\vert_F .$ On en d\'eduit que 
   $$ T_c(w_k\vert_F) = exp(-2i\pi.\mathcal{N}_c)(w_k\vert_F) =  \sum_0^{k-1} \frac{(-1)^j}{j!} (2i\pi)^j w_{k-j}\vert_F .$$
   On obtient alors  facilement (par r\'ecurrence) la relation d\'esir\'ee.
  
  \noindent La compatibilit\'e des applications \ $ r^n_c(k)$ \ avec les morphismes \ $ j_{k,k+k'} $ \ est imm\'ediate.
  
  \noindent Nos assertions sur le noyau se d\'eduisent facilement du lemme suivant.

  \subsection{\bf Lemme.}
    
 On utilise les notations introduites ci-dessus. Soit \ $ w \in \Gamma_c(Y, \mathcal{E}^n(k)) \cap Ker\,\delta $. Supposons que l'on ait \ $ r^n_c(k)([w]) = 0 $ \ et qu'il existe \ $ v \in \Gamma_c(Y, \mathcal{E}^{n-1}(k)) $ \ v\'erifiant \ $ _k\mathcal{N}(w) = \delta v $. Alors
  $$ j _{k,k+1}(w) = \begin{pmatrix}{w_k} \\ {\vdots} \\ {w_1}\\{0} \end{pmatrix} \in \delta \Gamma_c(Y, \mathcal{E}^{n-1}(k+1)) .$$
  
  \smallskip
  
  \noindent \textit{Preuve.} Commen{\c c}ons par traiter le cas \ $k =1$. La condition sur \ $\mathcal{N}_c(w)$ \ est vide, et on a seulement  \ $ dw_1 = 0 $ \ et \ $ [w_1\vert_F ] = 0 $. Mais alors \ $w_1$ \ est une section horizontale du fibr\'e de Gauss-Manin \`a support propre de  \ $f$ \ qui est nulle en \ $ s_0$. Elle est donc nulle, et on peut \'ecrire
  $$ w_1 = \frac{df}{f}\wedge \alpha + d\beta $$
  o\`u \ $ \alpha, \beta \in \Gamma_{c/f}(X, \mathcal{E}^{n-1}(1)) $. On a alors \ $ \frac{df}{f}\wedge d\alpha = 0 $, ce qui montre que \ $\alpha$ \ est une section du fibr\'e de Gauss-Manin \`a support propre de \ $f$ de degr\'e \ $n-1$. Mais ce fibr\'e (m\'eromorphe) est nul, car on a \ $ H^{n-1}_c(F) = 0 $ \ puisque \ $F$ \ est une vari\'et\'e de Stein de dimension \ $n$. On peut donc \'ecrire
  $$ \alpha = \frac{df}{f}\wedge \gamma + d\xi $$
  o\`u \ $ \gamma, \xi \in \Gamma_{c/f}(X, \mathcal{E}^{n-2}(1)) $. On a donc
  $$ w_1 = \frac{df}{f} \wedge d\xi + d\beta  = d(\beta - \frac{df}{f}\wedge\xi ).$$
  On obtient alors \ $ w_1 = \delta \tilde{\beta}$ \ o\`u \ $ \tilde{\beta} \in \Gamma_{c/f}(X, \mathcal{E}^{n-1}(1)) .$
   
  \noindent Passons au cas g\'en\'eral. La relation \ $_k\mathcal{N}(w) = \delta v $ \ donne 
  $$ w_{k-1} = dv_k -\frac{df}{f}\wedge w_{k-2} $$
   et \ $ \delta w = 0 $ \ donne \ $ dw_k = \frac{df}{f}\wedge w_{k-1} $. On aura donc \ $ dw_k = \frac{df}{f}\wedge dv_k $. La $n$-forme semi-m\'eromorphe \ $ w_k + \frac{df}{f}\wedge v_k$ \ est donc $d-$ferm\'ee \`a support  $f-$propre et induit  $0$  dans \ $ H^n_c(0) $. D'apr\`es le cas \ $k = 1 $ \ on peut \'ecrire  
    $$ w_k =  -\frac{df}{f}\wedge v_k  +  d\tilde{\beta} .$$
  On aura alors
  $$ j_{k,k+1}(w) = \begin{pmatrix}{w_k} \\{w_{k-1}}\\ {\vdots} \\ {w_1}\\{0} \end{pmatrix} = \delta \begin{pmatrix}{\tilde{\beta}}\\ {v_k }\\ {\vdots} \\ {v_2}\\ {v_1} \end{pmatrix} $$
  ce qui ach\`eve la preuve du lemme et de la proposition (1.9). $\hfill \blacksquare $
  
  \bigskip
  
  \subsection {\bf Remarque.}

   \noindent Comme on a \ $ j_{k,k+k'} \circ \tau_k = 0 $ \ pour  \ $k' \geq k$ \ il est inutile de quotienter par \ $\tau_k \big(\mathbb{H}^{n-1}_c(Y, (\mathcal{E}^{\bullet}(k), \delta^{\bullet}))\big)$ \ avant de prendre la limite inductive sur \ $k$ .

\bigskip

\subsection{}
 Nous  consid\'ererons \'egalement  le complexe de faisceaux sur \ $Y$ \ not\'e \ $(\mathcal{F}^{\bullet}(k), \delta_0^{\bullet} )$, d\'efini pour \ $n \geq 3$ \ comme suit :
$$(\mathcal{F}^{\bullet}(k), \delta_0^{\bullet} ) : = \  \mathcal{E}^0(k) \overset{\delta_0}{ \rightarrow } Ker\, \delta_0^1 \rightarrow 0 \rightarrow 0 \cdots  $$
On a une inclusion naturelle de complexes \ $ (\mathcal{F}(k)^{\bullet},\delta^{\bullet}_0) \rightarrow (\mathcal{E}^{\bullet}(k),\delta^{\bullet}_0 )$

\noindent Nous noterons par \ $(\mathcal{\check{E}}^{\bullet}(k),\delta_0^{\bullet} )$\ le complexe quotient ; ce qui donne de fa{\c c}on explicite :
$$ (\mathcal{\check{E}}^{\bullet}(k),\delta^{\bullet}_0 ) : = \  0 \rightarrow \mathcal{E}^1(k)/Ker\,\delta_0^1  \overset{\delta_0}{\rightarrow} \mathcal{E}^2(k) \overset{\delta_0}{\rightarrow} \mathcal{E}^3(k) \overset{\delta_0}{\rightarrow}\cdots $$
Le complexe \ $(\mathcal{\check{E}}^{\bullet}(k),\delta^{\bullet}_0 )$ \  n'aura  sur \ $Y^*$ \ que deux faisceaux de cohomologie donc nuls, comme c'\'etait le cas pour le complexe \ $(\mathcal{E}^{\bullet}(k),\delta^{\bullet}_u )$ \ 
pour une  valeur propre diff\'erente de  1. Ceci nous permettra d'utiliser sur \ $S^*$ \  la proposition 2 p. 435  de [B.91].

\noindent Pour \ $n = 2$ \ nous d\'efinirons le complexe \ $(\mathcal{F}^{\bullet}(k),\delta_0^{\bullet} )$ \ en rempla{\c c}ant dans la d\'efinition pr\'ec\'edente le faisceau \ $Ker\, \delta_0^1$ \ par le sous-faisceau 
$$\delta_0 \mathcal{E}^{0}(k) + \tau_k h^{0}(k) \hookrightarrow Ker\,\delta_0^{1}(k) .$$
De m\^eme nous remplacerons dans la d\'efinition de \ $(\mathcal{\check{E}}^{\bullet}(k),\delta^{\bullet}_0 )$ \ le quotient \ $\mathcal{E}^1(k)/Ker\,\delta_0^1  $ \ par 
$$\mathcal{E}^1(k)\big/ \delta_0 \mathcal{E}^{0}(k) + \tau_k h^{0}(k) .$$
Ceci conduit \`a nouveau au fait que le complexe \ $(\mathcal{\check{E}}^{\bullet}(k),\delta^{\bullet}_0 )$ \ ait deux faisceaux de cohomologie non nuls sur \ $Y^*$  isomorphes  \`a \ $h^{1}(k)/\tau_k (h^{0}(k) $ \ et \ $h^2(k)$ \ en degr\'es \ $ n-1 = 1 $ \ et \ $ n = 2$ \ respectivement. 

\noindent Le quotient \ $ \mathcal{E}^1(k)\big/ \delta_0 \mathcal{E}^{0}(k) + \tau_k h^{0}(k) $ \ n' a pas  la cohomologie sur \ $S^*$\ en degr\'es positifs. En effet, on a
$$ \delta_0\mathcal{E}^0(k)\,\cap\, \tau_k(h^0(k)) \simeq 0 $$
 On en d\'eduit facilement que \ $ \forall i \geq 1$
 $$H^i(S^*, \mathcal{E}^1(k)\big/ \delta_0\mathcal{E}^{0}(k) + \tau_k h^{0}(k) ) \simeq H^{i+1}(S^*, h^0(k)) \oplus H^{i+1}(S^*, \delta_0\mathcal{E}^{0}(k)) \simeq 0 .$$
 On  pourra donc encore utiliser sur \ $S^*$ \ la proposition 2 p. 435 de [B.91] pour d\'ecrire  l'hypercohomologie sur \ $S^*$ \ du complexe \ $(\mathcal{\check{E}}^{\bullet}(k),\delta^{\bullet}_0 )$.

\bigskip

\section{\bf  Le cas d'une singularit\'e isol\'ee pour la valeur propre 1 .}

\noindent Consid\'erons ici le cas o\`u  on a \ $n \geq 1 $ \ et o\`u  l'on suppose que les faisceaux \ $H^p(0)$\ sont tous nuls pour \ $p  \not= n $. Quitte \`a choisir le repr\'esentant de Milnor de notre germe assez petit on peut \'egalement supposer que le faisceau \ $H^n(0) $ \ est concentr\'e \`a l'origine. On a alors les r\'esultats suivants.

\bigskip

\subsection{\bf {Th\'eor\`eme.}\rm (voir [B.90]).}

\noindent {\it Soit \ $\tilde{f} : (\mathbb{C}^{n+1}, 0) \rightarrow (\mathbb{C}, 0)$ \ un germe de fonction holomorphe pr\'esentant une singularit\'e isol\'ee \`a l'origine pour la valeur propre  1 de la monodromie\footnote{Donc \ $H^p(0) = 0 $ \ pour \ $ p \not= n$ \ puisque, par perversit\'e, seul le faisceau \ $H^n(0)$ \ peut avoir un support de dimension \ $0$.}.
Il existe une forme hermitienne  \ $h$ \ sur \ $H^n(0)$, non d\'eg\'en\'er\'ee, invariante par la monodromie, d\'efinie de la fa{\c c}on suivante :}

\noindent {\it Soient  \ $ e, e' \in H^n(0) $ \ v\'erifiant \ $\mathcal{N}^k(e) =\mathcal{N}^k(e') = 0 $ \ o\`u on a not\'e par \ $-2i\pi.\mathcal{N}$ \ le logarithme nilpotent de la monodromie (unipotente)  agissant sur \ $H^n(0)$. Soient \ $ w,w' \in H^n(\Gamma(X, \mathcal{E}^{\bullet}(k)), \delta_0^{\bullet}) $ \ v\'erifiant \ $ r^n(k)(w) = e$ \ et \ $ r^n(k)(w') = e'$. Alors on a }
$$ (2i\pi)^{n+1}. h(e,e') = P_2\big(\lambda = 0, \int_X \vert f \vert^{2\lambda} \rho.\frac{df}{f}\wedge w_k \wedge \frac{\bar{df}}{\bar{f}} \wedge \bar{w'}_k \big) $$
{\it pour \ $\rho \in \mathcal{C}^{\infty}_c(X) $ \ v\'erifiant \ $ \rho \equiv 1 $ \ pr\`es de l'origine.}

\bigskip

\noindent On a not\'e \ $ P_2(\lambda = 0, F(\lambda)) $ \ le coefficient de \ $\frac{1}{\lambda^2}$ \ dans le d\'eveloppement de Laurent en \ $ \lambda = 0 $ \ de la fonction m\'eromorphe \ $F.$

\noindent On utilisera dans cet article, comme on l'a d\'ej\`a fait dans [B.97],  les complexes \ $(\Omega^{\bullet}(k), \delta^{\bullet}_0)$ \ pour le cas m\'eromorphe et \ $(\mathcal{E}^{\bullet}(k), \delta^{\bullet}_0 )$ \ pour le cas semi-m\'eromorphes au lieu de ceux correspondants \`a la differentielle \ $ \delta^{\bullet}_1$ \ de la diff\'erentielle (utilis\'es dans [B.90]) ce qui explique que l'\'enonc\'e ci-dessus diff\`ere de celui de [B.90] (en particulier par " $ \lambda = 0$ "). Nous renvoyons \`a [B.91], aux rappels donn\'es dans [B.97] paragraphe 2 ou au paragraphe 1 ci-dessus pour plus de pr\'ecisions sur ces complexes et sur l'application \ $r^n(k)$.

\noindent La forme hermitienne d\'efinie ci-dessus dans le cas d'une singularit\'e  isol\'ee pour la valeur propre 1 de la monodromie sera appel\'ee la forme hermitienne canonique.

\bigskip

\subsection{}
Dans le cas d'une singularit\'e isol\'ee, la forme hermitienne canonique a \'et\'e introduite dans [B.85]. Fran{\c c}ois  Loeser [Loe.86]  a montr\'e dans ce cas, en utilisant des r\'esultats de th\'eorie de Hodge, comment construire la forme hermitienne canonique \`a partir de la variation et de la dualit\'e (hermitienne) de Poincar\'e sur la fibre de Milnor. Ceci a \'et\'e g\'en\'eralis\'e au cas d'une singularit\'e isol\'ee pour la valeur propre 1  (sans th\'eorie de Hodge, car elle ne semble pas facilement applicable dans ce contexte plus g\'en\'eral) dans [B.97]. Rappelons les \'enonc\'es correspondants.

\bigskip

\subsection{\bf {Th\'eor\`eme.} \rm{(voir [B.97]).}}

\noindent \textit{ Soit \ $\tilde{f} : (\mathbb{C}^{n+1}, 0) \rightarrow (\mathbb{C}, 0)$ \ un germe de fonction holomorphe pr\'esentant une singularit\'e isol\'ee \`a l'origine pour la valeur propre  1 de la monodromie.
Il existe une application lin\'eaire "naturelle" de variation \ $ var : H^n(0) \rightarrow H^n_c(0) $ \ v\'erifiant les propri\'et\'es suivantes : }
\begin{itemize}
\item{0)} \it{Si la singularit\'e de \ $f$ \ est isol\'ee\footnote{Au sens "usuel", c'est \`a dire que l'on a \ $ df_x \not= 0 $ \ pour \ $ x \not= 0 .$}, on retrouve la variation "classique" (voir [A.V.G.]).}
\item{1)}\it{ L'application "var"  est topologique et on a }
\begin{align*}
can \circ var = T - 1 \quad {\rm sur} \quad  H^n(0)  \\
var \circ can = T_c - 1 \quad {\rm sur} \quad    H^n_c(0) 
\end{align*}
\textit{o\`u \ $ can : H^n_c(0) \rightarrow H^n(0) $ \ est l'application d'oubli de support et o\`u \ $T$ \ (resp. \ $T_c$) d\'esigne la monodromie agissant sur \ $H^n(0)$ \ (resp. sur \ $H^n_c(0)$ ).
\item{2)} L'application "var"  est bijective, $T-$invariante (c'est \`a dire que \ $ T_c \circ var = var \circ T$), et auto-adjointe pour la dualit\'e (hermitienne) de Poincar\'e \ $ \mathcal{I} : H^n(0) \times H^n_c(0) \rightarrow \mathbb{C} $ \ donn\'ee par}
$$ \mathcal{I}(a,b) = \frac{1}{(2i\pi)^n} \int_F a \wedge \bar{b} . $$
\textit{o\`u \ $F$ \ d\'esigne la fibre de Milnor de} \ $f$ \ en \ $0$.
\item{3)} Soit \ $ \Theta : = \sum_{k=0}^{\infty}  (-1)^k \frac{(T-1)^k}{k+1} .$\it{ C'est un automorphisme de \ $H^n(0)$ \ puisque \ $T$ \ est unipotent. Posons} \ $ \tilde{var} : = var \circ \Theta $. \it{Alors on a }
$$ \forall e, e' \in H^n(0), \quad h(e,e') = \mathcal{I}(\tilde{var}(e), e') .$$
\end{itemize}

\noindent On d\'eduit facilement des propri\'et\'es ci-dessus que l'on a
\begin{align*}
 \quad & Im\ can = Im\ (T-1) \qquad  &{\rm dans} \quad H^n(0)  \\
{\rm et} \qquad & Ker\ can = Ker\ (T_c -1) \qquad  &{\rm dans} \quad H^n_c(0) .
\end{align*}

\bigskip

\subsection{} 
Comme nous aurons besoin plus loin de construire une g\'en\'eralisation de l'application de variation, nous allons rappeler ici bri\`evement l'expression  de l'application \ $ \tilde{var}$ \ donn\'ee dans [B.97] en terme de formes diff\'erentielles semi-m\'eromorphes.

\noindent Soit donc \ $e \in H^n(0)$ \ v\'erifiant \ $\mathcal{N}^k(e) = 0 $. Soit \ $ w \in \Gamma(Y, \mathcal{E}^n(k)) $ \ satisfaisant \ $\delta_0(w) = 0 $ \ et \ $ r^n(k) = e .$ L'hypoth\`ese de singularit\'e isol\'ee pour la valeur propre 1  permet alors de trouver \ $ u \in \Gamma(Y^*, \mathcal{E}^{n-1}(k)) $ \ v\'erifiant
$$ _kN\, w \vert_{Y^*} = \delta_0\, u  .$$
Soit \ $\rho \in \mathcal{C}^{\infty}_{c/f}(X) $ \ v\'erifiant \ $ \rho \equiv 1 $ \ pr\`es de l'origine. Posons \ $ \xi = (1-\rho).u $ . Alors on a \ $ \xi \in \Gamma_c(Y, \mathcal{E}^{n-1}(k))$ \ et \ $ v =\   _kN\, w - \delta_0\,\xi $ \ est dans \ $\Gamma_c(Y,  Ker\, \delta_0^n(k)) .$ On a donc, pour \ $j \in [1,k]$,  avec la convention \ $w_0 = 0 = \xi_0 $ :
$$  w_{j-1} = v_j + d\xi_j - \frac{df}{f}\wedge \xi_{j-1} \ ;$$
et comme \ $ \delta_0\, w = 0 $ \ implique \ $ dw_k = \frac{df}{f} \wedge w_{k-1} $ \ on aura
$$  d(w_k + \frac{df}{f} \wedge \xi_k ) = \frac{df}{f} \wedge v_k .$$
Choisissons maintenant une fonction \ $ \sigma \in \mathcal{C}^{\infty}_{c/f}(X) $ \ v\'erifiant \ $ \sigma \equiv 1 $ \ sur un ouvert  de \ $Y$ \ contenant \ $ (Supp\, v)\, \cap\, Y $. Alors le th\'eor\`eme de Leray g\'en\'eralis\'e au cas d'une hypersurface non n\'ecessairement lisse mais ne pr\'esentant pas la valeur propre 1 pour la monodromie\footnote{Ce point est d\'etaill\'e dans [B.97] ; il sera repris et g\'en\'eralis\'e plus loin; voir le 5.}, permet d'\'ecrire au voisinage du bord de la boule de Milnor
$$ d\sigma \wedge (w_k + \frac{df}{f} \wedge \xi_k ) = \frac{df}{f}\wedge \eta + \omega + d\alpha $$
avec \ $\eta$ \ et \ $\omega$ \ d-ferm\'ees et \ $ \mathcal{C}^{\infty}$ \ \`a support \ $f-$propres de degr\'es respectifs \ $n$ \ et \ $n+1$, et o\`u \ $\alpha$ \ est semi-m\'eromorphe de degr\'e \ $n$ \ \`a p\^oles dans \ $Y$ \ et \`a support \ $f-$propre. On obtient alors que la section \ $V(w) \in \Gamma_c(Y, \mathcal{E}^n(k))$ \ d\'efinie par
\begin{align*}
 &  V(w)_k  = v_k + \eta  \\
 &   V(w)_j =  v_j  \quad \quad  {\rm pour} \quad j \in [1,k-1]
\end{align*}
est \ $\delta_0-$ferm\'ee et on a \ $ r^n(k)(V(w)) = \tilde{var}(e)$. Pour plus de d\'etails, on consultera [B.97], ou bien la construction de la variation g\'en\'eralis\'ee donn\'ee au paragraphe 6.

\bigskip

\section{}
 L'objectif de ce troisi\`eme paragraphe est de g\'en\'eraliser \`a notre situation la proposition 11 de [B.91]. 

\bigskip

\subsection {\bf La forme hermitienne canonique sur le syst\`eme local \ $H^{n-1}(0)\vert_{S^*}.$}

\subsubsection{}
 Consid\'erons donc le syst\`eme local \ $H^{n-1}(0)$ \ sur \ $S^*$. Sa fibre en chaque point \ $\sigma \in S^*$ \ s'identifie au sous-espace spectral pour la valeur propre 1 du  \ $(n-1)-$i\`eme groupe de cohomologie de la fibre de Milnor de la restriction de \ $f$ \`a une section hyperplane transverse \`a \ $S^*$ \ en \ $\sigma$. Comme la restriction de \ $f$ \ \`a un tel hyperplan a une singularit\'e isol\'ee pour la valeur propre 1  en \ $\sigma$, ce sous-espace spectral pour la valeur propre 1 de ce \ $(n-1)-$i\`eme groupe de cohomologie est muni d'une forme hermitienne canonique qui est non d\'eg\'en\'er\'ee et invariante par la monodromie de \ $f$. Ceci permet, en raisonnant comme dans [B.91] p.456, de munir le syst\`eme local (d'espaces vectoriels monodromiques) \ $H^{n-1}(0)\vert_{S^*}$ \ d'une forme hermitienne localement constante, non d\'eg\'en\'er\'ee et invariante par la monodromie\footnote{On prendra garde qu'\`a l'origine, cette construction n'a, \`a priori, pas de sens.}. Nous la noterons par
$$ h :  H^{n-1}(0)\vert_{S^*} \times H^{n-1}(0)\vert_{S^*} \rightarrow \underline{\mathbb{C}}_{S^*} .$$
On montre comme au lemme 1 de [B.91] p.456 que cette forme hermitienne est intrins\`eque.
La proposition suivante se d\'emontre de fa{\c c}on analogue \`a la proposition 1 de [B.91] p.457.

\bigskip

\subsubsection{\bf Proposition.}

 On se place sous les hypoth\`eses standards pour la valeur propre 1. Soit \ $V$ \ un ouvert de \ $X^* : = X \setminus \lbrace 0 \rbrace$ et soient \ $w$ \ et \ $w'$ \ deux sections sur \ $V$ \ du faisceau \ $\mathcal{E}^{n-1}(k)$ \ v\'erifiant \ $ \delta(w) = \delta(w') \equiv 0 $. Notons par \ $e$ \ et \ $e'$ \ les sections correspondantes de \ $H^{n-1}(0)$ \ sur \ $U : = V\, \cap \, S^* $ (via le morphisme \ $r^{n-1}(k)$). Alors la fonction \ $h(e,e')$ \ est localement constante sur \ $U$. Pour chaque forme diff\'erentielle \ $ \varphi \in \mathcal{C}^{\infty}(V)$ \ de degr\'e  2  v\'erifiant :
\begin{itemize}
\item{i)} $ Supp\,\varphi \,\cap\, S^* $ \ est compact,
\item{ii)}  \ $d\varphi = 0 $ \ au voisinage de \ $U$ \ dans \ $X^*$,
\end{itemize}
on aura la formule :
$$  (2i\pi)^{n} \int_U h(e,e').\varphi  =  P_2 \big( \lambda = 0, \int_X \vert f \vert^{2\lambda} w_k \wedge \bar{w'}_k \wedge \varphi \wedge \frac{df}{f} \wedge \frac{\bar{df}}{\bar{f}} \big). $$

\bigskip

\subsubsection{} 
Quitte \`a remplacer la formule (1) p.409 de [B.91] donnant l'expressions de la forme (hermitienne) d'intersection comme r\'esidu dans le prolongement analytique de la distribution  \ $ \vert f \vert^{2\lambda}$ \ par la formule du th\'eor\`eme (0.6), la preuve de cette proposition est analogue \`a celle de la proposition 1 de [B.91] p.457-460 .

\bigskip

\subsubsection{Remarque.}
Pour \ $n = 2$ \ il n'est pas \`a priori automatique de pouvoir repr\'esenter une section sur \ $U : = V\, \cap \, S^* $ du faisceau \ $H^1(0)$ \ par un \'el\'ement de \ $\Gamma(V, Ker\,\delta^1(k))$ \ pour \ $ k \gg 1$, car cela n\'ecessite la surjectivit\'e des applications 
$$ H^0(V,  Ker\,\delta^1(k)) \rightarrow H^0(V, h^1(k)) \rightarrow H^0(V, H^1(0)) .$$
L'isomorphisme \ $ h^1(k)\big/\tau_k(h^0(k)) \rightarrow H^1(0) $ \ pour \ $ k \gg 1$, montre que l'on a une obstruction \`a la surjectivit\'e de la seconde fl\`eche. En fait le d\'ecalage sur \ $k$ \ suffit \`a lever cette obstruction puisque \ $ \underset{k \rightarrow \infty}{lim} h^1(k) \rightarrow H^1(0) $ \ est un isomorphisme.

\noindent La surjectivit\'e de la premi\`ere fl\`eche est realis\'ee d\`es que l'on a \ $ H^1(V, \delta\mathcal{E}^0(k)) = 0 $, ce qui sera toujours vrai quitte \`a prendre pour \ $V$ \ un voisinage ouvert assez petit de \ $U$, car on a  \ $H^1(V, \delta\mathcal{E}^0(k)) \simeq H^2(V, h^0(k)) \simeq H^2(U, \mathbb{C}) = 0 .$ 

\subsubsection{}
 En utilisant le cycle fondamental de \ $S^*$ \ qui est l'\'el\'ement de \ $H_1(S^*, \mathbb{Z})$ \ obtenu en coupant \ $S^*$ \ par une sph\`ere  centr\'ee \`a l'origine de rayon assez petit, on obtient,  un accouplement sesquilin\'eaire non d\'eg\'en\'er\'e
$$ \tilde{h} : H^0(S^*, H^{n-1}(0)) \times H^1(S^*, H^{n-1}(0) \rightarrow \mathbb{C} .$$
Pour donner une formule analytique donnant le nombre \ $\tilde{h}(e,e')$ \ analogue \`a la proposition 11 de [B.91], nous devons commencer par d\'efinir l'analogue des applications  \ $\tilde{\theta}_k$ \ et \ $\tilde{ob}_k $ \ de loc. cit.

\bigskip

\subsection{\bf Les applications \ $\tilde{\theta}_k$ \ et \ $\tilde{ob}_k $.}

\bigskip

\subsubsection{}
 D\'efinissons l'entier \ $k_0$ \ comme l'ordre de nilpotence de la monodromie agissant sur le syst\`eme local \ $H^{n-1}(0)\vert_{S^*}$\footnote{Comme l'application de restriction \ $ \Gamma(S, H^{n-1}(0)) \rightarrow \Gamma(S^*, H^{n-1}(0))$ \ est injective (perversit\'e), \ $k_0$ \ majore \'egalement l'ordre de nilpotence de la monodromie agissant sur la fibre \`a l'origine de \ $H^{n-1}(0)$.}.

\bigskip

\subsubsection{} {\bf D\'efinition.}

\noindent {\it Sous les hypoth\`eses standards pour la valeur propre 1 et  pour \ $n \geq 3$, notons par \ $-2i\pi\mathcal{N}_n$\ le logarithme nilpotent de la monodromie agissant sur \ $H^n(0)$. Pour chaque \ $k \geq k_0 $ \ on a deux applications lin\'eaires naturelles }
\begin{align*}
    \tilde{ob}_k  : & \ Ker\, \mathcal{N}_n^k \rightarrow  H^1_{\lbrace 0 \rbrace}(S, H^{n-1}(0)) \\
    \tilde{\theta}_k  :  & \ Ker\, \mathcal{N}_n^k \rightarrow  H^1(S^*, H^{n-1}(0)) \simeq  H^2_{\lbrace 0 \rbrace}(S, H^{n-1}(0))
\end{align*}
{\it qui sont d\'efinies de la fa{\c c}on suivante :}
\begin{itemize}
\item{} {\it La compos\'ee}
$$ \mathbb{H}^n(Y, (\mathcal{E}^{\bullet}(k), \delta^{\bullet})) \rightarrow  \mathbb{H}^n(S^*, (\mathcal{\check{E}}^{\bullet}(k), \delta^{\bullet})) \overset{ob_k}{\rightarrow} H^0(S^*, H^{n-1}(0)) \overset{\partial}{\rightarrow} H^1_{\lbrace 0 \rbrace}(S, H^{n-1}(0))  $$
{\it passe au quotient par} \ $\tau_k \ \mathbb{H}^{n-1}(Y, (\mathcal{E}^{\bullet}(k), \delta^{\bullet}))$ \ {\it et donne} \ $\tilde{ob_k}$\footnote{Grace aux isomorphismes rappel\'es plus haut ; voir (1.8).}.
\item{} {\it La compos\'ee}
$$ \mathbb{H}^n(Y, (\mathcal{E}^{\bullet}(k), \delta^{\bullet})) \rightarrow  \mathbb{H}^n(S^*, (\mathcal{\check{E}}^{\bullet}(k), \delta^{\bullet})) \overset{\theta_k}{\rightarrow} H^1(S^*, H^{n-1}(0))  $$
{\it passe \'egalement  au quotient par} \ $\tau_k \ \mathbb{H}^{n-1}(Y, (\mathcal{E}^{\bullet}(k), \delta^{\bullet}))$ \ {\it et donne }\ $\tilde{\theta_k}.$
\end{itemize}

\bigskip

\noindent Ces r\'esultats sont analogues \`a ceux des pages 410 et 411 de [B.91], le complexe \ $(\mathcal{\check{E}}^{\bullet}(k), \delta^{\bullet}))$ \ ayant les m\^emes propri\'et\'es dans notre cas que le complexe "usuel" pour une valeur propre diff\'erente de  1 (puisqu'il n'a que deux faisceaux de cohomologie non nuls sur \ $Y^*$, et on a l'acyclicit\'e sur \ $S^*$ \ des faisceaux \ $\mathcal{\check{E}}^{\bullet}(k)$.) 

\noindent On prendra garde qu'ici \`a nouveau les applications \ $\tilde{\theta}_k$ \ sont compatibles quand \ $k$ \ augmente, alors que l'on a \ $\tilde{ob}_{k+1} = \tilde{ob}_k\circ \mathcal{N}_n.$

\noindent Nous d\'efinirons alors l'application \ $ \theta : H^n(0) \rightarrow  H^2_{\lbrace 0 \rbrace}(S, H^{n-1}(0)) $ \ comme \'etant \ $\tilde{\theta_k}$ \ pour \ $k \gg 1$.

\bigskip

\subsubsection{\bf Remarque.}

Toujours en supposant \ $n \geq 3$, le faisceau \ $ \check{\mathcal{E}}^1(k) \simeq \mathcal{E}^1(k) \big/ Ker\,\delta^1 $ \ n'a pas de cohomologie en degr\'e positif sur \ $S^*$. En effet on a les suites exactes
\begin{align*}
0 \rightarrow Ker\,\delta^1 \rightarrow \mathcal{E}^1(k) \rightarrow \check{\mathcal{E}}^1(k) \rightarrow 0 \\
0 \rightarrow Im\,\delta^0 \rightarrow  Ker\,\delta^1 \rightarrow h^1(k) \rightarrow 0 \\
0 \rightarrow h^0(k) \rightarrow \mathcal{E}^0(k) \rightarrow  Im\,\delta^0 \rightarrow 0
\end{align*}
qui donnent successivement les annulations pour tout \ $k \geq 1$
\begin{align*}
& H^i(S^*, Im\,\delta^0(k)) = 0 , & \forall i \geq 1  \\
& H^j(S^*, Ker\,\delta^1(k)) = 0 , & \forall j \geq 2 \\
& H^l(S^*, \check{\mathcal{E}}^1(k) ) = 0 , & \forall l \geq 1
 \end{align*}
 On a donc la d\'eg\'en\'erescence de la premi\`ere suite spectrale donnant l'hypercohomologie sur \ $S^*$ \ du complexe \ $ (\check{\mathcal{E}}^{\bullet}(k) ,\delta^{\bullet}) $ \ et donc un isomorphisme
 $$ \mathbb{H}^n(S^*, (\mathcal{\check{E}}(k)^{\bullet}, \delta^{\bullet})) \simeq H^n(\Gamma(S^*, \mathcal{\check{E}}(k)^{\bullet}), \delta^{\bullet}) .$$
 Nous pouvons \'etendre la d\'efinition des applications  \ $\tilde{\theta}_k$ \ et \ $\tilde{ob}_k $ \ aux \'el\'ements \ $\delta-$ferm\'es de \ $ \Gamma(S^*, \mathcal{E}^n(k))$ \ pour \ $k \geq k_0$ : ils induisent des classes dans \ $\mathbb{H}^n(S^*, (\mathcal{\check{E}}(k)^{\bullet}, \delta^{\bullet}))$ \ qui s'envoie dans  \ $\mathbb{H}^n(Y^*, (\mathcal{\check{E}}(k)^{\bullet}, \delta^{\bullet})) \simeq \mathbb{H}^n(S^*, (\mathcal{\check{E}}(k)^{\bullet}, \delta^{\bullet})) $ \ puisque les faisceaux de cohomologie non nuls du complexe \ $(\mathcal{\check{E}}(k)^{\bullet}, \delta^{\bullet})$ \ sont \`a supports dans \ $S$. Pour \'eviter toute confusion, nous noterons ces applications  par 
$$ \widehat{ob}_k : \Gamma(S^*, \mathcal{E}^n(k)) \cap Ker\,\delta \rightarrow H^0(S^*, H^{n-1}(0))$$
 et par  
$$ \widehat{\theta}_k : \Gamma(S^*, \mathcal{E}^n(k)) \cap Ker\,\delta \rightarrow H^1(S^*, H^{n-1}(0)).$$
On notera que \ $\widehat{ob}_k $ \ est simplement donn\'ee par la compos\'ee des morphismes de faisceaux sur \ $S^*$ 
$$ (\tau_k)^{-1} : h^n(k) \rightarrow h^{n-1}(k) \hookrightarrow H^{n-1}(0) .$$

\subsubsection{\bf Le cas \ $n = 2$.}

La d\'efinition des applications \ $\tilde{ob}_k $ \ et \ $\theta$ \ ne pose pas de probl\`eme s\'erieux pour \ $n = 2$. Une mani\`ere de le voir directement (sans utiliser [B.91]) est de transformer l'\'ecriture locale
$$ w\vert_{X_{\sigma}} = \tau_k(\alpha_{\sigma}) + \delta\beta_{\sigma} $$
correspondant \`a la surjection sur \ $S^*$ \ $\tau_k ; h^1(k) \rightarrow h^2(k) $, en
$$ j_{k,2k}(w)\vert_{X_{\sigma}} =  \delta \begin{pmatrix}\beta_{\sigma}\\ \alpha_{\sigma}\end{pmatrix} .$$
Alors 
 $$\begin{pmatrix} \beta_{\sigma} - \beta_{\sigma'} \\ \alpha_{\sigma} - \alpha_{\sigma'} \end{pmatrix}$$ induira un \'el\'ement de \ $H^1(S^*, h^1(2k))$ \ dont l'image dans \ $H^1(S^*, H^1(0))$ \ (pour \ $k \gg1$ ) sera \ $\theta([w])$. La cocha\^ine \ $\alpha_{\sigma}$ \ donnant une section globale sur \ $S^*$ du faisceau quotient \ $h^1(k)\big/\tau_k(h^0(k) \simeq H^1(0)$ \ pour \ $ k \geq k_0 $, ce qui d\'efinit \ $\tilde{ob}_k([w]).$
 
 \bigskip

\noindent Donnons maintenant l'analogue de la proposition 11 de [B.91].

 \subsubsection {\bf Th\'eor\`eme.}

\noindent {\it Sous les hypoth\`eses standards pour la valeur propre 1, consid\'erons des entiers }\ $k \geq k_0$ \ et \ $k' \geq k_0$. {\it Soient }\ $ v\in \Gamma(S^*, \mathcal{E}^n(k)) \cap Ker\, \delta $ \ et \ $ w\in \Gamma(S^*, \mathcal{E}^n(k)) \cap Ker\, \delta $. {\it Notons par} \ $ e : = \widehat{ob}_k(v)\in H^0(S^*, H^{n-1}(0))$ \ et par \ $ \eta : = \widehat{\theta}_k (w)\in H^1(S^*, H^{n-1}(0))$. {\it Alors on a }
$$ (2i\pi)^{n+1}.\tilde{h}(e,\eta) = \sum_{a=1}^{k_0 } (-1)^a P_{a+1} \big(\lambda = 0, \int_X \vert f \vert^{2\lambda} \frac{df}{f}\wedge \bar{w'}_{k'}\wedge v_a \wedge \gamma \big) $$
{\it o\`u \ $\gamma$ \ est une forme \ $\mathcal{C}^{\infty}$ \ de degr\'e  1  au voisinage de \ $S^*$ \ v\'erifiant}\footnote{Comme dans [B.91], si \ $\rho \in \mathcal{C}^{\infty}_c(X)$ \ vaut identiquement   1  pr\`es de l'origine, \ $ \gamma = d\rho$ \ convient.}
\begin{itemize}
\item{i)}  $Supp\, \gamma \cap \, S^* $ \ {\it est compact ;}
\item{ii)} \ $ d\gamma = 0$ \ {\it au voisinage de} \ $S^*$ ;
\item{iii)}  \ $\gamma$ \ {\it induit la classe du cycle fondamental de \ $S^*$ \ dans} \ $H^1_c(S^*, \mathbb{C}).$
\end{itemize}

\bigskip

\subsubsection{\bf Corollaire 1.}

\noindent Sous les hypoth\`eses du th\'eor\`eme pr\'ec\'edent, supposons que \ $v$ \ et \ $w$ \  soient  restrictions \`a \ $S^*$ \ d' \'el\'ements de \ $ \Gamma(S, \mathcal{E}^n(k)) \cap Ker\, \delta$. Alors on a
$ \tilde{h}(e,\eta) = 0 .$

\smallskip

\noindent  Ce corollaire est analogue au corollaire 1 de la proposition 11  de [B.91].

\bigskip

\subsubsection{\bf Corollaire 2.}

 Soient \ $k,k' \geq k_0 $ \ et soient \ $w \in H^n((\Gamma(Y, \mathcal{E}^{\bullet}(k)), \delta))$ \ et \ $ w' \in H^n((\Gamma(Y, \mathcal{E}^{\bullet}(k')), \delta))$. Posons \ $ e : = r^n(k)(w)$ \ et \ $e' : = r^n(k')(w')$. Alors on a
$$ (2i\pi)^{n+1}.\tilde{h}(\tilde{ob}_k(e), \theta (e')) = (-1)^k P_{k+2}\big(\lambda = 0, \int_X \vert f \vert^{2\lambda} \frac{df}{f}\wedge w_k \wedge \frac{\bar{df}}{\bar{f}} \wedge \bar{w'}_{k'} \big).$$

\smallskip

\noindent Ce corollaire est analogue au corollaire 2 de la proposition 11  de [B.91].

\bigskip

\subsubsection {\bf Corollaire 3.}

Sous les hypoth\`eses standards, l'image de \ $\theta$ \ est contenue dans l'orthogonal pour \ $\tilde{h}$ \ du sous-espace \ $ H^0(S, H^{n-1}(0))$ \ de \ $H^0(S^*, H^{n-1}(0)).$

\smallskip

\noindent Ce dernier corollaire est analogue \`a la premi\`ere partie de la proposition 12 de [B.91]. Il se d\'eduit imm\'ediatemment du corollaire 1 (voir la preuve du corollaire 1 de la proposition 11 de [B.91] p.466).

\bigskip

\subsubsection{\bf Remarque.}

 Soient \ $k \geq k_0$ , et \ $e \in Ker \mathcal{N}_n^k $ \ tels que \ $\tilde{ob}_k(e) \not= 0 $. Pour conclure grace au corollaire 2  que l'on aura alors un p\^ole d'ordre \ $k+2$ \  aux entiers n\'egatifs dans le prolongement m\'eromorphe de \ $ \vert f \vert ^{2\lambda} $ \  (voir le th\'eor\`eme 13  de [B.91] et notre th\'eor\`eme 3), il est essentiel d'avoir l'\'egalit\'e entre l'image de \ $\theta $ \ dans \ $ H^1(S^*, H^{n-1}(0)) $ \ et l'orthogonal pour la forme hermitienne \ $\tilde{h}$ \ de \ $H^0(S, H^{n-1}(0))$. C'est \`a dire de savoir que l'inclusion donn\'ee au corollaire 3  est en fait une \'egalit\'e. Ce qui revient \`a montrer que  \ $\tilde{h}$ \ \'etablit une dualit\'e hermitienne entre \ $H^1_{\lbrace 0 \rbrace}(S, H^{n-1}(0)) $ \ et \ $ Im \ \theta $. Ceci sera notre objectif dans ce qui suit. Il demandera encore beaucoup  de travail puisqu'il ne sera atteint  qu' \`a la fin du paragraphe 8. 

\smallskip

\noindent Ce r\'esultat sera \'egalement la clef de la non-d\'eg\'en\'erescence de la forme hermitienne canonique "g\'en\'eralis\'ee" que nous allons introduire plus loin.

\bigskip

\section{\bf L'espace vectoriel monodromique  \ $H^n_{c\,\cap\,S}(0)$.}

\bigskip

\subsection{}

 Introduisons sur \ $X$ \ et \ $Y$ \ les familles paracompactifiantes de supports suivants 
\begin{itemize}
\item{} \ $ c/f$ \  la famille des ferm\'es \ $f-$propres de \ $X$. Elle donne par restriction \`a \ $Y$ \ la famille des compacts de \ $Y$ \ que nous noterons simplement  "c".
\item{}  \ $ c \cap S$ \ la famille des ferm\'es de  \ $X$ \ (ou de \ $Y$ ) qui rencontrent \ $S$ \ suivant un compact.
\item{} \ $mod\,S$\ la famille des ferm\'es de \ $X$ \ (ou de \ $Y$) qui ne rencontrent pas \ $S$. 
\end{itemize}

\noindent Nous utiliserons \'egalement ces familles de supports pour un ouvert de \ $X$ \ ou de \ $Y$.

\bigskip

\subsection{}

Comme les faisceaux \ $  \mathcal{E}^{\bullet}(k)$ \ sont fins et la famille \ $c\,\cap\,S $ \ paracompactifiante, le \ $j-$i\`eme groupe de cohomologie du complexe 
$$\big( \Gamma_{c\,\cap\,S}(Y, \mathcal{E}^{\bullet}(k)), \delta^{\bullet}\big) $$
est isomorphe \`a l'hypercohomologie $$\mathbb{H}_{c\,\cap\,S}^j(Y, \mathcal{E}^{\bullet}(k)), \delta^{\bullet}).$$
Comme \ $\tau_k $ \ est un morphisme de complexe de degr\'e +1 il induit une fl\`eche
$$ \tau_k : \mathbb{H}_{c\,\cap\,S}^{j-1}(Y, \mathcal{E}^{\bullet}(k)), \delta^{\bullet}) \rightarrow \mathbb{H}_{c\,\cap\,S}^j(Y, \mathcal{E}^{\bullet}(k)), \delta^{\bullet}) .$$
Mais  l'\'egalit\'e\footnote{En fait le lemme 2 p. 442 de [B.91] donne une homotopie \ $t_k : \mathcal{E}^{\bullet}(k) \rightarrow \mathcal{E}^{\bullet}(k)$ \ v\'erifiant \ $\tau_k \circ _k\mathcal{N} - _k\mathcal{N}\circ \tau_k = \delta \circ t_k - t_{k}\circ \delta $ \ qui, combin\'ee avec l'\'egalit\'e \'evidente (dans \ $ Hom(\mathcal{E}^{\bullet}(k), \mathcal{E}^{\bullet}(k+k')$ ) :
$$ j_{k,k+k'} \circ \tau_k =  ( _{k+k'}\mathcal{N}^{k'})\circ \tau_{k+k'}\circ  j_{k,k+k'} $$
permet de raisonner directement et redonne (@).}
\begin{equation*}
 j_{k,k+k'}\circ\tau_k = \tau_{k+k'}\circ ( _{k+k'}\mathcal{N}^{k'})\circ j_{k,k+k'}  ,  \tag{@}
 \end{equation*}
 entre morphismes de  faisceaux  \ $ h^{\bullet}(k) \rightarrow h^{\bullet}(k+k') $ \  permet de voir qu'il est inutile de passer au quotient par \ $ \tau_k  \big(\mathbb{H}_{c\,\cap\,S}^{j-1}(Y, \mathcal{E}^{\bullet}(k)), \delta^{\bullet})\big) $ \ et l'on considerera donc simplement le syst\`eme inductif
$$ \big( \mathbb{H}_{c\,\cap\,S}^j(Y, \mathcal{E}^{\bullet}(k)), \delta^{\bullet}) ,  j_{k,k+k'} \big) $$
 De plus les endomorphismes de complexes \ $ _k\mathcal{N}$ \ d\'eduit des endomorphismes \ $ _kN \in End(\mathbb{C}^k) $ \ sont compatibles avec les \ $ j_{k,k+k'}$. 

\bigskip

\subsection {\bf D\'efinition.}

\noindent {\it Pa{\c c}ons-nous sous les hypoth\`eses standards pour la valeur propre 1. Nous d\'efinirons l'espace vectoriel \ $ H^j_{c\,\cap\,S}(0) $ \ comme la limite inductive quand \ $ k \rightarrow \infty$ \ du syst\`eme inductif ci-dessus.
L'action de \ $ _k\mathcal{N} $ \ sur ce syst\`eme inductif  permet de d\'efinir une action monodromique (unipotente) sur sa limite inductive qui sera donn\'ee par \ $T_{c\,\cap\,S} = exp(-2i\pi. \mathcal{N}_{c\,\cap\,S}).$ }

\bigskip

\subsection{}

 Les limites inductives
$$ \underset{k \rightarrow \infty}{lim} \ \mathbb{H}^n_{c}(Y, \mathcal{E}^{\bullet}(k), \delta^{\bullet})$$
et 
$$ \underset{k \rightarrow \infty}{lim} \ \mathbb{H}^n(Y, \mathcal{E}^{\bullet}(k), \delta^{\bullet}))$$
sont respectivement isomorphes \`a \ $H^n_c(0) $ \ et \ $H^n(0)$ \  (voir (1.8) et (1.9) ) ce qui nous fournit des applications \ $\mathbb{C}-$lin\'eaires monodromiques  "naturelles"
$$ can_c^{c\,\cap\,S} : H^n_c(0) \rightarrow H^n_{c\,\cap\,S}(0) \quad {\rm et} \quad can_{c\,\cap\,S} : H^n_{c\,\cap\,S} \rightarrow H^n(0) $$
dont la compos\'ee est l'application usuelle d'oubli de support 
 $$ can : H^n_c(0) \rightarrow H^n(0) .$$
 On remarquera que si la singularit\'e de \ $f$ \ est isol\'ee pour la valeur propre 1 (et donc \
 $S = \lbrace 0 \rbrace$), alors \ $H^n_{c\,\cap\,S}(0)$ \ est isomorphe \`a \ $H^n(0)$ \ via \ $can_{c\,\cap\,S}$.

\bigskip

\subsection{}

 Notre premier r\'esultat va montrer que pour \ $ n \geq 3$ \ le noyau \ $Ker\,\theta$ \ est simplement l'image de l'application \ $ H^n_{c\,\cap\, S}(0) \overset{can_{c\,\cap\,S}}{\longrightarrow} H^n(0)$ \ que nous avons introduite plus haut.
 
 \noindent Le cas \ $n = 2$ \ sera trait\'e s\'epar\'ement au (4.12).

\bigskip

\subsection{\bf Th\'eor\`eme 1.}

\noindent {\it Sous les hypoth\`eses standarts, on a, pour \ $n \geq 2$, la suite exacte : }

$$ 0 \rightarrow H^1_{\lbrace 0 \rbrace}(S, H^{n-1}(0)) \overset{i}{\rightarrow} H^n_{c\,\cap\, S}(0) \overset{can_{c\,\cap\,S}}{\longrightarrow} H^n(0) \overset{\theta}{\rightarrow} H^1(S^*, H^{n-1}(0)) .$$

\bigskip

\noindent \textit{\underline{Preuve}.} 

\smallskip

\noindent Nous utiliserons les lemmes suivants :

\bigskip

\subsection{\bf Lemme.}

 Supposons \ $n\geq 2$ \ et les hypoth\`eses standards v\'erifi\'ees pour la valeur propre 1 de la monodromie de \ $f$. Alors on  a
\begin{align}
\quad & H^j(S^*, \delta\mathcal{E}^i(k)) = 0 , \quad \quad \forall i \in [0,n-2] , \forall j \geq 1 \\
\quad & H^j(S, \delta\mathcal{E}^i(k)) \simeq H^j(Y, \delta\mathcal{E}^i(k)) = 0 , \quad \quad  \forall i \geq 0 , \forall j \geq 1 
\end{align}
\smallskip
\noindent \textit{Preuve.} Nous utiliserons les annulations suivantes pour les faisceaux de cohomologie
\ $ h^j(k)$ \ du complexe \ $( \mathcal{E}^{\bullet}(k),\delta^{\bullet})$ \ qui sont non nuls seulement pour \ $ j = 0,1,n-1,n,n+1$ et qui sont d\'ecrits au (1.2) :
\begin{align}
 H^q(S^*, h^j(k)) = 0 , \qquad \forall j \geq 0 , \forall q \geq 2  \\
 H^p(S, h^j(k)) \simeq H^p(Y, h^j(k)) = 0 , \qquad \forall j \geq 0 , \forall p \geq 1 .
\end{align}
 Pour \ $i = 0$ \ consid\'erons la suite exacte de faisceaux :
\begin{equation*}
 0 \rightarrow h^0(k) \rightarrow \mathcal{E}^0(k) \rightarrow \delta\mathcal{E}^0(k) \rightarrow 0 . \tag{$A_0(k)$} 
 \end{equation*}
Comme \ $ \mathcal{E}^0(k)$ \ est un faisceau fin, on obtient  (3) et (4)  pour \ $ i = 0 .$

\noindent Pour \ $i = 1$ \ on consid\`ere la suite exacte de faisceaux :
\begin{equation*}
 0 \rightarrow \delta\mathcal{E}^0(k) \rightarrow Ker\,\delta^1(k) \rightarrow h^1(k) \rightarrow 0 . \tag{$B_1(k)$}
 \end{equation*}
Elle donne \ $ H^q(S^*, Ker\,\delta^1(k)) = 0 , \forall q \geq 2 $ \ ainsi que l'annulation de  \\ $H^p(S, Ker\,\delta^1(k))$ \ et de \ $H^p(Y, Ker\,\delta^1(k)), \ \forall p \geq 1, $ \ grace \`a (3) et (4) \\ pour \ $i = 0$.

\noindent La suite exacte de faisceaux :
\begin{equation*}
  0 \rightarrow  Ker\,\delta^1(k) \rightarrow \mathcal{E}^1(k) \rightarrow \delta\mathcal{E}^1(k) \rightarrow 0  \tag{$A_1(k)$}
  \end{equation*}
donne alors les annulations (3) et (4) pour \ $i = 1$ \ puisque \ $\mathcal{E}^1(k) $ \ est fin.

\noindent Supposons que l'on a montr\'e les annulations (3) et (4)  pour \ $ i \in [1, n-3]$ ; nous allons en d\'eduire ces m\^emes annulations pour \ $i+1$, ce qui prouvera par r\'ecurrence les annulations (3) et (4) pour \ $i\in [0,n-2]$. 

\noindent Consid\`erons la suite exacte  de faisceaux :
$$ 0 \rightarrow \delta \mathcal{E}^i(k) \rightarrow \mathcal{E}^{i+1}(k) \rightarrow  \delta \mathcal{E}^{i+1}(k) \rightarrow 0 .$$
Elle permet de conclure aux annulations (3) et (4) pour \ $i+1$ \ grace \`a la finesse du faisceau
 \ $ \mathcal{E}^{i+1}(k) .$

\noindent Pour prouver les annulations (4) pour \ $i \geq n-1 $ \ on utilise successivement les suites exactes :
\begin{align*}
&  0 \rightarrow \delta\mathcal{E}^{q-1}(k) \rightarrow Ker\,\delta^{q}(k) \rightarrow h^{q}(k) \rightarrow 0 \tag{$B_q(k)$} \\
& 0 \rightarrow Ker\,\delta^{q}(k) \rightarrow \mathcal{E}^{q}(k) \rightarrow \delta\mathcal{E}^{q}(k) \rightarrow 0 \tag{$A_q(k)$}
\end{align*}
pour \ $q \geq n-1 $.  $\hfill \blacksquare $

\bigskip

\subsection {\bf Lemme.}

 Sous les hypoth\`eses standarts avec \ $n \geq 4$ \  et \ $k \gg1$ \ la fl\`eche naturelle
$$e_k : = \frac{H^0(S^*, Ker\,\delta^{n-1}(k))}{H^0(S, Ker\,\delta^{n-1}(k)) + \delta H^0(S^*, \mathcal{E}^{n-2}(k))} \rightarrow H^1_{\lbrace 0 \rbrace}(S, H^{n-1}(0)) $$
est un isomorphisme d'espaces vectoriels monodromiques.

\noindent Pour \ $n = 3 $ \ la fl\`eche analogue est surjective et de noyau isomorphe \`a \ $H^1(S^*, h^1(k)).$ Donc on obtient encore un isomorphisme monodromique, quitte \`a passer \`a la limite inductive sur \ $k$\footnote{En fait, de fa{\c c}on plus pr\'ecise, l'application \ $ j_{k,k+1}$ \ induit \ $0$ \ sur \ $h^1(k)$.}.

\noindent Pour \ $n = 2 $ \ la fl\`eche analogue, apr\`es passage \`a la limite inductive, est encore surjective et a un noyau isomorphe \`a \ $H^1(S^*, \mathbb{C})$.

\bigskip

\noindent{\it Preuve.} Comme on a suppos\'e \ $n \geq  3 $ \ l'application
$$ r^{n-1}(k) : h^{n-1}(k) \rightarrow  H^{n-1}(0)  $$
est un isomorphisme pour \ $ k \gg 1 $ \ d'apr\`es  (1.6). La suite exacte
\begin{equation*}
  0 \rightarrow \delta\mathcal{E}^{n-2}(k) \rightarrow Ker\,\delta^{n-1}(k) \rightarrow h^{n-1}(k) \rightarrow 0  \tag{$B_{n-1}(k)$}
  \end{equation*}
donne, puisque \ $ H^1(S^*,  \delta \mathcal{E}^{n-2}(k)) = 0 $ \ d'apr\`es (4.7), la surjectivit\'e de la fl\`eche \ $q : H^0(S^*, Ker\,\delta^{n-1}(k)) \rightarrow H^0(S^*, h^{n-1}(k)) $, d'o\`u une fl\`eche naturelle et surjective
$$ H^0(S^*, Ker\,\delta^{n-1}(k)) \rightarrow  H^1_{\lbrace 0 \rbrace}(S, H^{n-1}(0)) $$
puisque l'on a
$$ H^1_{\lbrace 0 \rbrace}(S, H^{n-1}(0)) \simeq \frac{H^0(S^*, H^{n-1}(0))}{H^0(S, H^{n-1}(0))} $$
gr\^ace \`a l'annulation de \ $ H^1(S, H^{n-1}(0)) $\footnote{Rappelons que \ $H^{n-1}(0)$ \ est un syst\`eme local sur \ $S^*$.}.

\noindent Il nous reste donc \`a identifier le noyau de cette surjection pour conclure la preuve. Comme les espaces \ $ H^0(S, Ker\,\delta^{n-1}(k))$ \ et \ $\delta H^0(S^*, \mathcal{E}^{n-2}(k))$ \ ont des images qui sont clairement dans le noyau, il suffit de montrer que tout \'el\'ement \ $ \alpha \in H^0(S^*, Ker\,\delta^{n-1}(k))$ \ dont l'image dans \ $ H^0(S^*, h^{n-1}(k)) $ \ se prolonge \`a \
$H^0(S, H^{n-1}(0))$ \ est de cette forme. Comme l'application 
$$ H^0(S, Ker\,\delta^{n-1}(k)) \rightarrow H^0(S, h^{n-1}(k)) $$ 
est surjective\footnote{Grace \`a l'annulation de \ $ H^1(S, \delta\mathcal{E}^{n-2}(k))$ \ qui est prouv\'ee au (4.7).} il existe \ $\beta \in H^0(S,  Ker\,\delta^{n-1}(k))$ \ dont l'image dans \ $H^0(S, h^{n-1}(k)) $ \ prolonge \`a \ $S$ \ l'image de \ $\alpha$ \ dans \ $H^0(S^*, h^{n-1}(k)) $. Donc la restriction de \ $\beta$ \ \`a \ $S^*$ \ est \'egale \`a \ $\alpha$ \ modulo le noyau \ $H^0(S^*, \delta\mathcal{E}^{n-2}(k))$. On conclut  pour \ $n \geq 4$ \ en remarquant que l'\'egalit\'e \ $ Ker\,\delta^{n-2}(k) = \delta \mathcal{E}^{n-3}(k)$ \ (qui est vraie pour \ $n \geq 4$) \ donne, grace \`a l'annulation de \ $ H^1(S^*, \delta \mathcal{E}^{n-3}(k)) $ \ montr\'ee au (4.7), la surjectivit\'e de l'application naturelle
$$\delta :  H^0(S^*, \mathcal{E}^{n-2}(k)) \rightarrow H^0(S^*, \delta \mathcal{E}^{n-2}(k))$$
d'o\`u notre assertion.

\noindent Pour \ $n = 3 $ \ la suite exacte de cohomologie de la suite exacte de faisceaux :
\begin{equation*}
 0 \rightarrow Ker\,\delta^1(k) \rightarrow \mathcal{E}^{1}(k)) \rightarrow \delta \mathcal{E}^{1}(k))\rightarrow 0 \tag{$A_1(k)$}
 \end{equation*}
montre que le conoyau de la fl\`eche 
$$\delta :  H^0(S^*, \mathcal{E}^{1}(k)) \rightarrow H^0(S^*, \delta \mathcal{E}^{1}(k))$$
est \'egal \`a \ $ H^1(S^*, Ker\,\delta^1(k))$ \ qui est isomorphe \`a \  $H^1(S^*, h^1(k))$ \ comme on le constate en utilisant la suite exacte longue de cohomologie de la suite exacte :
\begin{equation*}
 0 \rightarrow \delta\mathcal{E}^0(k) \rightarrow Ker\,\delta^1(k)) \rightarrow h^1(k) \rightarrow 0 \tag{$B_1(k)$}
 \end{equation*}
et (4.7). D'o\`u notre assertion pour \ $n = 3$, puisque \ $ \underset{k \rightarrow \infty} {\lim}h^1(k) = 0 $ \ dans ce cas.

\noindent Pour \ $n = 2 $ \  on a \ $ \underset{k \rightarrow \infty}{lim} \ h^1(k) \simeq H^1(0) $ \ et la suite exacte 
 $$ 0 \rightarrow \delta\mathcal{E}^0(k) \rightarrow Ker\,\delta^1(k)\rightarrow h^1(k) \rightarrow 0 $$
 donne la surjectivit\'e de \ $ H^0(S^*, Ker\,\delta^1(k)) \rightarrow H^0(S^*, h^1(k)) $ \ puisque 
  $$ H^1(S^*, \delta\mathcal{E}^0(k)) \simeq H^2(S^*, h^0(k)) \simeq 0 .$$
  Donc la fl\`eche \ $\underset{k \rightarrow \infty}{lim} \ (e_k)$ \ est surjective.
  
  \noindent Son noyau va \^etre \'egal \`a la limite inductive du quotient
  
  $$ \frac{H^0(S^*, \delta\mathcal{E}^0(k))}{\big(H^0(S, Ker\,\delta^1(k)) + \delta H^0(S^*, \mathcal{E}^0(k)) \big)\cap H^0(S^*, \delta\mathcal{E}^0(k)) } .$$
  Mais on a \ $ H^0(S, Ker\,\delta^1(k)) \cap  H^0(S^*, \delta\mathcal{E}^0(k)) \simeq H^0_{\lbrace 0 \rbrace}(S, h^1(k)) $ \ qui a une limite inductive nulle\footnote{Car \ $H^1(0)$ \ n'a pas de section non nulle \`a support l'origine.};  la limite inductive du d\'enominateur co{\"i}ncide donc avec celle de \ $\delta H^0(S^*, \mathcal{E}^0(k))$. La suite exacte
  $$ 0 \rightarrow h^0(k) \rightarrow \mathcal{E}^0(k)) \rightarrow \delta\mathcal{E}^0(k)) \rightarrow 0 $$
  montre alors que le noyau est isomorphe \`a
  $$ \underset{k \rightarrow \infty} \lim \  \frac{ H^0(S^*, \delta\mathcal{E}^0(k))}{\delta H^0(S^*, \mathcal{E}^0(k))} \simeq  \underset{k \rightarrow \infty} \lim H^1(S^*, h^0(k)) \simeq H^1(S^*, \mathbb{C}).$$
  D'o\`u notre assertion pour \ $n = 2 $.  $ \hfill \blacksquare$

\bigskip

\subsection{Construction de l'application \ $i$.}
Comme nous supposons \ $n \geq 3$ \ il nous suffit, d'apr\`es le lemme 4.8 , de construire une application lin\'eaire monodromique
$$ H^0(S^*, Ker\,\delta^{n-1}(k)) \rightarrow H^n_{c\,\cap \,S}(0) $$
et de v\'erifier que son noyau co{\"i}ncide avec \ $ \delta H^0(S^*, \mathcal{E}^{n-2}(k)) + H^0(S, Ker\,\delta^{n-1}(k))$\footnote{quitte \`a passer \`a la limite inductive sur \ $k$ \ pour \ $n = 3$.} \ et que son image  co{\"i}ncide avec \ $ Ker\,can_{c\,\cap\,S} $.

\noindent Soit donc \ $\alpha \in H^0(S^*, Ker\,\delta^{n-1}(k))$ \ et fixons une fonction  \ $\chi \in \mathcal{C}^{\infty}(X) $ \ qui v\'erifie \ $ \chi \equiv 1 $ \ pr\`es de \ $ S \setminus B(0,r) $ \ et qui est nulle d\`es que l'on s'\'eloigne de \ $ S \setminus B(0,r) $ \ de sorte que \ $ Supp\,( d\chi )\cap S^* $ \ soit compact. Posons alors \ $ i(\alpha) : = d\chi\wedge \alpha = \delta(\chi.\alpha)$. Il est imm\'ediat de v\'erifier que la classe ainsi d\'efinie dans \ $H^n_{c\,\cap \,S}(0) $ \ ne d\'epend ni de \ $r \in ]0,\varepsilon[$\footnote{O\`u ici \ $\varepsilon$ \ est le rayon de la boule de Milnor \ $X$ \ que l'on consid\`ere.} ni du choix de \ $\chi$. De plus, si on a \ $\alpha = \delta\beta$, avec \ $\beta \in  H^0(S^*, \mathcal{E}^{n-2}(k))$, on aura 
$$ \delta (\chi.\alpha) = \delta(-d\chi\wedge \beta) $$
et \ $ d\chi\wedge \beta $ \ est \`a support dans \ $ c\,\cap\,S $.

\noindent Si \ $\alpha$ \ est la restriction \`a \ $S^*$ \ de \ $\tilde{\alpha} \in H^0(S, Ker\,\delta^{n-1}(k))$ \ on aura \ $ \delta(\chi.\alpha) = \delta((1-\chi).\tilde{\alpha}) $ \ et \ $(1-\chi).\tilde{\alpha} $ \ est \`a support dans \ $ c\,\cap\,S $. Ceci montre que l'application \ $i$ \ est bien d\'efinie pour \ $n \geq 3$. Elle est monodromique car on a
$$ d\chi \wedge \mathcal{N}(\alpha) = \mathcal{N}(d\chi\wedge\alpha).$$
Pour montrer l'injectivit\'e de \ $i$, consid\'erons \ $\alpha \in H^0(S^*, Ker\,\delta^{n-1}(k))$ \ tel que\ $i(\alpha) = \delta(\chi.\alpha) = \delta \beta$ \ avec \ $\beta \in \Gamma_{c\,\cap\,S}(Y, \mathcal{E}^{n-1}(k)) $\footnote{Quitte \`a choisir \ $k \gg1 $.}. Alors \  $ \chi.\alpha -\beta$ \ est dans \ $ \Gamma(Y, Ker\,\delta^{n-1}(k)) $ \ et on obtient ainsi un prolongement de \ $\alpha$ \ \`a \ $H^0(S, Ker\,\delta^{n-1}(k))$.

\noindent Pour identifier l'image de \ $i$ \ commen{\c c}ons par remarquer que \ $ can_{c\,\cap\,S}\circ i = 0 $. En effet, par d\'efinition, \ $i(\alpha) = \delta(\chi.\alpha)$ \ et \ $ \chi.\alpha \in \Gamma(Y, \mathcal{E}^{n-1}(k))$. D'o\`u l'inclusion de l'image de \ $i$ \ dans \ $ Ker\,can_{c\,\cap\,S}$. R\'eciproquement, soit \ $ u$ \ dans \ $ \Gamma(Y, \mathcal{E}^{n-1}(k)) $ \  telle que le support de \ $\delta u $ \ soit dans \ $ c\, \cap\,S $. On a donc un compact \ $ K$ \  de \ $S$ \ tel qu'au voisinage de \ $S \setminus  K$ \   on ait \ $\delta u = 0$. Alors \ $u$ \ d\'efinit un \'el\'ement de 
 $$ H^0(S\setminus K, h^{n-1}(k)) \simeq H^0(S^*, H^{n-1}(0)) $$ 
 pour \ $ k \gg 1 $.
 L'image par \ $i$ \ de cette classe est repr\'esent\'ee par \ $ \delta(\chi.u) .$ Mais \ $ \delta u - \delta( \chi.u) = \delta((1-\chi).u) $ \ et on a \ $ (1-\chi).u \in \Gamma_{c\,\cap\,S}(Y, \mathcal{E}^{n-1}(k))$ , ce qui montre que la classe initiale \ $\delta u$ \ est bien dans l'image de \ $i$.

\bigskip

\subsection{Le noyau de l'application \ $\theta$.}

Rappelons maintenant comment calculer de l'application \ $\theta.$ D'apr\`es [B.91] p. 428 on a la suite exacte de faisceaux suivante, pour \ $k \gg1$
$$h^{n-2}(k) \overset{\tau_k}{\rightarrow}  h^{n-1}(k) \overset{\tau_k}{\rightarrow} h^{n}(k) \rightarrow H^n(0) \rightarrow 0$$
Soit  \ $ w \in \Gamma(Y, Ker\,\delta^n(k))$. On peut donc \'ecrire, localement le long de \ $S^*$, puisque le faisceau \ $H^n(0)$ \ est concentr\'e \`a l'origine,
$$ w\vert_{X_{\sigma}} = \tau_{k}(\alpha_{\sigma}) + \delta \beta_{\sigma} $$
avec \ $ \delta \alpha_{\sigma}= 0 $ (comparer avec [B.91] p.439). On obtient, comme dans loc. cit. que les \ $\alpha_{\sigma}$ \ se recollent\footnote{Pour \ $n \geq 4$ ;  pour \ $n = 3$ \ le recollement a lieu seulement dans \ $H^2(0)$.}  pour donner une section sur \ $S^*$ \ du faisceau  \ $h^{n-1}(k)$ \ dont l'image dans \ $H^1_{\lbrace 0 \rbrace}(S, h^{n-1}(k)) $ \ est ind\'ependante des choix effectu\'es. 
En utilisant ce qui pr\'ec\`ede, on peut trouver \ $ \alpha \in \Gamma(S^*, Ker\,\delta^{n-1}(k)) $ \ telle que \ $w - \tau_k(\alpha) $ \ soit localement \ $\delta-$exacte le long de \ $S^*$. A partir des \'ecritures locales
$$ w - \tau_k(\alpha)\vert_{X_{\sigma}} = \delta\beta_{\sigma} $$
on d\'efinit  \ $\theta(w)$ \ par la classe dans $$ H^1(S^*, h^{n-1}(k))\simeq H^1(S^*, H^{n-1}(0)) $$
pour \ $ k \gg 1$, du  1-cocycle \ $ \beta_{\sigma} - \beta_{\sigma'} $. On v\'erifie facilement que ceci est ind\'ependant des choix effectu\'es.

\smallskip

\noindent Nous voulons montrer maintenant que le noyau de \ $\theta$ \ est l'image de \ $can_{c\,\cap\,S}$. D'abord si \ $w$ \ a son support dans \ $c\,\cap\,S$ \ on peut choisir pour chaque \ $ \sigma \in S \setminus K , \alpha_{\sigma} = 0 $ \ ainsi que \  $ \beta_{\sigma} = 0 .$ La section \ $\alpha$ \ sera donc nulle sur \ $S\setminus K$ \ et comme \ $H^{n-1}(0) $ \ est un syst\`eme local sur \ $S^*$, on aura \ $\alpha = 0 $ \ et donc \ $\theta[w]$ \ sera nul dans \ $ H^1(S\setminus K, H^{n-1}(0)) $. A nouveau le fait que \ $H^{n-1}(0)$ \ soit un syst\`eme local sur \ $S^*$ \ donne la nullit\'e de \ $\theta[w]$.

\noindent R\'eciproquement, supposons que \ $\theta[w] = 0 $. Utilisons maintenant un \ $k$ \ assez grand pour assurer que \ $\tilde{ob}_k \equiv 0 $\footnote{Comme on l'a d\'eja fait remarquer plus haut la formule \ $ \tilde{ob}_{k+k'} = \tilde{ob}_k \circ \mathcal{N}^{k'} $ \ de [B.91]  est encore valable ici ; donc pour \ $k$ \ assez grand on a \ $\tilde{ob}_k \equiv 0 .$}. Cela signifie que l'on peut \'ecrire globalement sur \ $S^*$ 
$$ w = \delta\beta $$
avec \ $ \beta \in \Gamma(S^*, \mathcal{E}^{n-1}(k))$. Soit  $\chi\in\mathcal{C}^{\infty}(Y)$ \ v\'erifiant \ $\chi \equiv 1 $ \ pr\`es de \ $S \setminus B(0,r)$ \ et s'annulant identiquement d\`es que l'on s'\'eloigne de \ $S \setminus B(0,r)$. Alors \ $ w - \delta(\chi.\beta) $ \ est un repr\'esentant de la classe de \ $w$ \ qui est \`a support dans \ $ c\,\cap\,S $. Donc cette classe est bien dans l'image de  \ $can_{c\,\cap\,S}$. 

\noindent Ceci ach\`eve la preuve du th\'eor\`eme 4.6. $\hfill \blacksquare$

\bigskip

\subsection{}

Une cons\'equence importante mais imm\'ediate pour \ $n \geq 3$ \ du th\'eor\`eme 4.6 et du corollaire 3  du th\'eor\`eme 3.3.1 est l'\'egalit\'e de la codimension de \ $Im (\theta) $ \ dans l'orthogonal pour \ $\tilde{h}$ \ de \ $H^0(S, H^{n-1}(0)) $ avec \ $\dim H^n_{c\, \cap\,S}(0) - \dim H^n(0)$. Ceci montre d\'eja que la dimension de l'espace vectoriel \ $ H^n_{c\, \cap\,S}(0)$ \ est au moins \'egale \`a celle de \ $H^n(0)$.
Notre objectif va \^etre maintenant de montrer que l'on a en fait \'egalit\'e de ces deux dimensions. Ceci n'est pas \'evident et sera obtenu via la construction d'une application de variation injective (pour \ $n \geq 3$ \ au moins) commutant aux monodromies :
$$  var : H^n_{c\,\cap\,S}(0) \rightarrow H^n_c(0) .$$
On en d\'eduira alors que \ $ \dim  H^n_{c\, \cap\,S}(0) \leq \dim H^n_c(0)  = \dim H^n(0) $ \ pour \ $n \geq 3$ \  cette derni\`ere \'egalit\'e r\'esultant de la dualit\'e de Poincar\'e sur la fibre de Milnor de \ $f$ \ \`a l'origine.
On en conclura, toujours pour \ $n \geq 3$, que, de plus, cette application de variation est bijective, ce qui sera la clef de la non-d\'eg\'en\'erescence de la forme hermitienne canonique "g\'en\'eralis\'ee" :
$$ \mathcal{H} :   H^n_{c\,\cap\,S}(0) \times H^n(0) \rightarrow \mathbb{C}  $$
qui sera d\'efinie plus loin.\\
Ceci montrera \'egalement que \ $Im(\theta)$ \ co{\"i}ncide avec l'orthogonal dans \ $ H^1(S^*, H^{n-1}(0))$ \ pour \ $\tilde{h}$ \ de \ $H^0(S, H^{n-1}(0))$.

\subsection{Le cas \ $n =2$.}

\subsubsection{}

D\'efinissons sous les hypoth\`eses standards pour la valeur propre 1 et pour \ $n =2$ 
$$ \mathcal{K} : = \underset{k \rightarrow\infty}{lim} \ \frac{H^0(S^*, Ker\,\delta^{1}(k))}{H^0(S, Ker\,\delta^{1}(k)) + \delta H^0(S^*, \mathcal{E}^{0}(k))} .$$
Il r\'esulte du Lemme (4.8) que l'on a une suite exacte monodromique :

\begin{equation*}
0 \rightarrow H^1(S^*, \mathbb{C}) \overset{j}{\rightarrow} \mathcal{K} \overset{e}{\rightarrow} H^1_{\lbrace 0 \rbrace}(S, H^1(0)) \rightarrow 0
\end{equation*}
obtenue par passage \`a la limite inductive sur \ $k$. 

\noindent On remarquera que le sous-espace \ $H^1(S^*, \mathbb{C}) $ \ est dans la partie invariante par la monodromie de \ $\mathcal{K}$. Ceci r\'esulte du fait que la monodromie agit comme l'identit\'e sur les faisceaux \ $h^0(k)$.
 
 \subsubsection{Th\'eor\`eme 1, cas \ $n = 2$.}
 
\textit{ Sous les hypoth\`eses standarts pour la valeur propre 1, on a pour \ $n = 2$ \ la suite exacte monodromique :
 $$ 0 \rightarrow \mathcal{K} \overset{i}{\rightarrow} H^2_{c\,\cap\,S}(0)\overset{can_{c\,\cap\,S}}{ \longrightarrow} H^2(0) \overset{\theta}{\rightarrow} H^1(S^*, H^1(0)) $$
 o\`u l'application \ $i$ \ est d\'eduite du lemme (4.8).}
 
 \bigskip
 
 \noindent La preuve est tout \`a fait analogue \`a celle donn\'ee dans le cas \ $n \geq 3$ \ modulo les adaptations que l'on vient de faire pour tenir compte du r\'esultat diff\'erent du lemme (4.8) dans ce cas, et modulo les modifications suivantes :
 \begin{itemize}

 \item{} Pour montrer que \ $ Ker\,can_{c\,\cap\,S} \subset Im(i) $ \ on consid\`ere  \ $ u \in \Gamma(Y, \mathcal{E}^{1}(k)) $ \ qui est telle que le support de \ $\delta u $ \ soit dans \ $ c\, \cap\,S $. On a donc un compact \ $ K$ \  de \ $S$ \ tel qu'au voisinage de \ $S \setminus  K$ \   on ait \ $\delta u = 0$. Alors \ $u$ \ d\'efinit un \'el\'ement de \ $ H^0(S\setminus K,  Ker\,\delta^{1}(k))$. Comme on a
 $$H^0(S\setminus K, H^1(0)) \simeq H^0(S^*, H^1(0))$$ 
 son image dans \ $H^0(S\setminus K, H^1(0))$ \  donne une section de \ $H^0(S^*, H^1(0))$ \ ce qui nous fournit une section \ $ v \in H^0(S^*, h^1(k))$ \ qui prolonge l'image de \ $u$ \ dans \ $H^0(S\setminus K, h^1(k))$ \ pour \ $k \gg1 $.  Comme on a \ $H^1(S^*, \delta\mathcal{E}^0(k)) = 0 $ \ d'apr\`es le lemme (4.8), on en d\'eduit l'existence de \ $u_1 \in H^0(S^*, Ker\, \delta^1(k))$ \ dont la restriction \`a \ $ S\setminus K $ \ est \'egale \`a \ $u$. Comme on a \ $H^0(S\setminus K, h^0(k)) \simeq H^0(S^*, h^0(k))$ \ on peut modifier \ $u_1$ \ pour que sa restriction \`a \ $ S\setminus K $ \ soit exactement  \'egale \`a \ $u$.
 On conclut alors comme dans le cas \ $n \geq 3.$

 \item{} Pour l'\'etude du noyau de \ $\theta$, on doit remplacer la suite exacte consid\'er\'ee pour \ $ n \geq 3$ \  par la suite exacte :
 $$ 0 \rightarrow \frac{h^1(k)}{\tau_k h^0(k)} \rightarrow h^2(k) \rightarrow H^2(0) \rightarrow 0 .$$
 et utiliser l'isomorphisme sur \ $S^* : \frac{h^1(k)}{\tau_k h^0(k)} \simeq H^1(0) $.
 
 \noindent Localement pr\`es de \ $\sigma \in S^*$ \ on a encore une \'ecriture
 $$ w\vert_{X_{\sigma}}  - \tau_k(\alpha_{\sigma}) = \delta\beta_{\sigma} $$
avec \ $\delta(\alpha_{\sigma}) = 0$. Mais maintenant le recollement des \ $\alpha_{\sigma}$ \ n'a lieu sur \ $S^*$ \ que dans le quotient \ $\frac{h^1(k)}{\tau_k h^0(k)}$. On ne peut donc, en g\'en\'eral, trouver un \'element \ $\alpha \in \Gamma(S^*, Ker\,\delta^1(k)) $ \ tel que \ $ w - \tau_k(\alpha)$ \ soit localement \ $\delta-$exact le long de \ $S^*$. Mais ceci ne change rien \`a la preuve de l'inclusion \ $  Im(can_{c\,\cap\,S})\subset Ker\,\theta  .$

 \item{} Pour voir que \ $ Ker\,\theta \subset Im(can_{c\,\cap\,S})$ \ l'argument permettant d'\'ecrire
 $ w = \delta \beta$ \ globalement pr\`es de \ $S^*$ \ avec \ $\beta \in \Gamma(S^*, \mathcal{E}^1(k))$  n'est plus correct. L'annulation de \ $\tilde{ob}_k $ \ pour \ $k \gg 1$ \ reste valable, mais donne seulement l'existence, pour  \ $k \gg 1$,  de \ $\gamma_{\sigma} \in \Gamma(X_{\sigma}, Ker\,\delta^0(k))$, v\'erifiant \ $\alpha_{\sigma} = \tau_k( \gamma_{\sigma})$ \ dans \ $\Gamma(X_{\sigma}, Ker\,\delta^1(k))$ \ quitte \`a modifier le choix des \ $\beta_{\sigma}$. Mais quitte \`a changer \ $k$ \ en \ $k+1$ \ en appliquant \ $j_{k,k+1}$, on rend \ $\tau_k(\gamma_{\sigma}) \ \delta-$exact, et donc \ $j_{k,k+1}(w)$ \ devient localement \ $\delta-$exacte le long de \ $S^*$. Alors l'argument du cas \ $n \geq 3$ \ s'applique et permet de conclure \`a l'\'egalit\'e  \ $ Ker\,\theta = Im(can_{c\,\cap\,S}).$
  \end{itemize}
 
 \subsubsection{}
 Montrons que dans le cas \ $n = 2$ \ le th\'eor\`eme 1 permet  de prouver une in\'egalit\'e de dimension analogue \`a celle du cas \ $ n \geq 3$\footnote{En fait, dans le cas \ $n = 2 $ \ on doit syst\'ematiquement remplacer l'espace \ $H^2_{c\,\cap\,S}(0) $ \ par son quotient par \ $H^1(S^*, \mathbb{C})$ \ et donc \ $ \mathcal{K}$ \ par son quotient par \ $H^1(S^*, \mathbb{C})$ \ pour avoir les m\^emes r\'esultats que pour \ $n \geq 3 $.}
 $$ \dim_{\mathbb{C}} \ H^2_{c\,\cap\,S}(0) \geq \dim_{\mathbb{C}}\  H^2(0) + \dim_{\mathbb{C}}\ H^1(S^*, \mathbb{C}).$$
 En effet la suite exacte du th\'eoreme donne
 $$ \dim\,  H^2_{c\,\cap\,S}(0) + \dim\, Im(\theta) = \dim\, \mathcal{K} + \dim\, H^2(0) $$
 et la suite exacte du (4.12.1) donne
$$  \dim\, H^1(S^*, \mathbb{C}) + \dim\, H^1_{\lbrace 0 \rbrace}(S, H^1(0))  =  \dim\, \mathcal{K} .$$
L'inclusion de \ $Im(\theta)$ \ dans l'orthogonal pour \ $\tilde{h}$ \ de \ $ H^0(S, H^1(0))$ \ donne, puisque \ $H^0(S^*, H^1(0))$ \ est le dual de \ $  H^1(S^*, H^1(0))$, l'in\'egalit\'e 
 $$ \dim\, Im(\theta) +  \dim\, H^0(S, H^1(0)) \leq \dim\, H^1(S^*, H^1(0)) .$$
 D'o\`u notre assertion, puisque 
 $$ H^1_{\lbrace 0 \rbrace}(S, H^1(0)) \simeq \frac{H^0(S^*, H^1(0))}{H^0(S, H^1(0))}.$$

\bigskip


\section{La suite exacte longue de J. Leray g\'en\'eralis\'ee.}

\bigskip

 \noindent Le premier ingr\'edient dans la construction de notre variation sera une g\'en\'eralisation encore un peu plus sophistiqu\'ee que celle de [B.97] du r\'esidu de J. Leray.

\bigskip

\subsection{ Th\'eor\`eme.}

\noindent\textit{ Soit \ $Z$ \ une vari\'et\'e complexe connexe de dimension \ $n+1$  \ et soit \ $Y $ \ une hypersurface ferm\'ee d'int\'erieur vide dans \ $Z$.  Soit \ $S$ \ un sous-ensemble analytique ferm\'e et d'int\'erieur vide de \ $Y$ \ tel qu'en chaque point \ $y$ \ de \ $Y\setminus S$,  l'hypersurface \ $Y$ \ admette une \'equation locale r\'eduite \ $f_y$ \ dont la monodromie n'admet pas la valeur propre 1 (dans son action sur la cohomologie r\'eduite de la fibre de Milnor de \ $f_y$ \ en \ $y$ ).}

\noindent \textit{Notons par \ $mod\,S $ \ la famille (paracompactifiante) des ferm\'es de \ $Z$ \ qui ne rencontrent pas  \ $S$.}

\noindent \textit{ Alors on a la suite exacte longue " de Leray" \`a supports \ $mod\,S$ :}
\begin{equation*}
\cdots \rightarrow H^{q}_{mod\,S}(Z, \mathbb{C}) \overset{i^q}{\rightarrow} H^{q}_{mod\,S}(Z\setminus Y, \mathbb{C}) \overset{Res^q}{\rightarrow} H^{q-1}_{mod\,S}(Y, \mathbb{C}) \overset{\partial^q}{\rightarrow} H^{q+1}_{mod\,S}(Z, \mathbb{C}) \rightarrow \cdots 
\end{equation*}
\textit{ Si on dispose d'une \'equation glabale r\'eduite \ $f \in \mathcal{O}(Z)$ \ de \ $Y$ \ dans \ $Z$, pour toute classe \ $[w] \in H^q_{mod\,S}(Z\setminus Y, \mathbb{C})$ \ repr\'esent\'ee par une forme semi-m\'eromorphe d-ferm\'ee  \ $w$, \`a p\^oles dans \ $Y$ \ et \`a support dans \ $mod\,S$, il existe un voisinage ouvert \ $\mathcal{W}$ \ de \ $Y$ \ dans \ $Z$ \ et des formes \ $\eta \in \mathcal{C}^{\infty}(\mathcal{W})^{q-1} $ \ et \ $\omega \in \mathcal{C}^{\infty}(\mathcal{W})^q$ \  d-ferm\'ees sur \ $\mathcal{W}$ \ \`a supports dans \ $mod\,S \vert_{\mathcal{W}}$ \ et une forme semi-m\'eromorphe \ $\alpha$ \ de degr\'e \ $q-1$ \ \`a p\^oles dans \ $Y$ \ et \`a support dans  \ $mod\,S \vert_{\mathcal{W}}$ \ v\'erifiant  :}
\begin{equation*}
w =\frac{df}{f} \wedge \eta + \omega + d\alpha \quad {\rm sur} \quad \mathcal{W}\setminus Y  \tag{@}
\end{equation*}
\textit{ La classe \ $ [\eta\vert_{Y}] \in H^{q-1}_{mod\,S}(Y, \mathbb{C})$ \ est alors l'image par \ $Res^q$ \ de la classe }
$$[w] \in H^q_{mod\,S}(Z \setminus Y, \mathbb{C}).$$

\bigskip

\subsection{Remarque.}

Le seul point non trivial dans la suite exacte ci-dessus est bien sur l'identification du groupe de cohomologie \`a support \ $ H^{q+1}_{Y, mod\,S}(Z, \mathbb{C}) $ \ avec le groupe \ $H^{q-1}_{mod\,S}(Y, \mathbb{C}) $ \ qui est analogue \`a ce qui se passe dans le cas d'une hypersurface lisse. La pr\'esence de la famille de supports \ $ mod\,S $ \ \ est relativement anodine puisqu'elle est paracompactifiante et que sa restriction \`a \ $Y$ \ est constitu\'ee des ferm\'es de \ $Y$ \ qui sont dans  $ mod\,S $. On notera cependant que cet \'enonc\'e peut s'appliquer en pr\'esence d'un lieu  singulier  pour l'hypersurface \ $Y$ \ qui est bien plus gros que \ $S$.

\bigskip

\subsection{D\'emonstration.}

La preuve repose sur deux faits :
\begin{itemize}
\item 1)  Le quasi-isomorphisme, d\^u \`a A. Grothendieck (voir [Gr. 65]) :
$$ Rj_*j^* \underline{\mathbb{C}}_Z \simeq ( \mathcal{C}^{\infty}_Z[*Y]^{\bullet}, d^{\bullet} )$$
qui montre que le complexe des formes semi-m\'eromorphes \`a p\^oles dans \ $Y$ \ donne une incarnation "fine"  du complexe \ $ Rj_*j^* \underline{\mathbb{C}}_Z$ .
\item 2) La d\'eg\'en\'erescence de la suite spectrale :
$$ E_2^{p,q} : = H^q_{mod\,S}(Z, \underline{H}_Y^p(\mathbb{C})) $$
qui r\'esulte du fait que, sur \ $ Z\setminus S $, on a \ $ \underline{H}^p_Y(\mathbb{C}) = 0 $ \  pour \ $p \not= 2 $ \ et \ $ \underline{H}^2_Y(\mathbb{C}) \simeq \underline{\mathbb{C}}_Y .\frac{df}{f} $ , combin\'e avec le fait que si le faisceau \ $\mathcal{F}$ \ est nul sur \ $Z\setminus S$ \ alors on a \ $ H^p_{mod\,S}(Z, \mathcal{F}) = 0 , \forall p \geq 0 $.
\end{itemize}
Le calcul de  \ $ \underline{H}^p_Y(\mathbb{C})  $ \ sur \ $Z\setminus S$ \ vient de l'absence de la valeur propre 1  pour la monodromie locale (voir la remarque finale de [B.97] ) d'une \'equation locale r\'eduite \ $f$ \ de \ $Y$, qui donne l'annulation des groupes \ $ H^i(D^*, R^j f_*(\mathbb{C}_Z)) $ \ pour tout \ $ i \geq 0 $ \  quand \ $ j \not= 0$.

\noindent La suite exacte longue de Leray en r\'esulte alors d'apr\`es [Go.58] th. 4.10.1 puisque la d\'eg\'en\'erescence de la suite spectrale  donne les isomorphismes
$$ H^{q+1}_{Y, mod\,S}(Z, \mathbb{C}) \simeq H^{q-1}_{mod\,S}(Y, \mathbb{C}) .$$
Pour prouver \ $(@)$ \ repr\'esentons \ $Res^q[w] \in H^{q-1}_{mod\,S}(Y, \mathbb{C}) $ par une forme \ $\eta \in \mathcal{C}^{\infty}$ \ de degr\'e \ $q-1$ \ d-ferm\'ee \`a support dans \ $mod\,S$\footnote{Ceci est possible grace \`a [Go.58] th. 4.11.1 et \`a de Rham.} sur un voisinage ouvert \ $\mathcal{W}$ \ de \ $Y$ \ assez petit. Alors la classe \ $ [w - \frac{df}{f} \wedge \eta ] \in H^q_{mod\,S}(\mathcal{W}\setminus Y, \mathbb{C}) $ \ aura un r\'esidu nul  d'apr\`es la suite exacte de Leray sur \ $\mathcal{W}$.  On en d\'eduit l'existence de \ $ \omega \in \mathcal{C}^{\infty}_{mod\,S}(\mathcal{W})^q $ \ v\'erifiant \ $ d\omega = 0 $ \ et dont la restriction \`a \ $Z\setminus Y$ \ induit la classe \ $ [w - \frac{df}{f} \wedge \eta ] \in H^q_{mod\,S}(\mathcal{W}\setminus Y, \mathbb{C}) $. On en alors d\'eduit l'existence de \ $\alpha$ \ donnant \ $(@)$. $\hfill \blacksquare$

\bigskip

\section{ Construction de \ $ \widetilde{var} : H^n_{c\,\cap\,S}(0) \rightarrow H^n_c(0) .$}

\subsection{}

 Consid\'erons donc, pour \ $ k \geq 1 $ \ donn\'e, \ $w\in\Gamma_{c\,\cap\,S}(Y, \mathcal{E}^n(k)) \cap Ker\,\delta $. Notons par \ $X'$ \ la trace sur \ $X$ \ d'une boule centr\'ee \`a l'origine et de rayon \ $\varepsilon' < \varepsilon$ \ mais assez proche de \ $\varepsilon$ \ pour avoir \ $ Supp(\,w) \cap\, S \subset Y' $ \ o\`u \ $ Y' : = Y \cap X' $. Notons alors par \ $ Z$ \ un voisinage ouvert de \ $ Y \setminus \bar{Y}' $. Sur \ $Z$,  le morphisme de complexe
$$( \mathcal{F}^{\bullet}_{mod\,S}(k), \delta^{\bullet}) \rightarrow (\mathcal{E}^{\bullet}_{mod\,S}(k), \delta^{\bullet}) $$
est un quasi isomorphisme\footnote{On a not\'e ici  par \ $\mathcal{G}_{mod\,S}$ \ le sous faisceau du faisceau \ $\mathcal{G}$ \ des sections qui sont nulles au voisinage de \ $S$. On a donc, par d\'efinition la suite exacte courte :
$$ 0 \rightarrow \mathcal{G}_{mod\,S} \rightarrow \mathcal{G} \rightarrow j_! j^*(\mathcal{G}) \rightarrow 0  $$
o\`u \ $j : X \setminus S \hookrightarrow X $ \ d\'esigne l'inclusion.}. Comme la monodromie agit comme l'identit\'e\footnote{L'op\'erateur \ $_k\mathcal{N}$ \ est nul sur \ $h^0(k)$ \ et, en dehors de \ $S$,  sur \ $h^1(k)$. Attention, pas sur \ $S$ \ pour \ $n=2$ !} sur les faisceaux de cohomologies du complexe  \ $( \mathcal{F}^{\bullet}_{mod\,S}(k), \delta^{\bullet}) $ \ on peut trouver \ $ u\in\Gamma_{mod\,S}(Z, \mathcal{E}^{n-1}(k))$ \ v\'erifiant 
$$ _k\mathcal{N}(w) = \delta u \quad {\rm sur} \quad Z .$$
Soit \ $\rho\in\mathcal{C}^{\infty}_{c/f}(X)$ \ valant identiquement  1  au voisinage de \ $\bar{X}'$. Alors posons  \ $\xi : = (1-\rho).u$ \ et \ $v : = _k\mathcal{N}(w) - \delta \xi $. On a alors \ $v \in \Gamma_{c/f}(Y, \mathcal{E}^n(k)) \cap Ker\, \delta $. Cela donne explicitement 
$$ w_{j-1} = v_j + d\xi_j - \frac{df}{f}\wedge \xi_{j-1} \qquad \forall j\in [1,k]  $$
avec la convention \ $ w_0 = 0 = \xi_0 $. On aura en particulier, puisque \ $\delta w = 0 $
$$ d(w_k + \frac{df}{f}\wedge \xi_k) = \frac{df}{f}\wedge v_k .$$
Choisissons maintenant une boule centr\'ee \`a l'origine de rayon \ $ \varepsilon'' \in ]\varepsilon', \varepsilon[$ \ et assez proche de \ $\varepsilon$ \ pour que la trace \ $X''$ \ de cette boule sur \ $X$ \ contienne \ $Supp\, (v)$. Posons alors \ $\tilde{Y} : = Y \setminus\bar{Y}'' $. Alors la forme semi-m\'eromorphe  \ $w_k + \frac{df}{f}\wedge \xi_k$ \ de degr\'e \ $n$ \ est \ $d-$ferm\'ee sur un voisinage ouvert assez petit \ $\tilde{Z}$ \ de \ $\tilde{Y}$ \ et  a un support qui ne rencontre pas \ $S \cap \tilde{Y} $. On peut donc lui appliquer le th\'eor\`eme  (5.1)  et l'\'ecrire sur \ $\tilde{Z}$, quitte \`a localiser autour de \ $f^{-1}(0)$ :
\begin{equation*}
 w_k + \frac{df}{f}\wedge \xi_k =  \frac{df}{f}\wedge\eta + \omega + d\alpha   \tag{@@}
 \end{equation*}
o\`u les formes \ $\mathcal{C}^{\infty}$ \ $\eta$ \ et \ $\omega$ \ sont  \ $\mathcal{C}^{\infty}$ \ sur un voisinage ouvert \ $\mathcal{V}(\tilde{Y})$ \ de \ $\tilde{Y}$ \ dans \ $\tilde{Z}$, de degr\'es respectifs \ $n-1$ \ et \ $n$, v\'erifient :
\begin{itemize}
\item{i)} \ $d\eta = 0 , \  d\omega = 0 $
\item{ii)}\ $ Supp\, \eta \cap S = \emptyset , \  Supp\, \omega \cap S = \emptyset .$
\end{itemize}
et o\`u la \ $(n-1)-$forme semi-m\'eromorphe \ $\alpha$ \ sur \ $\mathcal{V}(\tilde{Y})$, est \`a p\^oles dans \ $\tilde{Y}$ \ et \`a support \ $mod\,S$.
Soit maintenant \ $\sigma\in\mathcal{C}^{\infty}_{c/f}(X)$ \ valant identiquement  1  au voisinage de \ $\bar{X}'' $. Alors la forme semi-m\'eromorphe de degr\'e \ $n+1$ \ \`a poles dans \ $\tilde{Y}$ \ et \`a support \ $f-$propre dans \ $\tilde{Z}$
$$  W : = d\sigma \wedge (w_k +  \frac{df}{f}\wedge \xi_k) $$
est \'egale sur \ $\tilde{Z}$, grace \`a (@@), \`a 
$$W = -\frac{df}{f}\wedge(d\sigma\wedge\eta) + d\sigma\wedge \omega - d(d\sigma\wedge\alpha) $$
o\`u les formes \  $d\sigma\wedge\eta$ \ et \ $d\sigma \wedge\omega$ \ sont  \ $\mathcal{C}^{\infty}$ \ sur un voisinage ouvert \ $\mathcal{V}(\tilde{Y})$ \ de \ $\tilde{Y}$ \ dans \ $\tilde{Z}$, de degr\'es respectifs \ $n$ \ et \ $n+1$, \`a supports \ $f-$propres dans  \ $\mathcal{V}(\tilde{Y})$ \ v\'erifient "\`a fortiori"
\begin{itemize}
\item{i)} \ $d(d\sigma\wedge\eta) = 0 , \ d(d\sigma\wedge\omega) = 0 $
\item{ii)}\ $ Supp(d\sigma\wedge\eta) \cap S = \emptyset , \  Supp(d\sigma\wedge\omega) \cap S = \emptyset .$
\end{itemize}
et o\`u la \ $n-$forme semi-m\'eromorphe \ $d\sigma\wedge\alpha$ \ sur \ $\mathcal{V}(\tilde{Y})$, est \`a p\^oles dans \ $\tilde{Y}$ \ et \`a support \ $f-$propre $\cap\,mod\,S$ \ \'egalement.

\noindent {\bf Nous d\'efinirons alors}  
$$2i\pi. \widetilde{var}([w]) : = r^n_c(k)(\tilde{v}) \quad {\rm avec} $$ 
$$ \tilde{v}_k : = v_k + d\sigma\wedge\eta, \quad {\rm et} \quad  \tilde{v}_j : = v_j \quad \forall j\in[1,k-1].$$
Il reste \'evidemment \`a v\'erifier que tout ceci d\'efinit bien une application lin\'eaire
$$ \widetilde{var} : H^n_{c\,\cap\,S}(0) \rightarrow H^n_c(0) $$
commutant aux monodromies.

\bigskip

\subsection{}

 Commen{\c c}ons par remarquer que si l'on change \ $w$ \ par son image par \ $j_{k,k+k'}$ \ le r\'esultat de notre construction ne change pas. 

\noindent Supposons maintenant que \ $ w = \delta (\beta)$ \ avec \ $ \beta \in \Gamma_{c\,\cap\,S}(Y, \mathcal{E}^{n-1}(k)) $. Alors on peut choisir \ $ u = \mathcal{N}( \beta)$ \ et on aura
$$ w_k + \frac{df}{f}\wedge \xi_k = d\beta_k - \rho \frac{df}{f}\wedge \beta_{k-1} .$$
Sur \ $\tilde{Z}$ \ on aura donc \ $ w_k + \frac{df}{f}\wedge \xi_k = d\beta_k $. Ceci montre que l'on peut prendre \ $ \eta = 0 , \omega = 0 \ {\rm et} \  \alpha = \beta_k $. Comme on a \ $ v = \delta(\rho.\mathcal{N}(\beta)) $ \ qui induit la classe nulle dans \ $H^n_c(0)$, ceci montre bien que la modification de \ $w$ \ par un cobord ne change pas le r\'esultat de notre construction.

\noindent On en d\'eduit imm\'ediatement que si \ $ w = \tau_k (\gamma) $ \ avec \ $ \gamma \in \Gamma_{c\,\cap\,S}(Y, Ker \delta^{n-1}(k)) $, alors on trouve \ $0$ \ dans \ $H^n_c(0)$. En effet, on peut remplacer \ $w_k$ \ par \ $ j_{k,k+k'}(w_k) $ \ d'apr\`es ce qui pr\'ec\`ede, et alors on aura d'apr\`es l'\'egalit\'e (@) du (4.2)
$$ j_{k,k+k'}(\tau_k(\gamma)) = \tau_{k+k'}(_{k+k'}\mathcal{N}^{k'}(j_{k,k+k'}(\gamma)))  $$
ce qui est dans \ $ \delta \big(\Gamma_{c\,\cap\, S}(Y, \mathcal{E}^{n-1}(k+k')\big)$ \ pour  \ $ k' \geq k $.

\noindent Les changements de choix des fonctions \ $\rho$ \ et \ $\sigma$ \ sont laiss\'es en exercice au lecteur.

\noindent Le changement du choix de \ $\eta $ \ consiste, d'apr\`es le th\'eor\`eme (5.1)  \`a remplacer \ $\eta$ \ par \ $ \eta + d\theta $ \ o\`u \ $ \theta \in \mathcal{C}^{\infty}_{mod\,S}(\tilde{Z})^{n-2} $. Mais on ajoute alors \`a \ $ \tilde{v}$ \ le cobord 
$$ \delta (j_{1,k}(d\sigma \wedge \theta )= j_{1,k}(-d\sigma \wedge d\theta) $$
ce qui ne change pas l'image dans \ $H^n_c(0)$.

\noindent Il nous reste \`a voir l'action de la monodromie. Si on part de \ $_k\mathcal{N}(w) $ \ on peut prendre \ $_k\mathcal{N}(u)$ \ et donc \ $_k\mathcal{N}(\xi)$ \ dans notre construction. On doit  alors prendre l'image de \ $_k\mathcal{N}(\tilde{v}) =\  _k\mathcal{N}(v) $ \ dans \ $H^n_c(0)$ \ ce qui donne bien le r\'esultat attendu. En effet, l'\'egalit\'e
$$ w_{k-1} + \frac{df}{f} \wedge \xi_{k-1} = v_k + d\xi_k  $$
permet de voir que sur \ $\tilde{Z}$, on a \ $ w_{k-1} + \frac{df}{f} \wedge \xi_{k-1} = d\xi_k $. On peut donc choisir \ $ \eta = 0 , \omega = 0 \ {\rm et} \ \alpha = \xi_k $. Ceci termine nos v\'erifications sur cette construction.

\bigskip

\newpage

\section{ La forme hermitienne canonique \\ $ \mathcal{H} : H^n_{c\,\cap\,S}(0)\times H^n(0) \rightarrow \mathbb{C}.$}

\bigskip

\noindent Nous adopterons ici le point de vue de [B.90]  qui est repris dans [B.97].

\bigskip

\subsection{Th\'eor\`eme-D\'efinition.}

\noindent \textit{Sous les hypoth\`eses standards d\'efinissons la forme hermitienne canonique 
$$ \mathcal{H} : H^n_{c\,\cap\,S}(0) \times H^n(0) \rightarrow \mathbb{C} $$
de la fa{\c c}on suivante :  pour \ $ e \in H^n_{c\,\cap\,S}(0) $ \ et \ $e' \in H^n(0)$ \ repr\'esent\'es par \ $w \in \Gamma_{c\,\cap\,S}(Y, \mathcal{E}^n(k)) \cap Ker\,\delta $ \ et \ $w' \in \Gamma(Y,  \mathcal{E}^n(k)) \cap Ker\,\delta $ \ v\'erifiant \ $ [w] = e $ \ et \ $ r^n(k)(w') = e' $ \ posons
$$ \mathcal{H}(e,e') = \frac{1}{(2i\pi)^{n+1}} P_2\big(\lambda = 0 , \int_X \vert f \vert^{2\lambda} \rho.\frac{df}{f} \wedge w_k \wedge \frac{\bar{df}}{\bar{f}} \wedge \bar{w}'_k \big) $$
o\`u \ $ \rho \in \mathcal{C}_c^{\infty}(X) $ \ vaut identiquement  1  au voisinage du compact  \ $ Supp\,w \,\cap\,S $.}

\noindent \textit{Notons par \ $ \mathcal{I} : H^n_c(0) \times H^n(0) \rightarrow \mathbb{C} $ \ la dualit\'e (hermitienne) de Poincar\'e sur la fibre de Milnor \ $F$ \ de \ $f$ \ \a l'origine, donn\'ee par
$$ \mathcal{I}(a,b) : = \frac{1}{(2i\pi)^n} \int_F a \wedge \bar{b} .$$
Alors on  a \ $ \mathcal{H}(e,e') = \mathcal{I}(\widetilde{var}(e), e') $ \ o\`u l'application
$$ \widetilde{var} : H^n_{c\,\cap\,S}(0) \rightarrow H^n_c(0)  $$
a \'et\'e construite plus haut.}

\bigskip

\subsection{Remarque.}

 Une cons\'equence imm\'ediate du th\'eor\`eme (7.1)  sera l'\'equivalence entre les deux propri\'et\'es suivantes  :
\begin{itemize}
\item{} La variation est bijective (ou bien \ $\widetilde{var}$ \ est bijective).
\item{}  La forme hermitienne canonique est non d\'eg\'en\'er\'ee.
\end{itemize}
\noindent Nous prouverons la premi\`ere de ces assertions pour tout  \ $n \geq 3$ \ \`a la fin  du paragraphe 8, ce qui donnera donc \'egalement  la non d\'eg\'en\'erescence de \ $\mathcal{H}$.

\noindent Pour \ $ n = 2$ \ la forme hermitienne canonique \ $ \mathcal{H}$ \ passe au quotient par le sous-espace (T-invariant) \ $j(H^1(S^*, \mathbb{C})$ \ de \ $H^2_{c\,\cap\,S}(0)$ \ et induit une dualit\'e hermitienne entre \  $H^2_{c\,\cap\,S}(0) \big/ j(H^1(S^*, \mathbb{C})$ \ et \ $H^2_{c\,\cap\,S}(0)$. Donc, \`a nouveau, on a un r\'esultat analogue au cas \ $n \geq 3$ \ quitte \`a remplacer \ $H^2_{c\,\cap\,S}(0)$ \ par son quotient par  \ $j(H^1(S^*, \mathbb{C})$.

\bigskip

\noindent La d\'emonstration du th\'eor\`eme (7.1) occupera le reste de ce paragraphe 7.

\subsection{}

Commen{\c c}ons par v\'erifier que le nombre \ $\mathcal{H}(e,e')$ \ est bien d\'efini, c'est \`a dire est ind\'ependant des divers choix effectu\'es.
\begin{itemize}
\item{Ind\'ependance du choix de \ $\rho$} \\
Soit donc \ $\rho'\in \mathcal{C}_c^{\infty}(X) $ \ valant aussi  identiquement  1  au voisinage du compact  \ $ Supp\,w \,\cap\,S $. Alors on a \ $Supp\,(\rho-\rho')\,\cap\,S = \emptyset $. Donc le prolongement m\'eromophe de 
$$\int_X \vert f \vert^{2\lambda}( \rho- \rho').\frac{df}{f} \wedge w_k \wedge \frac{\bar{df}}{\bar{f}} \wedge \bar{w}'_k $$
n'a, au pire, qu'un p\^ole simple en \ $\lambda = 0 $, ce qui donne l'ind\'ependance d\'esir\'ee.
\item{Ind\'ependance du choix de \ $w'$ \ repr\'esentant \ $e'$}\\
Il s'agit en fait de montrer que si on a \ $ w' = \delta v' $ \ avec \ $ v'\in \Gamma(Y, \mathcal{E}^{n-1}(k))$ \ alors on a 
$$P_2\big(\lambda = 0 , \int_X \vert f \vert^{2\lambda} \rho.\frac{df}{f} \wedge w_k \wedge \frac{\bar{df}}{\bar{f}} \wedge \bar{\delta v}'_{k} \big) = 0  .$$
Comme on a \ $ \frac{df}{f} \wedge w'_k = - d(\frac{df}{f}\wedge v'_k )$ \ ceci va r\'esulter de la formule de Stokes et du fait que \ $Supp\,d\rho\,\cap\,Supp\,w\,\cap \, S = \emptyset $, le prolongement m\'eromorphe de \ $\vert f \vert^{2\lambda} $ \ sur \ $ X \setminus S$ \ n'ayant, au pire, que des p\^oles simples aux entiers n\'egatifs.

\item{Ind\'ependance du choix de \ $w$ \ repr\'esentant \ $e$}\\
Comme ci-dessus, il s'agit en fait de montrer que si on a \ $ w= \delta v$ \ avec \ $ v\in \Gamma_{c\,\cap\, S}(Y, \mathcal{E}^{n-1}(k))$ \ alors on a 
$$P_2\big(\lambda = 0 , \int_X \vert f \vert^{2\lambda} \rho.\frac{df}{f} \wedge\delta v_k \wedge \frac{\bar{df}}{\bar{f}} \wedge \bar{ w}'_{k} \big) = 0  .$$
La preuve est analogue au cas pr\'ec\'edent, en choisissant la fonction \ $\rho$ \ de fa{\c c}on qu'elle soit identiquement \'egale \`a  1  sur un  voisinage ouvert de  \ $Supp\,v\,\cap\,S $. Ceci est possible grace \`a la condition de support sur \ $v$ \ et  \`a l'ind\'ependance du choix de \ $\rho$ \ prouv\'ee ci-dessus.

\end{itemize}

\subsection{}

Montrons maintenant l'invariance de \ $\mathcal{H}$ \ par la monodromie, c'est \`a dire que l'on a
$$ \mathcal{H}(\mathcal{N}(e), e') = \mathcal{H}(e,\mathcal{N}(e')) \quad \forall e\in H^n_{c\,\cap\,S}(0), \quad \forall e'\in H^n(0) $$
o\`u \ $-2i\pi.\mathcal{N}$ \ d\'esigne le logarithme (nilpotent) de la monodromie agissant sur \ $H^n_{c\,\cap\,S}(0)$ \ ou bien sur \ $ H^n(0) $. Compte tenu des ind\'ependances de choix prouv\'es ci-dessus, ceci revient \`a \'etablir la formule :
\begin{align*}
& P_2\big(\lambda = 0 , \int_X \vert f \vert^{2\lambda} \rho.\frac{df}{f} \wedge w_{k-1} \wedge \frac{\bar{df}}{\bar{f}} \wedge \bar{w}'_k \big) = \\
& P_2\big(\lambda = 0 , \int_X \vert f \vert^{2\lambda} \rho.\frac{df}{f} \wedge w_k \wedge \frac{\bar{df}}{\bar{f}} \wedge \bar{w}'_{k-1} \big) 
\end{align*}
puisque \ $_k\mathcal{N}(w)$ \ et respectivement \ $_k\mathcal{N}(w')$ \ induisent \ $\mathcal{N}(e)$ \ et
 \ $\mathcal{N}(e')$ \ dans \ $H^n_{c\,\cap\,S}(0)$ \ et \ $H^n(0)$.

\noindent Mais pour \ $\Re(\lambda) \gg 0 $ \ on a 
\begin{align*}
d\big(\vert f \vert^{2\lambda}\rho.(\frac{df}{f} + \frac{\bar{df}}{\bar{f}})\wedge w_k \wedge \bar{w}'_k \big) = \\
\vert f \vert^{2\lambda}.d\rho \wedge(\frac{df}{f} + \frac{\bar{df}}{\bar{f}})\wedge w_k \wedge \bar{w}'_k \  + \\
-\vert f \vert^{2\lambda}\rho.(\frac{df}{f} + \frac{\bar{df}}{\bar{f}})\wedge dw_k \wedge \bar{w}'_k  \ + \\
(-1)^{n+1} \vert f \vert^{2\lambda}\rho.(\frac{df}{f} + \frac{\bar{df}}{\bar{f}})\wedge w_k \wedge d\bar{w}'_k .
\end{align*} 
Comme on a \ $Supp\,d\rho\,\cap\,Supp\,w\,\cap \, S = \emptyset $ \ le prolongement m\'eromorphe du premier terme du membre de droite de l'\'egalit\'e ci-dessus n'aura, au pire, qu'un p\^ole simple en \ $\lambda = 0$. Par ailleurs le fait que \ $w$ \ et \ $w'$ \ soit \ $\delta-$ferm\'ees donne
$$ dw_k = \frac{df}{f}\wedge w_{k-1} \quad {\rm et} \quad dw'_k = \frac{df}{f}\wedge w'_{k-1} .$$
La formule de Stokes et le prolongement analytique donnent  alors la formule d\'esir\'ee.

\bigskip

\subsection{}

Il nous reste \`a prouver la formule 
 $$ \mathcal{H}(e,e') = \mathcal{I}(\widetilde{var}(e), e') $$
pour achever la preuve du th\'eor\`eme (7.1).

\noindent Reprenons les notations utilis\'ees en (6.1) lors de la construction de \ $\widetilde{var}(e)$.
Pour \ $\Re{\lambda} \gg  0$ \ on aura
\begin{align*}
 d\big((\vert f \vert^{2\lambda}\sigma.(w_k +  \frac{df}{f}\wedge\xi_k) \wedge \frac{\bar{df}}{\bar{f}} \wedge \bar{w}'_k \big) = & \  \lambda.\vert f \vert^{2\lambda}\sigma.\frac{df}{f} \wedge w_k \wedge \frac{\bar{df}}{\bar{f}} \wedge \bar{w}'_k  \\
&  + \vert f \vert^{2\lambda} d\sigma \wedge (w_k +  \frac{df}{f}\wedge\xi_k) \wedge \frac{\bar{df}}{\bar{f}} \wedge \bar{w}'_k \\
 & +  \vert f \vert^{2\lambda} \frac{df}{f} \wedge v_k \wedge \frac{\bar{df}}{\bar{f}} \wedge \bar{w}'_k  
\end{align*}
puisque \ $ d(w_k +  \frac{df}{f}\wedge\xi_k) = \frac{df}{f} \wedge v_k $ \ et \ $\sigma.v_k = v_k $.

\noindent D'apr\`es la formule de Stokes, l'int\'egrale du membre de gauche n'aura aucun p\^ole, et, par prolongement analytique, la nullit\'e du r\'esidu en \ $\lambda = 0 $ \ du membre de droite donnera la relation
\begin{align*}
- P_2\big(\lambda = 0, \int_X \vert f \vert^{2\lambda}\sigma.\frac{df}{f} \wedge w_k \wedge \frac{\bar{df}}{\bar{f}} \wedge \bar{w}'_k\big) = & \  Res\big(\lambda = 0, \int_X  \vert f \vert^{2\lambda}W \wedge \frac{\bar{df}}{\bar{f}} \wedge \bar{w}'_k\big) \\
+ &Res\big(\lambda = 0, \int_X  \vert f \vert^{2\lambda}\frac{df}{f} \wedge v_k \wedge \frac{\bar{df}}{\bar{f}} \wedge \bar{w}'_k \big) .
\end{align*}
Mais on a 
$$W = -\frac{df}{f}\wedge(d\sigma\wedge\eta) + d\sigma\wedge \omega - d(d\sigma\wedge\alpha) .$$
Montrons qu'alors
\begin{align*}
Res\big(\lambda = 0, \int_X  \vert f \vert^{2\lambda}d\sigma\wedge \omega \wedge \frac{\bar{df}}{\bar{f}} \wedge \bar{w}'_k\big) = 0 \\
Res\big(\lambda = 0, \int_X  \vert f \vert^{2\lambda}d(d\sigma\wedge\alpha) \wedge \frac{\bar{df}}{\bar{f}} \wedge \bar{w}'_k\big)  = 0 .
\end{align*}
Pour le premier r\'esidu c'est cons\'equence du fait que la forme \ $\omega$ \ est \ $\mathcal{C}^{\infty}$ \ et donc qu'il n'y a pas de puissance n\'egative de \ $f$ \ dans cette int\'egrale. On conclut en utilisant que les racines du polyn\^ome de Bernstein de \ $f$ \ sont strictement n\'egatives ([K.76]).

\noindent Pour montrer que le second r\'esidu est nul, on utilise \`a nouveau la formule de Stokes et le prolongement analytique pour obtenir
$$Res\big(\lambda = 0, \int_X  \vert f \vert^{2\lambda}d(d\sigma\wedge\alpha) \wedge \frac{\bar{df}}{\bar{f}} \wedge \bar{w}'_k\big)  =  P_2\big(\lambda = 0, \int_X  \vert f \vert^{2\lambda}\sigma. \frac{df}{f}\wedge d\alpha\wedge \frac{\bar{df}}{\bar{f}} \wedge \bar{w}'_k\big). $$
Mais comme le support de \ $d\alpha$ \ ne rencontre pas \ $S$ \ le membre de droite est nul.
Il nous reste donc finalement l'\'egalit\'e :
$$(2i\pi)^{n+1}. \mathcal{H}(e,e') = - Res\big(\lambda = 0, \int_X  \vert f \vert^{2\lambda}\frac{df}{f} \wedge (v_k + d\sigma \wedge \eta)\wedge \frac{\bar{df}}{\bar{f}} \wedge \bar{w}'_k \big) .$$
Mais la fonction 
$$ s \rightarrow \big(\frac{i}{2\pi}\big)^n \int_{f=s} (v_k + d\sigma\wedge\eta)\wedge\bar{w}'_k $$
est un polyn\^ome de degr\'e au plus \ $k-1$ \ en \ $Log (\vert\frac{s}{s_0}\vert^2) $ \ dont le terme constant est \ $ \mathcal{I}(\widetilde{var}(e),e')$. 

\noindent Par transformation de Mellin on obtient la formule d\'esir\'ee. \ $\hfill \blacksquare$

\bigskip

\section{Injectivit\'e de la variation.}

\subsection{}

Pour montrer que la variation est injective il nous suffit, d'apr\`es le th\'eor\`eme (7.1) de montrer que pour toute classe non nulle  \ $ e \in H^n_{c\,\cap\,S}(0)$ \ il existe \ $e'\in H^n(0)$ \ telle que l'on ait \ $\mathcal{H}(e,e') \not= 0 .$ En fait, comme la variation commute \`a la monodromie, on peut se contenter de montrer cela pour une classe \ $e$ \ qui est invariante par la monodromie.

\noindent Soit donc \ $v \in\Gamma_{c\,\cap\,S}(Y, \mathcal{E}^n(k)) \cap Ker\,\delta $ \ d\'efinissant la classe \ $e$ \ consid\'er\'ee.
Comme on suppose que \ $T(e) = e$, ce qui \'equivaut \`a \ $\mathcal{N}(e) = 0 $, on pourra trouver, quitte \`a choisir \ $k$ \ assez grand, \ $ u\in\Gamma_{c\,\cap\,S}(Y, \mathcal{E}^{n-1}(k)) $ \ v\'erifiant \ $ _k\mathcal{N}(v) = \delta u $ \ ce qui se traduit par les \'egalit\'es :
$$ v_{j-1} = du_j - \frac{df}{f}\wedge u_{j-1}\quad \forall j\in[1,k] \qquad {\rm avec} \quad  v_0 = u_0 = 0 .$$
On en d\'eduit que l'on a \ $ d\hat{w} = 0 $ \ o\`u l'on a pos\'e \ $ \hat{w} = v_{k} + \frac{df}{f}\wedge u_{k-1} .$ De plus on a l' \'egalit\'e
$$ j_{1,k+1}(\hat{w}) = j_{k,k+1}(v) -\delta \tilde{u} $$
o\`u l'on d\'efinit \ $ \tilde{u}_{k+1} = 0$ \ et \ $ \tilde{u}_j = u_j \quad \forall j \in [1,k] .$ Ceci montre que la classe \ $e$ \ est repr\'esent\'ee par la forme semi-m\'eromorphe \ $d-$ferm\'ee \ $\hat{w}$ \ qui a son support dans \ $c\,\cap\,S$.

\noindent Remarquons par ailleurs que l'image dans \ $H^n(0)$ \ de la classe \ $e$ \  peut \'egalement \^etre repr\'esent\'ee par une forme m\'eromorphe $d-$ferm\'ee \ $ w \in \Gamma(Y, \Omega^n)$ \ gr\^ace au quasi-isomorphisme \ $(\Omega^{\bullet}(k), \delta^{\bullet}) \simeq (\mathcal{E}^{\bullet}(k), \delta^{\bullet}) $ \ pour \ $k = 1$. On aura alors l'existence de \ $\gamma \in \Gamma(Y, \mathcal{E}^{n-1}(1))$ \ v\'erifiant \ $ w = \tilde{w} + d\gamma $ \ au voisinage de \ $Y$.

\noindent La g\'en\'eralisation suivante du r\'esultat de contribution "sureffective" d\'emontr\'ee dans  [B.84 b] va nous permettre de conclure dans le cas o\`u l'image de la classe \ $e$ \ dans \ $H^n(0)$ \  est non nulle.
Nous traiterons le cas o\`u cette image est nulle \`a la fin du paragraphe 8 (voir (8.10) \`a (8.22)). 

\bigskip

\subsection{ Th\'eor\`eme.}

\noindent \textit{Soit \ $f : X \rightarrow D $ \ un repr\'esentant de Milnor d'un germe non constant de fonction holomorphe \`a l'origine de \ $\mathbb{C}^{n+1}$. Soit \ $ Y = f^{-1}(0)$ \ et soit  \ $S_1$ \ un  sous-ensemble analytique ferm\'e de \ $Y$ \ tel qu'en chaque point \ $y$ \ de \ $Y\setminus S_1$ \ la monodromie locale de \ $f$ \ en \ $y$ \ agissant sur la cohomologie (r\'eduite) de la fibre de Milnor de \ $f$ \ en \ $y$ \  ne pr\'esente pas la valeur propre 1.}

\noindent \textit{On suppose que l'on a \ $ H^n(X\setminus S_1, \mathbb{C}) = 0 $.}

\noindent \textit{Soit \ $w$ \ un n-forme m\'eromorphe d-ferm\'ee sur \ $X$ \ \`a p\^oles dans \ $Y$ \ induisant une classe non nulle dans \ $H^n(F, \mathbb{C}) $ \ o\`u \ $F$ \ d\'esigne la fibre de Milnor de \ $f$ \ en \ $0$. Supposons qu'il existe des formes semi-m\'eromorphes \ $\tilde{w}$ \ et \ $\gamma$ \ sur \ $X$ \ \`a p\^oles dans \ $Y$ \ de degr\'es respectifs \ $n$ \ et \ $n-1$ \ v\'erifiant les conditions suivantes
\begin{itemize}
\item{i)} \quad  Supp  $\tilde{w} \, \cap S_1 = K$ \ est un compact de \ $X$ ;
\item{ii)}  \quad On a  sur \ $X\setminus Y  \quad  w = \tilde{w} + d\gamma $.
\end{itemize}
Alors il existe \ $\omega \in \Gamma(X, \Omega^{n+1}) \ {\rm et} \quad  j \in [1, n] $ \ tels que le prolongement m\'eromorphe de l'int\'egrale
$$ \int_X \vert f \vert^{2\lambda}\bar{f}^{-j} \frac{df}{f}\wedge \tilde{w} \wedge \rho.\bar{\omega}  $$
ait en \ $\lambda = 0 $ \ un p\^ole d'ordre \ $\geq 2$ \ o\`u \ $ \rho \in \mathcal{C}^{\infty}_c(X) $ \ vaut identiquement \ $1$ \ au voisinage de \ $K$.}

\bigskip

\subsection{ Remarques.}
\begin{itemize}
\item{1)} La conclusion du th\'eor\`eme est, bien sur, ind\'ependante du choix de la fonction \ $\rho \in \mathcal{C}^{\infty}_c(X)$ \ valant identiquement \ $1$ \ au voisinage de \ $K$ \ puisque si \ $\sigma \in \mathcal{C}^{\infty}_c(X)$ \ vaut identiquement \ $0$ \ au voisinage de \ $K$ \ le support de \ $\sigma.\tilde{w} $ \ ne rencontre plus \ $S_1$ \ ce qui implique que le prolongement m\'eromorphe de \ $ \vert f \vert^{2\lambda}.\sigma,\tilde{w} $ \ ne pr\'esente que des p\^oles d'ordre \ $\leq 1 $ \ aux entiers.
\item{2)} Si la dimension du sous-ensemble analytique (ferm\'e) \ $S_1$ \ de Y est \ $\leq p$ \ avec \ $ 2p < n+1 $ \ alors la condition \ $H^n(X\setminus S_1, \mathbb{C}) = 0 $ \ est automatiquement r\'ealis\'ee\footnote{ On peut s'en convaincre en montrant par r\'ecurrence sur \ $p$ \ que pour un sous-ensemble analytique ferm\'e de dimension \ $p$ \ de \ $\mathbb{C}^{n+1}$ \  les faisceaux de cohomologie \`a support \ $\underline{H}^q_{S}(\mathbb{C}) $ \ sont nuls pour \ $ q < 2n -2p +1 $.}. 

\noindent Nous utiliserons essentiellement le cas \ $ p = 1$ \ et \ $n \geq 2 $ \ de ce th\'eor\`eme.
\end{itemize}

\bigskip

\subsection{}

 \textit{D\'emonstration.} Posons pour \ $ j \in \mathbb{N} $ :
\begin{equation*}
\mathcal{T}^{n,0}_j : =  \ Res(\lambda = 0, \int_X \vert f \vert^{2\lambda} \bar{f}^{-j} w \wedge \Box) 
\  + \  P_2(\lambda = 0, \int_X \vert f \vert^{2\lambda} \bar{f}^{-j} \frac{df}{f}\wedge \gamma^{n-1,0}\wedge \Box )
\end{equation*}
et montrons d\'eja que l'on a \ $ d'\mathcal{T}^{n,0}_j = 0 $ \ sur \ $ V : = X \setminus K $ \ pour tout \ $j \in \mathbb{N}$. On a
\begin{equation*}
 d'  Res(\lambda = 0, \int_X \vert f \vert^{2\lambda} \bar{f}^{-j} w \wedge \Box )\ = \
 P_2(\lambda = 0, \int_X \vert f \vert^{2\lambda} \bar{f}^{-j}\frac{df}{f}\wedge w \wedge \Box )
\end{equation*}
puisque \ $ dw = d'w = 0 $. Par ailleurs
\begin{equation*}
 d'P_2(\lambda = 0, \int_X \vert f \vert^{2\lambda} \bar{f}^{-j} \frac{df}{f}\wedge \gamma^{n-1,0}\wedge \Box ) \ = \
 - P_2(\lambda = 0, \int_X \vert f \vert^{2\lambda} \bar{f}^{-j} \frac{df}{f}\wedge d' \gamma^{n-1,0}\wedge \Box )
\end{equation*}
ce qui donne notre assertion puisque les p\^oles aux  entiers n\'egatifs du prolongement m\'eromorphe de \ $\vert f \vert^{2\lambda} $ \ sont simples au plus le long de \ $Y\setminus S_1$ \ et que l'on a \ $ w - d'\gamma^{n-1,0} = \tilde{w}^{n,0} = 0 $ \ au voisinage de \ $ V \cap S_1 $.

\smallskip

\noindent Supposons maintenant que pour  \ $ j = n $ \ le courant \ $\mathcal{T}^{n,0}_n$ \  soit \ $d'-$exact sur \ $V$. C'est \`a dire qu'il existe un courant \ $ U^{n-1,0}_n $ \ de type \ $(n-1,0)$ \ sur \ $V$ \ v\'erifiant
$$  \mathcal{T}^{n,0}_n = d' U^{n-1,0}_n \quad {\rm sur} \quad V .$$
Alors, comme on a 
$$ \bar{f}. \mathcal{T}^{n,0}_{j+1} =  \mathcal{T}^{n,0}_j \quad \forall j \in \mathbb{N} $$
on obtient, en posant \ $ U^{n-1,0}_j = \bar{f}^{n-j}.U^{n-1,0}_n \quad \forall j\in [1,n] $ 
$$ \mathcal{T}^{n,0}_j = d' U^{n-1,0}_j \quad {\rm sur} \quad V  \quad \forall j\in [1,n].$$
Montrons qu'alors les courants sur \ $V$ \ d\'efinis pour \ $ j \in [1,n-1] $
\begin{align*}
 \mathcal{T}^{n-1,1}_j : = & \  d''U^{n-1,0}_j + j.d\bar{f}\wedge U^{n-1,0}_{j+1} \ + \\
 & - P_2(\lambda = 0, \int_X \vert f \vert^{2\lambda} \bar{f}^{-j} \frac{\bar{df}}{\bar{f}}\wedge \gamma^{n-1,0}\wedge \Box )\  + \\
 & +  P_2(\lambda = 0, \int_X \vert f \vert^{2\lambda} \bar{f}^{-j} \frac{df}{f}\wedge \gamma^{n-2,1}\wedge \Box )
\end{align*}
sont \ $d'-$ferm\'es sur \ $V$. Pour cela calculons \ $ d''\mathcal{T}^{n,0}_j + j.d\bar{f}\wedge \mathcal{T}^{n,0}_{j+1} $. On trouve
\begin{align*}
 d'' \mathcal{T}^{n,0}_j  +  j.d\bar{f}\wedge \mathcal{T}^{n,0}_{j+1} = & \  P_2(\lambda = 0, \int_X \vert f \vert^{2\lambda} \bar{f}^{-j}\frac{\bar{df}}{\bar{f}}\wedge w \wedge \Box )  + \\
  & +  P_3(\lambda = 0, \int_X \vert f \vert^{2\lambda} \bar{f}^{-j}\frac{\bar{df}}{\bar{f}}\wedge \frac{df}{f}\wedge \gamma^{n-1,0} \wedge \Box ) + \\
  & - P_2(\lambda = 0, \int_X \vert f \vert^{2\lambda} \bar{f}^{-j}\frac{df}{f}\wedge d'' \gamma^{n-1,0} \wedge \Box )
\end{align*}
Utilisons maintenant les relations \ $ w = d'\gamma^{n-1,0} $ \  et \ $ d'' \gamma^{n-1,0} + d' \gamma^{n-2,1} = 0 $ \ valables au voisinage de \ $ V \cap S_1 $. Elles donnent, puisque l'on a
\begin{align*}
& d'  P_2(\lambda = 0, \int_X \vert f \vert^{2\lambda} \bar{f}^{-j}\frac{\bar{df}}{\bar{f}}\wedge  \gamma^{n-1,0} \wedge \Box ) =  \\
& - P_2(\lambda = 0, \int_X \vert f \vert^{2\lambda} \bar{f}^{-j}\frac{\bar{df}}{\bar{f}}\wedge d' \gamma^{n-1,0} \wedge \Box )  + \\
& - P_3(\lambda = 0, \int_X \vert f \vert^{2\lambda} \bar{f}^{-j}\frac{\bar{df}}{\bar{f}}\wedge \frac{df}{f}\wedge \gamma^{n-1,0} \wedge \Box ) 
\end{align*}
ainsi que
\begin{equation*}
 d' P_2(\lambda = 0, \int_X \vert f \vert^{2\lambda} \bar{f}^{-j} \frac{df}{f}\wedge \gamma^{n-2,1}\wedge \Box )\  = \  - P_2(\lambda = 0, \int_X \vert f \vert^{2\lambda} \bar{f}^{-j} \frac{df}{f}\wedge d' \gamma^{n-2,1}\wedge \Box )
\end{equation*}
l'\'egalit\'e
\begin{align*}
 d'' \mathcal{T}^{n,0}_j  +  j.d\bar{f}\wedge \mathcal{T}^{n,0}_{j+1} = & - d' P_2(\lambda = 0, \int_X \vert f \vert^{2\lambda} \bar{f}^{-j}\frac{\bar{df}}{\bar{f}}\wedge  \gamma^{n-1,0} \wedge \Box )\  + \\
 & + d'P_2(\lambda = 0, \int_X \vert f \vert^{2\lambda} \bar{f}^{-j} \frac{df}{f}\wedge \gamma^{n-2,1}\wedge \Box )
 \end{align*}
 ce qui prouve bien la \ $d'-$fermeture sur \ $V$ \ des courants \ $\mathcal{T}^{n-1,1}_j $ \ pour \ $j \in [1, n-1]$.
 
 \noindent Comme on a l'annulation des groupes \ $ H^{p+1}(V, \Omega^q) $ \ pour tout \ $ p \in [0,n-2]$ \ et tout \ $q \in \mathbb{N}$, on peut trouver un courant \ $ U^{n-2,1}_{n-1} $ \ sur \ $V$ \ v\'erifiant \ $ d' U^{n-2,1}_{n-1}  = \mathcal{T}^{n-1,1}_{n-1} $. Posons \ $  U^{n-2,1}_{j} : = \bar{f}^{n-j-1}.U^{n-2,1}_{n-1}$ \ pour \ $j\in [1,n-1] $. Alors les hypoth\`eses du Lemme suivant sont v\'erifi\'ees pour \ $p = 0 , a = 1 , b = n-1 $.

\bigskip

\subsection{Lemme.}

 On suppose donn\'es sur \ $V$ \ des courants \ $U_j^{n-p-1,p}$ \ pour \ $j \in [a,b+1]$ \ et \ $U_j^{n-p-2,p+1} $ \ pour \ $j \in [a,b]$ \   o\`u les \ $U_j^{n-p-1,p}$ \ v\'erifient \ $ \bar{f}.U_{j+1}^{n-p-1,p} = U_j^{n-p-1,p} $ \ pour \ $j \in [a,b]$ \ et o\`u les courants \ $U_j^{n-p-2,p+1}$ \ v\'erifient \ $ \bar{f}.U_{j+1}^{n-p-2,p+1} = U_j^{n-p-2,p+1} $ \ pour \ $j \in [a,b-1]$.

\noindent Supposons \'egalement donn\'ee sur \ $V$ \ une forme semi-m\'eromorphe \ $\gamma$ \ de degr\'e \ $n-1$ \ \`a p\^oles dans \ $\lbrace f = 0 \rbrace$ \ telle que \ $d\gamma$ \ soit de type \ $(n,0)$.

\noindent Supposons que l'on ait sur \ $V$ \ les \'egalit\'es de courants suivantes pour \ $ j \in [a,b]$
\begin{align*}
\mathcal{T}^{n-p-1,p+1}_j : = &\  d''U_j^{n-p-1,p} + j.\bar{df}\wedge U_{j+1}^{n-p-1,p} + \\
 \quad & - P_2(\lambda = 0,\int_V \vert f \vert^{2\lambda} \bar{f}^{-j}\frac{\bar{df}}{\bar{f}}\wedge \gamma^{n-p-1,p}\wedge\Box )  + \\
 \quad  & \  P_2(\lambda = 0,\int_V \vert f \vert^{2\lambda} \bar{f}^{-j}\frac{df}{f}\wedge \gamma^{n-p-2,p+1}\wedge\Box ) \\
 \quad  =  & \  d'U_j^{n-p-2,p+1} 
  \end{align*}
   Alors les courants
  \begin{align*}
 \mathcal{T}^{n-p+2,p+2}_j : = & \  d''U_j^{n-p-2,p+1} + j.\bar{df}\wedge U_{j+1}^{n-p-2,p+1} + \\
\quad &  - P_2(\lambda = 0,\int_V \vert f \vert^{2\lambda} \bar{f}^{-j}\frac{\bar{df}}{\bar{f}}\wedge \gamma^{n-p-2,p+1}\wedge\Box )  + \\
 \quad & +   P_2(\lambda = 0,\int_V \vert f \vert^{2\lambda} \bar{f}^{-j}\frac{df}{f}\wedge \gamma^{n-p-3,p+2}\wedge\Box )
 \end{align*} 
 sont \ $d'-$ferm\'es sur \ $V$ \ pour \ $ j \in [a,b-1] $.
 
 \bigskip
 
 \subsection{Preuve du lemme.}
 
On a pour \ $j \in [a,b]$ :
 \begin{align*}
 d'\big(d'' U_j^{n-p-2,p+1} \big) & = - d'' \big(d' U_j^{n-p-2,p+1} \big) \\
 \quad & = - d'' \big(j.\bar{df}\wedge U_{j+1}^{n-p-1,p}\big) \ + \\
 \quad & + d''\big(P_2(\lambda = 0,\int_V \vert f \vert^{2\lambda} \bar{f}^{-j}\frac{\bar{df}}{\bar{f}}\wedge \gamma^{n-p-1,p}\wedge\Box )\big) \\
 \quad & - d''\big( P_2(\lambda = 0,\int_V \vert f \vert^{2\lambda} \bar{f}^{-j}\frac{df}{f}\wedge \gamma^{n-p-2,p+1}\wedge\Box )  \big)
 \end{align*}
 et donc
 \begin{align*}
  d'\big(d'' U_j^{n-p-2,p+1} \big) & =  j.\bar{df}\wedge d'' U_{j+1}^{n-p-1,p}  \qquad  \tag{1}  \\
  \quad & - P_2(\lambda = 0,\int_V \vert f \vert^{2\lambda} \bar{f}^{-j}\frac{\bar{df}}{\bar{f}}\wedge d''\gamma^{n-p-1,p}\wedge\Box ) \qquad \tag{2} \\
  \quad  & - P_3(\lambda = 0,\int_V \vert f \vert^{2\lambda} \bar{f}^{-j}\frac{\bar{df}}{\bar{f}}\wedge\frac{df}{f}\wedge\gamma^{n-p-2,p+1}\wedge\Box )  \quad \tag{3} \\
  \quad & j.P_2(\lambda = 0,\int_V \vert f \vert^{2\lambda} \bar{f}^{-j}\frac{\bar{df}}{\bar{f}}\wedge\frac{df}{f}\wedge \gamma^{n-p-2,p+1}\wedge\Box )  \quad \tag{4} \\
  \quad & + P_2(\lambda = 0,\int_V \vert f \vert^{2\lambda} \bar{f}^{-j}\frac{df}{f}\wedge d''\gamma^{n-p-2,p+1}\wedge\Box )  \quad \tag{5}
  \end{align*}
  puisque l'on a
 \begin{align*}
 & d'' \big(\int_V \vert f \vert^{2\lambda} \bar{f}^{-j}\frac{df}{f}\wedge\gamma^{n-p-2,p+1}\wedge\Box \big) =  \\
&\qquad \qquad \qquad (\lambda - j).\int_V \vert f \vert^{2\lambda} \bar{f}^{-j}\frac{\bar{df}}{\bar{f}}\wedge\frac{df}{f}\wedge \gamma^{n-p-2,p+1}\wedge\Box \quad + \\
&\qquad \qquad  \qquad - \int_V \vert f \vert^{2\lambda} \bar{f}^{-j}\frac{df}{f}\wedge d''\gamma^{n-p-2,p+1}\wedge\Box 
 \end{align*}
 
  Ensuite on a
  \begin{align*}
&  d'\big( j.\bar{df}\wedge U_{j+1}^{n-p-2,p+1}\big) =  -j.\bar{df}\wedge d''U_{j+1}^{n-p-1,p} \ + \tag{6} \\
 &  -j.\bar{df}\wedge  P_2(\lambda = 0,\int_V \vert f \vert^{2\lambda} \bar{f}^{-j-1}\frac{df}{f}\wedge \gamma^{n-p-2,p+1}\wedge\Box ) \quad \tag{7} ;
  \end{align*}
  puis
  \begin{align*}
 & - d' P_2(\lambda = 0,\int_V \vert f \vert^{2\lambda} \bar{f}^{-j}\frac{\bar{df}}{\bar{f}}\wedge \gamma^{n-p-2,p+1}\wedge\Box ) = \\
 & - P_3(\lambda = 0,\int_V \vert f \vert^{2\lambda} \bar{f}^{-j}\frac{df}{f}\wedge\frac{\bar{df}}{\bar{f}}\wedge\gamma^{n-p-2,p+1}\wedge\Box )  \quad \tag{8} \\
  & + P_2(\lambda = 0,\int_V \vert f \vert^{2\lambda} \bar{f}^{-j}\frac{\bar{df}}{\bar{f}}\wedge d'\gamma^{n-p-2,p+1}\wedge\Box ) \quad \tag{9}
  \end{align*}
  puisque 
  \begin{align*}
 & d' \int_V \vert f \vert^{2\lambda} \bar{f}^{-j}\frac{\bar{df}}{\bar{f}}\wedge \gamma^{n-p-2,p+1}\wedge\Box  = \\
 & \lambda.\int_V \vert f \vert^{2\lambda} \bar{f}^{-j}\frac{df}{f}\wedge\frac{\bar{df}}{\bar{f}}\wedge\gamma^{n-p-2,p+1}\wedge\Box \ + \\
 & - \int_V \vert f \vert^{2\lambda} \bar{f}^{-j}\frac{\bar{df}}{\bar{f}}\wedge d'\gamma^{n-p-2,p+1}\wedge\Box .
 \end{align*}
 Et enfin 
  \begin{align*}
&  d'\big(P_2(\lambda = 0,\int_V \vert f \vert^{2\lambda} \bar{f}^{-j}\frac{df}{f}\wedge \gamma^{n-p-3,p+2}\wedge\Box ) \big) = \\
& - P_2(\lambda = 0,\int_V \vert f \vert^{2\lambda} \bar{f}^{-j}\frac{df}{f}\wedge d' \gamma^{n-p-3,p+2}\wedge\Box ) \quad \tag{10}
\end{align*}
Maintenant les relations
\begin{align*}
 & d''\gamma^{n-p-1,p} + d'\gamma^{n-p-2,p+1} = 0  \\
 & d''\gamma^{n-p-2,p+1} + d'\gamma^{n-p-3,p+2} = 0
\end{align*}
qui r\'esultent de notre hypoth\`ese sur \ $\gamma$ \ ainsi que les \'egalit\'es
$$ (1) + (6) = 0 , (2)+ (9) = 0 , (3) + (8) = 0 , (4) + (7) = 0 , (5) + (10) = 0 $$
permettent de conclure.  \ $\hfill \blacksquare$

\bigskip

\subsection{Remarque.}

Pour peu que l'on ait \ $ H^{p+2}(V, \Omega^{n-p+2}) = 0 $ \ on peut trouver, sous l'hypoth\`ese du Lemme, un courant \ $ U^{n-p+1,p+2}_{b-1}$ \ sur \ $V$ \ v\'erifiant \ $ d'U^{n-p+1,p+2}_{b-1} = \mathcal{T}^{n-p+2,p+2}_{b-1} $ \ et poser \ $ U^{n-p+1,p+2}_{j} = \bar{f}^{b-j-1}. U^{n-p+1,p+2}_{b-1} $ \ pour \ $j\in [a,b-a -1] $ \ et se retrouver \`a nouveau dans l'hypoth\`ese du Lemme pour \ $ p +1, a \ {\rm et} \ b-1 $.

\bigskip

\subsection{Fin de la d\'emonstration du th\'eor\`eme (8.2).}

 Puisque l'on dispose de l'annulation de $ H^{p+1}(V, \Omega^q) $ \ pour tout \ $ p \in [0,n-2]$ \ et tout \ $q \in \mathbb{N}$, la remarque ci-dessus permet de continuer \`a appliquer  le Lemme \ $(n-1)-$fois (pour \ $p\in [0,n-2]$)\ et d'obtenir un courant $d'-$ferm\'e de type \ $(0,n)$
\begin{align*}
 \mathcal{T}^{0,n}_1 : = & d''U^{0,n-1}_1 + \bar{df}\wedge U^{0,n-1}_2 + \\
  \quad & - P_2(\lambda = 0,\int_V \vert f \vert^{2\lambda} \bar{f}^{-j}\frac{\bar{df}}{\bar{f}}\wedge \gamma^{0,n-1}\wedge\Box ) .  \\
  \end{align*}
  Donc \ $ \mathcal{T}^{0,n}_1$ \ est une forme antiholomorphe de degr\'e  \ $n$ \ sur \ $V$. De plus on a 
  $$ d'' \mathcal{T}^{0,n}_1 + \frac{\bar{df}}{\bar{f}} \wedge  \mathcal{T}^{0,n}_1 = 0 $$
  sur \ $V\setminus Y$ \ car \ $ d''\gamma^{0,n-1} = 0 $ \ montre que \ $ P_2(\lambda = 0,\int_V \vert f \vert^{2\lambda} \bar{f}^{-j}\frac{\bar{df}}{\bar{f}}\wedge \gamma^{0,n-1}\wedge\Box )$ \ est annul\'e par \ $\wedge \bar{df} $ \ et \ $d$. 
  \noindent Par Hartogs, elle se prolonge  en une forme antiholomorphe \ $\bar{\Omega_1} $ \ sur \ $X$ \ v\'erifiant \ $ d(f.\Omega_1) = 0 $ \ sur \ $X$.
  
  \noindent De plus on aura des courants \ $\mathcal{U}$ \ et \ $\mathcal{V}$ \ sur \ $V$ \ de degr\'e $n-1$ \ v\'erifiant  l'\'egalit\'e suivante \ sur \ $V\setminus S_1 = X \setminus S_1$\footnote{Rappelons que les supports des parties polaires d'ordre \ $\geq 2$ \ aux entiers du prolongement m\'eromorphes de \ $\vert f \vert^{2\lambda}$ \ sont contenus dans \ $S_1$.} :
  $$ \mathcal{T}^{n,0}_1 - \bar{\Omega}_1 =  d\mathcal{U} + \bar{df}\wedge \mathcal{V} . $$
  Remarquons d\'eja que la \ $n-$forme holomorphe $d-$ferm\'ee \ $\Omega_0 :  = f.\Omega_1$ \ est $d-$exacte sur \ $X$. Donc elle induit la classe nulle dans \ $H^n(F,\mathbb{C})$. Il en sera donc de m\^eme pour \ $ \Omega_1$ et pour  une $n-$forme holomorphe \ $A$ \ v\'erifiant \ $ dA = df \wedge \Omega_1 $, d'apr\`es le th\'eor\`eme de positivit\'e de Malgrange (voir [M.74] ou l'appendice de [B.84 b]).
  
  \noindent  Posons
  $$ T : = Pf(\lambda = 0 ,\int_X \vert f \vert^{2\lambda} w \wedge \Box ). $$
  Montrons que l'on a \ $ dT = \bar{df}\wedge \mathcal{T}^{n,0}_1 $ \ sur \ $X$. En effet, on a
  $$ d'T = Res(\lambda = 0 ,\int_X \vert f \vert^{2\lambda}\frac{df}{f}\wedge w \wedge \Box ) = 0 $$
  car l'absence de puissance de \ $\bar{f}$ \ en d\'enominateur montre que la fonction m\'eromorphe dont on prend le r\'esidu \`a l'origine n'a pas de p\^oles en ce point\footnote{On utilise \`a nouveau le fait que les racines du polyn\^ome de Bernstein de \ $f$ \ sont strictement n\'egatives ! Voir [K.76]}. De plus on a
  $$ d''T = Res(\lambda = 0 ,\int_X \vert f \vert^{2\lambda}\frac{\bar{df}}{\bar{f}}\wedge w \wedge \Box ) $$
  ce qui co{\"i}ncide bien avec \ $ \bar{df}\wedge \mathcal{T}^{n,0}_1 $ \ car on a
  \begin{align*}
&  P_2(\lambda = 0, \int_X \vert f \vert^{2\lambda} \frac{\bar{df}}{\bar{f}}\wedge \frac{df}{f}\wedge \gamma^{n-1,0}\wedge \Box )  =  \\
& \qquad \qquad d''Res(\lambda = 0, \int_X \vert f \vert^{2\lambda} \frac{df}{f}\wedge \gamma^{n-1,0}\wedge \Box ) \ + \\
 & \qquad  \qquad + Res(\lambda = 0, \int_X \vert f \vert^{2\lambda} \frac{df}{f}\wedge d'' \gamma^{n-1,0}\wedge \Box ) \\
\quad & \qquad \qquad  =  \ 0
 \end{align*}
 car \`a nouveau les residus \`a l'origine sont nuls pour la m\^eme raison que ci-dessus.
 
 \noindent On obtient donc l'\'egalit\'e suivante de courants sur \ $X \setminus S_1$
 $$ d\big(T - \bar{A} + \bar{df}\wedge \mathcal{U} \big) = 0 $$
 Comme on  a, par hypoth\`ese, \ $H^n(X\setminus S_1, \mathbb{C}) = 0 $ \ on en conclut \`a l'existence d'un courant \ $ \mathcal{W}$ \  de degr\'e $n-1$ \ sur \ $X \setminus S_1$ \ tel que
 $$ T = \bar{A} + \bar{df}\wedge \mathcal{U} + d\mathcal{W} .$$
 Alors le Lemme (D) de [B.84 a]  donne que \ $w$ \ induit la classe nulle dans \ $H^n(F, \mathbb{C})$ \ contredisant notre hypoth\`ese.
 
 \noindent Donc le courant \ $\mathcal{T}^{n,0}_n$ \ n'est pas \ $d'-$exact sur \ $V$. La dualit\'e entre \\ 
 $H^n(X\setminus K, \mathcal{O}) $ \ et \ $\Gamma(K,\Omega^{n+1})$\footnote{On peut choisir pour \ $K$ \ la trace sur \ $S_1$ \ d'une boule assez grosse.} et la densit\'e de l'image de la restriction \ $ \Gamma(X,\Omega^{n+1}) \rightarrow \Gamma(K, \Omega^{n+1})$ \ permettent  alors de trouver une forme holomorphe \ $\omega\in \Gamma(X,\Omega^{n+1})$ \ telle que l'on ait \ $< \mathcal{T}^{n,0}_n , d'\rho\wedge \bar{\omega} > \not= 0 $ \ ce qui donne
 \begin{align*}
&  \ Res(\lambda = 0, \int_X \vert f \vert^{2\lambda} \bar{f}^{-n} w \wedge d'\rho\wedge \bar{\omega} )  \ +\\
\qquad & P_2(\lambda = 0, \int_X \vert f \vert^{2\lambda} \bar{f}^{-n} \frac{df}{f}\wedge \gamma^{n-1,0}\wedge d'\rho\wedge \bar{\omega} ) \not= 0
\end{align*}
Mais on a
\begin{align*}
 d\big( \vert f \vert^{2\lambda} \bar{f}^{-n} .w \wedge \rho. \bar{\omega} ) = &  \lambda \vert f \vert^{2\lambda} \bar{f}^{-n}.\frac{df}{f}\wedge  w \wedge \rho. \bar{\omega} \  + \\
 \qquad \qquad & (-1)^n \vert f \vert^{2\lambda} \bar{f}^{-n} . w \wedge d'\rho\wedge \bar{\omega}
\end{align*}
ce qui donne :
\begin{align*}
&  \ Res(\lambda = 0, \int_X \vert f \vert^{2\lambda} \bar{f}^{-n} w \wedge d'\rho\wedge \bar{\omega} ) = \\
& \qquad \qquad (-1)^{n+1}  P_2(\lambda = 0, \int_X \vert f \vert^{2\lambda} \bar{f}^{-n} \frac{df}{f}\wedge w \wedge\rho. \bar{\omega} )
\end{align*}
Comme on a
\begin{align*}
&  P_2(\lambda = 0, \int_X \vert f \vert^{2\lambda} \bar{f}^{-n} \frac{df}{f}\wedge \gamma^{n-1,0}\wedge d'\rho\wedge \bar{\omega} ) = \\
&\qquad \qquad  (-1)^{n} P_2(\lambda = 0, \int_X \vert f \vert^{2\lambda} \bar{f}^{-n} \frac{df}{f}\wedge d' \gamma^{n-1,0}\wedge \rho\ . \bar{\omega} ) 
\end{align*}
on obtient finalement
$$ P_2(\lambda = 0, \int_X \vert f \vert^{2\lambda} \bar{f}^{-n} \frac{df}{f}\wedge(w - d' \gamma^{n-1,0})\wedge \rho\ . \bar{\omega} ) \not= 0 $$
c'est \`a dire, vus les types
$$P_2(\lambda = 0, \int_X \vert f \vert^{2\lambda} \bar{f}^{-n} \frac{df}{f}\wedge\tilde{w} \wedge \rho\ . \bar{\omega} ) \not= 0 .$$
Ceci ach\`eve la preuve du th\'eor\`eme (8.2). \ $\hfill \blacksquare$

\bigskip

\subsection{}

 Pour conclure \`a l'injectivit\'e de la variation  dans le cas o\`u l'image de la classe \ $ e \in H^n_{c\,\cap\,S}(0) $ \ dans \ $H^n(0)$ \ n'est pas nulle, il suffit de remarquer que si l'on \'ecrit \ $ \omega = \frac{df}{f}\wedge w' $ \ la d\'ecomposition dans le syst\`eme de Gauss-Manin (localis\'e) en degr\'e \ $n$ \  de \ $f$ \ de la forme m\'eromorphe \ $ \frac{w'}{f^n} $ \ va fournir une  classe \ $e' \in H^n(0)$ \ v\'erifiant \ $\mathcal{H}(e,e') \not= 0 $. En effet, les formes m\'eromorphes induisant des classes dans les sous espaces spectraux de la monodromie pour des valeurs propres diff\'erentes de \ $1$ \ ne donneront pas de p\^oles doubles aux entiers n\'egatifs, pas plus que les formes du type \ $ \frac{df}{f}\wedge du $ \ o\`u \ $u$ \ est m\'eromorphes de degr\'e \ $n-1$ \ dans le prolongement m\'eromorphe de
$$ \int_X \vert f \vert^{2\lambda} \rho \frac{df}{f}\wedge \tilde{w}\wedge \bar{\Box} .$$

\bigskip

\subsection{}

 Le cas plus d\'elicat est celui o\`u la classe \ $e$ \ consid\'er\'ee est dans le noyau de l'application \ $ can_{c\,\cap\,S} : H^n_{c\,\cap\,S}(0) \rightarrow H^n(0) .$ D'apr\`es le th\'eor\`eme 1 du (4.6) il existe \ $[\alpha] \in H^1_{\lbrace 0 \rbrace}(S, H^{n-1}(0)) $ \ v\'erifiant \ $ i([\alpha]) = e $. De plus, comme l'application \ $i$ \ est injective, on peut supposer que l'on a \ $ T([\alpha]) = [\alpha] $ \ dans \ $H^1_{\lbrace 0 \rbrace}(S, H^{n-1}(0)) .$

\noindent Les lemmes suivants vont nous fournir un repr\'esentant de la classe \ $[\alpha]$ \ consid\'er\'ee qui met en \'evidence la propri\'et\'e d'invariance par la monodromie.

\bigskip

\subsection{ Lemme.}

Supposons \ $n \geq 3$. Soit \ $ [\alpha] \in H^1_{\lbrace 0 \rbrace}(S, H^{n-1}(0))$ \ v\'erifiant \ $T(i([\alpha])) = i([\alpha])$ \ dans \ $ H^n_{c\,\cap\,S}(0) $. Alors, pour \ $k$ \ assez grand, il existe un repr\'esentant \\ \ $ \hat{\alpha}\in \Gamma(S^*, Ker\,\delta^{n-1}(k +1)) $ \  tel que l'on ait 
$$\mathcal{N}(\hat{\alpha}) \in \Gamma(Y,  Ker\,\delta^{n-1}(k +1)) $$
 et dont l'image dans \ $H^1_{\lbrace 0 \rbrace}(S, H^{n-1}(0)) $ \ est \ $[\alpha]$.

\bigskip

\noindent \textit{Preuve.} Pour \ $k \gg 1 $ \ on a un isomorphisme \ $ r^{n-1}(k) : h^{n-1}(k) \rightarrow H^{n-1}(0) .$ Puisque l'on a \ $H^1(S,H^{n-1}(0)) = 0 $ \  et l'annulation de \ $ H^1(S^*, \delta\mathcal{E}^{n-2}) $ \ prouv\'ee au lemme (4.7), on peut relever \ $[\alpha]$ \ en une section \ $\alpha' \in \Gamma(S^*, Ker\,\delta^{n-1}(k)) .$ Grace \`a  l'annulation de \ $ H^1(Y, \delta\mathcal{E}^{n-2}) $ \ prouv\'ee au lemme (4.7), l'hypoth\`ese \ $ T[\alpha] = [\alpha]$ \ dans \ $ H^1_{\lbrace 0 \rbrace}(S, H^{n-1}(0))$ \ permet  de trouver \ $ \beta \in \Gamma(Y,  Ker\,\delta^{n-1}(k))  $ \ et \ $ \Gamma \in \Gamma(S^*, \delta\mathcal{E}^{n-2}(k))$ \ v\'erifiant
$$ \mathcal{N} \alpha' = \beta + \Gamma \quad {\rm au \ voisinage \ de} \ S^* .$$
Ceci se d\'eduit d'une "chasse" facile sur le diagramme commutatif  exact  d'espaces vectoriels monodromiques :

$$\xymatrix{\quad & \quad &  0 & \quad \\
\quad & \quad &\quad H^1_{\lbrace 0 \rbrace}(S, H^{n-1}(0))\ar[u] & \quad \\
H^0(S^*, \delta\mathcal{E}^{n-2}(k)) \ar[r] & H^0(S^*, Ker\,\delta^{n-1}(k)) \ar[r] & H^0(S^*, H^{n-1}(0)) \ar[r] \ar[u]& 0 \\
H^0(Y, \delta\mathcal{E}^{n-2}(k)) \ar[r] \ar[u]& H^0(Y, Ker\,\delta^{n-1}(k)) \ar[r] \ar[u]& H^0(Y, H^{n-1}(0)) \ar[r] \ar[u] & 0 } $$
La suite exacte 
$$ 0 \rightarrow Ker\,\delta^{n-2}(k) \rightarrow \mathcal{E}^{n-2}(k) \rightarrow \delta\mathcal{E}^{n-2}(k) \rightarrow 0  $$
donne, pour \ $ n \geq 4$ \ la surjectivit\'e de la fl\`eche
$$ \delta : H^0(S^*, \mathcal{E}^{n-2}(k)) \rightarrow H^0(S^*, \delta\mathcal{E}^{n-2}(k)) $$
puisque, pour  \ $ n \geq 4$, on a \ $ Ker\,\delta^{n-2}(k) \simeq \delta\mathcal{E}^{n-3}(k)$ \ et que l'annulation du \ $H^1(S^*, \delta\mathcal{E}^{n-3}(k))$ \ est donn\'ee par (4.7). On peut donc trouver \\
 $ \gamma \in \Gamma(S^*,  \mathcal{E}^{n-2}(k)) $ \ v\'erifiant \ $ \Gamma = \delta\gamma $.
Pour \ $n = 3 $ \ on a \ $ H^1(S^*, Ker\,\delta^1(k)) \simeq H^1(S^*, h^1(k)) \not= 0 $. Mais l'application $$j_{k,k+1} : H^1(S^*, h^1(k)) \rightarrow H^1(S^*, h^1(k+1)) $$
est nulle (voir (1.6)) et quitte \`a changer \ $k$ \ en \ $k+1$ \ on peut lever l'obstruction et trouver \'egalement \ $ \gamma \in \Gamma(S^*,  \mathcal{E}^{n-2}(k+1)) $ \ v\'erifiant \ $ \Gamma = \delta\gamma $.

\noindent  D\'efinissons alors \ $ \hat{\alpha}\in \Gamma(S^*,  Ker\,\delta^{n-1}(k +1)) $ \  en posant 
$$ \hat{\alpha}_{k+1} = \alpha'_k + \frac{df}{f}\wedge \gamma_k \quad {\rm et} \quad \hat{\alpha}_j = \beta_j \quad \forall j\in [1,k] .$$
L'\'egalit\'e \ $\delta\hat{\alpha} = 0 $ \ est imm\'ediate et on a \ $ \mathcal{N}(\hat{\alpha}) = j_{k,k+1}(\beta)$ \ ce qui ach\`eve la preuve. $\hfill \blacksquare $

\bigskip

\subsection{Variante pour \ $n = 2 $.}

Si \ $\alpha \in \mathcal{K} $ \ v\'erifie \ $ \mathcal{N}([\alpha]) = 0 $ \ dans \ $ H^1_{\lbrace 0 \rbrace}(S, H^1(0))$, il existe \ $ k \gg 1 $ \ et \ $ \hat{\alpha} \in \Gamma(S^*, Ker\,\delta^1(k)) $ \ tel que \ $\mathcal{N}(\hat{\alpha}) $ \ soit la restriction \`a \ $ S^*$ \ d'un \'el\'ement de \ $  \Gamma(Y,, Ker\,\delta^1(k)) $, et induisant la classe de \ $\alpha$ \ dans  
$$ \frac{\mathcal{K}}{H^1(S^*, \mathbb{C})} \simeq H^1_{\lbrace 0 \rbrace}(S, H^1(0)) .$$

\bigskip

\noindent \textit{Preuve.} Pour \ $ k \gg 1 $ \ on peut repr\'esenter \ $\alpha$ \ par \ $ \tilde{\alpha} \in \Gamma(S^*, Ker\,\delta^1(k)$ ; de plus on peut \'egalement trouver \ $ \beta \in \Gamma(S, Ker\,\delta^1(k)) $ \ et \ $\Gamma \in \Gamma(S^*, \delta\mathcal{E}^0(k) $ \ v\'erifiant
$$ \mathcal{N}(\tilde{\alpha}) = \beta + \Gamma .$$
Comme \ $ \underset{k \rightarrow \infty}{\lim} h^1(k) \simeq H^1(0) $ \ est \`a support dans \ $S$, on peut supposer \'egalement que \ $\beta$ \ est la restriction \`a \ $S$ \ de \ $\tilde{\beta} \in \Gamma(Y, Ker\,\delta^1(k)$.

\noindent Comme la monodromie est l'identit\'e sur \ $ H^1(S^*, \mathbb{C}) $, on peut trouver \ $ \gamma \in \Gamma(S^*, \mathcal{E}^0(k)) $ \ v\'erifiant \ $ \mathcal{N}(\Gamma) = \delta \gamma .$
On aura alors
$$ \tilde{\Gamma} : = \begin{pmatrix} {\Gamma_k + \frac{df}{f} \wedge \gamma_k } \\ {0} \\ {\vdots} \\  {0 } \end{pmatrix} = \begin{pmatrix} {\Gamma_k} \\ {\Gamma_{k-1} } \\ {\vdots} \\ {0} \end{pmatrix} - \delta \begin{pmatrix} {0} \\ {\gamma_k} \\ {\vdots} \\ {\gamma_1} \end{pmatrix} .$$
Donc, quitte \`a ajouter un cobord  (et \`a changer \ $k$ \ en \ $k+1$) on peut supposer que
$$ \tilde{\alpha} = \tilde{\beta} + \tilde{\Gamma} $$
avec \ $ \Gamma_j = 0 $ \ pour \ $ j \not= k $.

\noindent Posons alors \ $ \hat{\alpha}_{k+1} = \alpha_k - \tilde{\Gamma_k} $ \ et \ $ \hat{\alpha}_j = \beta_j \quad \forall j \in [1,k] $. Alors \ $\hat{\alpha}$ \ v\'erifie bien les propri\'et\'es requises. $\hfill \blacksquare$

\subsection{ Lemme.}

 Dans la situation du lemme (8.11) ci-dessus, il existe (quitte \`a restreindre \ $X$ ) une $n-$forme semi-m\'eromorphe  \ $\hat{w}$ \ qui est d\'efinie et \ $d-$ferm\'ee  dans \ $X^* = X \setminus\lbrace 0 \rbrace$ \ et \`a support dans \ $mod\, S$, qui,\ restreinte \`a  \ $X\setminus \bar{X}'\setminus Y $, induit l'image de la restriction de \ $i[\alpha]$ \ dans \ $H^n_{mod\, S}(X \setminus \bar{X}' \setminus Y, \mathbb{C}) $.

\bigskip

\noindent \textit{Preuve.} Fixons un voisinage ouvert \ $ \mathcal{W}$ \ de \ $S^*$ \ dans \ $X^*$ \
et supposons que, quitte \`a restreindre le disque de centre \ $0$ \ dans \ $\mathbb{C}$ \ qui d\'efinit \ $X$, la section \ $\beta$ \ construite au lemme (8.11) est d\'efinie et \ $\delta-$ferm\'ee sur \ $X$. Soit \ $\chi \in \mathcal{C}^{\infty}(\mathcal{W})$ \ valant identiquement  1  pr\`es de \ $S^*$ \ et identiquement nulle pr\`es de \ $\partial\mathcal{W}$ . Nous la prolongerons par \ $0$ \ sur \ $X^*$. Posons  alors 
\begin{align*}
\hat{w} = d\chi\wedge \hat{\alpha}_{k+1} - (1 - \chi)\frac{df}{f}\wedge\beta_k  \\
\beta'_j = \beta_j \quad \forall j\in[1,k] \quad {\rm et} \quad \beta'_{k+1} = 0 .
\end{align*}
On a alors 
$$ d\chi\wedge \hat{\alpha} = j_{1,k+1}(\tilde{w}) + \delta ( (1 - \chi).\beta' ) .$$
Comme \ $(1 - \chi).\beta' $ \ est \`a support dans \ $mod\, S$ , on en d\'eduit que la $n-$forme semi-m\'eromorphe \ $\hat{w}$ \ v\'erifie les propri\'et\'es demand\'ees. $\hfill \blacksquare $

\bigskip

\noindent {\bf Remarque importante.}

\smallskip

\noindent Si la classe d\'efinie par \ $\hat{w}$ \ dans \ $H^n_{mod\,S}(X^*\setminus Y^*, \mathbb{C})$ \ est nulle, on peut trouver une $(n-1)-$forme semi-m\'eromorphe \ $\varepsilon$ \ sur \ $X^*$ \  \`a p\^oles dans \ $Y^*$ \ et \`a support dans \ $mod\,S$\footnote{ce qui signifie que l'adh\'erence de  son support dans \ $X^*$ \  ne rencontre pas \ $S^*$.} et v\'erifiant
$$ \hat{w} = d\varepsilon \quad {\rm sur} \ X^* .$$
Posons alors
$$ \gamma : = \chi.\hat{\alpha} - (1 - \chi).\beta' - j_{1,k+1}(\varepsilon) .$$
Comme \ $\delta \gamma = 0 $ \ cela d\'efinit une classe dans
$$ \mathbb{H}^{n-1}(X^*, (\mathcal{E}^{\bullet}(k+1), \delta^{\bullet})$$
qui s'identifie au  $(n-1)-$i\`eme groupe de cohomologie de la fibre de Milnor de \ $f$ \ \`a l'origine. On a donc un \'el\'ement  de \ $ H^0(S, H^{n-1}(0))$. Mais ceci montre, puisque \ $\varepsilon \ {\rm et} \ (1-\chi)$ \ sont nulles au voisinage de \ $S^*$ \ que la section correspondante sur \ $S$ \ du faisceau \ $H^{n-1}(0)$ \ prolonge la section "initiale" \ $\alpha \in H^0(S^*, H^{n-1}(0))$. Ceci montre que la nullit\'e de \ $\hat{w}$ \ dans l'espace  \ $H^n_{mod\,S}(X^*\setminus Y^*, \mathbb{C})$ \ implique celle de la classe de \ $\alpha$ \ dans \ $H^1_{\lbrace 0 \rbrace}(S, H^{n-1}(0))$.

\subsection{Variante pour \ $n =2$.}

L'\'enonc\'e est ici  le m\^eme qu'en (8.13) sauf que l'on obtient dans l'espace  \ $H^2_{mod\, S}((X \setminus \bar{X}' )\setminus Y, \mathbb{C}) $ \ qu'une \'egalit\'e {\bf modulo} l'image dans cet espace  du sous-espace \ $i\circ j(H^1(S^*, \mathbb{C}))$.

\noindent La preuve est analogue.

\smallskip

\noindent La remarque importante qui suit le (8.13) est valable pour \ $n = 2$. Dans ce cas, elle donne que la classe initiale de \ $\mathcal{K}$ \ dont on part, est dans \ $ j(H^1(S^*, \mathbb{C}))$ \ d\`es que la classe de \ $\hat{w}$ \ dans \ $H^n_{mod\,S}(X^*\setminus Y^*, \mathbb{C})$ \ est nulle.

\bigskip

\noindent On remarquera que l'image dans \ $H^2_{c\,\cap\,S}(0)$ \ de  \ $ j(H^1(S^*, \mathbb{C}))$ \  (modulo laquelle on travaille) est l'image par la restriction naturelle de \ $H^2_{mod\, S}(X^*, \mathbb{C})$. Ceci r\'esulte de notre construction et de l'isomorphisme
$$ H^2_{mod\, S}(X^*, \mathbb{C}) \simeq H^1(S^*, \mathbb{C}) $$
qui d\'ecoule facilement des annulations des \ $ H^i(X^*, \mathbb{C}) $ \ pour \ $i = 1,2 $.

\bigskip

\subsection{}

 Notre strat\'egie, pour terminer la preuve de l'injectivit\'e de la variation pour \ $n \geq 3$, consiste \`a consid\'erer le diagramme commutatif suivant dont les lignes sont exactes d'apr\`es le th\'eor\`eme (5.1)\footnote{Touts ces groupes de cohomologie sont \`a valeurs dans \ $ \mathbb{C}.$} :
$$\xymatrix{H^n_{mod\,S}(X^*) \ar[r] & H^n_{mod\,S}(X^*\setminus Y^*) \ar[d]^{r} \ar[r]^{Res} & H^{n-1}_{mod\,S}(Y^*) \ar[d]^{r'} \ar[r] & H^{n+1}_{mod\,S}(X^*)  \\
\cdots  \ar[r] & H^n_{mod\,S}(X\setminus \bar{X}'\setminus Y) \ar[r]^{Res'} & H^{n-1}_{mod\,S}(Y\setminus \bar{Y}') \ar[d]^{\partial} \ar[r] & \cdots \\
\quad & \quad & H^n_c(0) & \quad} $$
et \`a montrer les propri\'et\'e suivantes qui permettent facilement de conclure  que si \ $ \partial\circ Res' $ \ annule l'image de \ $\hat{w}$ \  dans \ $H^n_{mod\,S}(X\setminus \bar{X}'\setminus Y)$ \  alors c'est que cette image est d\'eja nulle ; mais alors l'image de \ $ d\chi\wedge \tilde{\alpha}$ \ est \'egalement nulle dans ce m\^eme espace. Le th\'eor\`eme  (4.6) donne alors  que la classe de \ $ \alpha \in H^1_{\lbrace 0 \rbrace}(S, H^{n-1}(0)) $ \ est nulle. Ceci ach\`eve donc la preuve de l'injectivit\'e de la variation pour \ $n \geq 3$, modulo la d\'emonstration des assertions suivantes:

\begin{itemize}
\item{1)} La restriction de \ $\hat{w}$ \ est dans l'image de l'application de restriction not\'ee  \ $r$ .
\item{2)} L'application de r\'esidu \ $Res$ \ est un isomorphisme, grace \`a l'annulation des groupes \ $H^i_{mod\,S}(X^*) \quad {\rm pour} \quad i = n, n+1$ \ pour \ $n \geq 3 $.
\item{3)} La compos\'ee  \ $\partial\circ r' $ \ est injective.
\end{itemize}

\noindent Les lemmes (8.11) et le (8.13) donnent d\'eja  la propri\'et\'e 1).

\noindent  Montrons donc les annulations qui donnent la seconde propri\'et\'e pour \ $n \geq 3$. 

\noindent Comme on a \ $ H^i(X^*) = 0 \quad {\rm pour} \quad i = n, n+1$ \ il nous suffit de montrer que si\ $ \psi \in \mathcal{C}^{\infty}(X^*) $ \ est de degr\'e \ $n-1$ \ ou \ $n$ \ et a une diff\'erentielle identiquement nulle au voisinage de \ $S^*$ \ il existe \ $ \chi \in \mathcal{C}_{mod\,S}^{\infty}(X^*)$ \ telle que \ $d\chi = d\psi$ \ sur \ $X^*$. Comme on a \ $n \geq 3$ \ on a une base de voisinages ouverts \ $\mathcal{W}$ \  de \ $S^*$ \ dans \ $X^*$ \ qui v\'erifient \ $ H^i(\mathcal{W}) = 0 $ \  pour \ $ i \geq 2 $. Ceci permet ais\'ement de conclure.

\bigskip

\subsection{}

Il nous reste \`a montrer l'injectivit\'e de la fl\`eche
$$ \partial\circ r' : H^{n-1}_{mod\,S}(Y^*) \rightarrow H^n_c(0) $$
qui est d\'efinie comme suit : \`a  \ $\psi \in \mathcal{C}^{\infty}_{mod\,S}(\mathcal{V}(Y^*))$ \ v\'erifiant \ $d\psi = 0 $ \ on associe la restriction \`a \ $\lbrace f=s \rbrace$ \ pour \ $s$ \ assez voisin de \ $0$ \ de la forme \ $ d\rho \wedge \psi $ \ o\`u \ $ \rho \in \mathcal{C}^{\infty}_c(X)$ \ vaut identiquement au voisinage de l'origine.

\noindent Cette injectivit\'e va r\'esulter d'un th\'eor\`eme g\'en\'eral.

\bigskip

\subsection{Th\'eor\`eme.}

 \textit{Notons respectivement par \ $ \mathbb{H}^p_c(k) $ \ et \ $ \mathbb{H}^p_{c\,\cap\,S}(k)$ \ les groupes d'hypercohomologie \`a supports \ $ \mathbb{H}^p_c(Y, \mathcal{E}(k)^{\bullet}, \delta^{\bullet})$ \ et \ $ \mathbb{H}^p_{c\,\cap\,S}(Y, \mathcal{E}(k)^{\bullet}, \delta^{\bullet})$. 
Alors on a, pour chaque \ $ k \geq 1 $ \ la suite exacte longue d'espaces vectoriels monodromiques :}
$$\rightarrow \mathbb{H}^{p-1}_{mod\,S}(\partial Y, \mathcal{E}^{\bullet}(k), \delta^{\bullet}) \overset{j}{\rightarrow}  \mathbb{H}^p_c(k)  \overset{can_{c,S}}{\longrightarrow} \mathbb{H}^p_{c\,\cap\,S}(k) \overset{\eta}{\rightarrow}  \mathbb{H}^{p}_{mod\,S}(\partial Y, \mathcal{E}^{\bullet}(k), \delta^{\bullet})\rightarrow$$
\textit{o\`u l'on a not\'e par \ $ \mathbb{H}^{p}_{mod\,S}(\partial Y, \mathcal{E}^{\bullet}(k), \delta^{\bullet})$ \ la limite inductive quand \ $ \varepsilon' \rightarrow \varepsilon$ \ de \ $ \mathbb{H}^{p}_{mod\,S}( Y\setminus\bar{Y}', \mathcal{E}^{\bullet}(k), \delta^{\bullet})$ \ avec \ $ Y' = Y \,\cap\,X' $ \ o\`u \ $X'$ \ est la trace sur \ $X$ \ de la boule de centre \ $0$ \ et de rayon} \ $\varepsilon'$\footnote{Rappelons que \ $X$ \ est l'intersection de la boule de centre $0$ \ et de rayon \ $\varepsilon$ \ avec \ $f^{-1}(D_{\eta}), \quad \eta \ll \varepsilon.$}.

\bigskip

\subsection{D\'emonstration.}

 Commen{\c c}ons par d\'efinir l'application \ $j$. Fixons \ $ X': = B(0, \varepsilon') \cap f^{-1}(D) \subset X = B(0, \varepsilon) \cap f^{-1}(D), \varepsilon - \varepsilon' \ll\varepsilon \ll 1 $. Soit \ $\varphi$ \ une \ section de \ $\mathcal{E}^{p-1}(k)$  \ sur \ $ X\setminus \bar{X'} $ \ identiquement nulle au voisinage de \ $S$ \ et \ $\delta-$ferm\'ee. Soit,  de plus,\ $\rho \in \mathcal{C}^{\infty}(X) $, une fonction identiquement \'egale \`a  1  au voisinage de \ $\bar{X'}$ et \`a support \ $f-$propre. Alors \ $ j(\varphi) = d\rho \wedge \varphi = \delta((\rho - 1).\varphi)$ \ est dans \ $ \Gamma_c(Y, Ker\,\delta \cap \mathcal{E}^p(k))$; elle d\'efinit donc un \'el\'ement de \ $ \mathbb{H}^p_c(k)$. Montrons que ceci d\'efinit bien une application lin\'eaire 
$$ j : \mathbb{H}^{p-1}_{mod\,S}(\partial Y, \mathcal{E}^{\bullet}(k), \delta^{\bullet}) \rightarrow  \mathbb{H}^p_c(k)   .$$
Si \ $ \varphi$ \ est l'image par \ $\delta$   d'une section \ $ \psi \in \Gamma_{mod\, S}(\partial Y, \mathcal{E}^{p-2}(k)) $ \ alors on aura  
$$ d\rho \wedge\delta \psi = \delta(-d\rho\wedge\psi) .$$ 
o\`u l'on a \ $ d\rho\wedge \psi \in \Gamma_c(Y, \mathcal{E}^{p-1}(k)) $, en prolongeant par \ $0$.

\noindent On a clairement \ $j\circ _k\mathcal{N} =   _k\mathcal{N}\circ j $ \ ce qui montre l'aspect monodromique de l'application lin\'eaire \ $j$.

\noindent On remarquera que pour \ $k = 1$ \ cette application est bien compatible avec  l'application \ $\partial$ \ consid\'er\'ee plus haut.

\noindent Identifions le noyau de \ $j$ ; si on a  \ $ d\rho \wedge \varphi = \delta\gamma$ \ o\`u \ $\gamma \in \Gamma_c(Y, \mathcal{E}^{p-1}(k)) $, on aura, quitte \`a augmenter la taille de \ $ X' $, que la restriction \ $\tilde{\gamma}$ \ de \ $\gamma$ \ \`a \ $X\setminus \bar{X'}$ \ sera dans \ $ \Gamma_{mod\,S}(\partial Y, \mathcal{E}^{p-1}(k)) $ \ et v\'erifiera \ $ d\rho \wedge \varphi = \delta\tilde{\gamma} .$ Alors \ $ \gamma + (1-\rho).\varphi $ \ d\'efinira un \'el\'ement de \ $ \mathbb{H}^{p-1}_{c\,\cap\,S}(Y, \mathcal{E}^{\bullet}(k), \delta^{\bullet}) .$ La construction de l'analogue\ $\eta$ \  en degr\'e \ $(p-1)$ \ de l'application \ $\eta$ \ que nous allons d\'efinir maintenant montrera que le noyau de \ $j$ \ est bien l'image de \ $ \mathbb{H}^{p-1}_{c\, \cap \,S}(k) $ \ par \ $\eta$.

\noindent La construction de l'application \ $\eta$ \ est tr\`es simple. 

\noindent Si \ $ w \in \Gamma_{c\,\cap\,S}(Y, Ker\,\delta \cap \mathcal{E}^{p}(k))$ \ pour un choix de \ $X'$ \ assez gros, la restriction de \ $w$ \ \`a \ $Y \cap (X\setminus\bar{X'})$ \ est \`a support dans \ $mod\,S$ \ et donne une classe dans \ $ \mathbb{H}^{p}_{mod\,S}(\partial Y, \mathcal{E}^{\bullet}(k), \delta^{\bullet}) $. On v\'erifie imm\'ediatement que ceci passe au quotient. L'application est donc bien d\'efinie et elle est clairement lin\'eaire et monodromique. 

\noindent De plus, il est clair que si \ $w$ \ est \`a support compact, son image par \ $\eta$ \ est nulle.

\smallskip

\noindent R\'eciproquement, si \ $\eta[w] = 0 $, alors on peut \'ecrire sur \ $\partial Y$ :
$$ w = \delta v $$
avec \ $v \in \Gamma_{mod\,S}(\partial Y, \mathcal{E}^{p-1}(k)) $. Alors pour \ $\rho \in \mathcal{C}^{\infty}_{c/f}(X) $ \ valant identiquement sur \ $\bar{X'}$ \ avec \ $X'$ \ assez gros, on aura \ $ w' : = w - \delta((1-\rho).v) $ \ qui repr\'esentera la classe \ $[w] \in  H^p_{c\,\cap\,S}(0)$ \ et sera \`a support compact sur \ $Y$. Donc \ $[w]$ \ est bien dans l'image de \ $can_{c,S} .$

\noindent Pour achever la d\'emonstration il nous reste \`a montrer que le noyau de l'application
$$ can_{c,S} : \mathbb{H}^p_c(k) \rightarrow  \mathbb{H}^p_{c\,\cap\,S}(k) $$
est bien l'image de  $j$. Soit donc \ $v \in \Gamma_{c\,\cap\,S}(Y, Ker\,\delta \cap \mathcal{E}^{p-1}(k)) $ \ tel que \ $w = \delta v $ \ soit \`a support compact dans \ $Y$. Alors \ $v$ \ d\'efinit par restriction \`a \ $\partial Y$ \ un \'el\'ement \ $[v]$ de \ $ \mathbb{H}^{p-1}_{mod\,S}(\partial Y, \mathcal{E}^{\bullet}(k), \delta^{\bullet})$. Pour calculer \ $j([v])$ \ choisissons \ $\rho$ \ comme ci-dessus\footnote{et tel que l'on ait  \ $Supp(d\rho) \cap Supp(\delta v) = \emptyset $.}. Alors on aura \ $ j([v]) = [d\rho\wedge v] = \delta((\rho - 1).v) $. Mais on a 
$$ \delta v - \delta(\rho .v) =  \delta((1-\rho).v) $$
avec \ $(1-\rho).v \in \Gamma_c(Y, \mathcal{E}^{p-1}(k)) $ \ et donc la classe de \ $j([v])$ \ dans \ $\mathbb{H}^p_c(k)$ \ co\"incide bien avec celle de \ $ w = \delta v . \hfill \blacksquare$

\bigskip

\noindent Pour terminer la preuve de l'injectivit\'e de la variation (c'est \`a dire pour d\'eduire l'assertion 3) du (8.15) du th\'eor\`eme (8.17)), il nous suffit de montrer que pour \ $n \geq 3$ \ la fl\`eche \ $\eta$ \ en degr\'e \ $n-1$ \ de la suite exacte longue de (8.17) est nulle, ce qui r\'esulte de la proposition suivante.

\bigskip

\subsection{Proposition.}

 Sous les hypoth\`eses standards pour \ $n \geq 3$ \ l'application canonique 
$$ can_{c,S} :  \underset{k \rightarrow \infty}{lim}  \mathbb{H}^{n-1}_c(k) \rightarrow  \underset{k \rightarrow \infty}{lim}  \mathbb{H}^{n-1}_{c\,\cap\,S}(k) $$
est un isomorphisme. De plus, pour \ $n \geq 4 $ \ ces deux groupes sont nuls.

\bigskip

\noindent \textit{Preuve.} Commen{\c c}ons par montrer le lemme suivant

\bigskip

\subsection{Lemme.}

Pour \ $ n \geq 4 $ \ on a \ $ H^{n-1}_{c\,\cap\,S}(Y, \mathbb{C}) = 0 $. Pour \ $ n = 2 , 3 $ \ on a l'isomorphisme \ $  H^{n-1}_{c\,\cap\,S}(Y, \mathbb{C}) \simeq H^{n-1}_{\lbrace 0 \rbrace}(S, \mathbb{C} )$.

\bigskip

\noindent \textit{Preuve.} Comme on peut remplacer \ $Y$ \ par la limite inductive de ses voisinages ouverts, on peut calculer \ $  H^{n-1}_{c\,\cap\,S}(Y, \mathbb{C})$ \ \`a l'aide  du complexe de de Rham des germes de formes diff\'erentielles \`a support dans \ $ c\,\cap\,S $. Soit donc \ $\varphi $ \ un tel germe d-ferm\'e de degr\'e \ $n-1$. Comme \  $H^{n-1}(Y, \mathbb{C}) = 0 $ \ puisque \ $Y$ \ est contractible, on peut trouver un germe \ $\psi$ \ de degr\'e \ $n-2$ \ v\'erifiant \ $ d\psi = \varphi $. Si \ $K$ \ est un compact de \ $S$ \ assez gros, on aura \ $ d\psi = 0 $ \ au voisinage de \ $S\setminus K$. On constate alors que si la classe de cohomologie d\'efinie par \ $\psi$ \ dans \ $H^{n-2}(S\setminus K, \mathbb{C})$ \  est nulle pour \ $n \geq 3$ \ ou prolongeable \`a \ $S$ \ pour \ $n = 2$ \ alors la classe d\'efinie par \ $\varphi$ \ est nulle dans \ $  H^{n-1}_{c\,\cap\,S}(Y, \mathbb{C})$. Ceci permet alors de construire un isomorphisme \ $ H^{n-1}_{c\,\cap\,S}(Y, \mathbb{C}) \simeq H^{n-1}_{ \lbrace 0 \rbrace}(S, \mathbb{C} )$ \ ce qui ach\`eve la preuve puisque pour \ $ n \geq 4 $ \ le groupe \ $ H^{n-1}_{\lbrace 0 \rbrace}(S, \mathbb{C} )$ \ est nul. $\hfill \blacksquare$

\bigskip

\subsection{Preuve de la proposition (8.19).}

Supposons maintenant \ $ n \geq 3 $ \ et notons par \ $\Phi $ \ une famille paracompactifiante de ferm\'es de \ $Y$ \ \'egale soit \`a \ $c$ \ ou bien \`a \ $c\,\cap\,S $. Comme les faisceaux \ $\mathcal{E}^{\bullet}(k) $ \ sont fins, on aura
$$ \mathbb{H}^{n-1}_{\Phi}(Y, \mathcal{E}^{\bullet}(k), \delta^{\bullet}) \simeq \Gamma_{\Phi}(Y, Ker\,\delta^{n-1}(k))\big/ \delta \Gamma_{\Phi}(Y, \mathcal{E}^{n-2}(k)) .$$
Soit maintenant \ $w \in  \Gamma_{\Phi}(Y, Ker\,\delta^{n-1}(k)) $. Alors si \ $K$ \ est un compact de \ $Y$ \ assez gros, $w$ \ est identiquement nulle au voisinage de \ $S\setminus K$. La suite exacte
$$ 0 \rightarrow \delta \mathcal{E}^{n-2}(k) \rightarrow Ker\,\delta^{n-1} \rightarrow H^{n-1}(0) \rightarrow 0 . $$
donne alors 
$$ 0 \rightarrow \Gamma_{\Phi}(Y, \delta \mathcal{E}^{n-2}(k)) \rightarrow  \Gamma_{\Phi}(Y, Ker\,\delta^{n-1} ) \overset{r}{\rightarrow} \Gamma_{\Phi}(Y, h^{n-1}(k)) \rightarrow H^1_{\Phi}(Y, \delta \mathcal{E}^{n-2}(k)) \cdots   $$
Et comme, pour \ $n \geq 3$ \ le faisceau \ $ h^{n-1}(k)$ \ est \`a support dans \ $S$, on aura 
$$ r(w)\vert_{S\setminus K} \equiv 0 \quad {\rm dans} \quad \Gamma_{\Phi}(S\setminus K,  h^{n-1}(k)) \simeq  \Gamma_{\Phi}(Y\setminus K,  h^{n-1}(k)). $$
Comme, de plus, \ $ h^{n-1}(k)$ \ est un syst\`eme local sur \ $S^*$ \ et n'admet pas de section non nulle \`a support l'origine, on aura \ $r(w) = 0 $ \ et on en conclut que  \ $ w \in \Gamma_{\Phi}(Y, \delta \mathcal{E}^{n-2}(k)).$

\noindent Les suites exactes
$$ 0 \rightarrow \delta\mathcal{E}^{q-1}(k) \rightarrow \mathcal{E}^{q}(k)\rightarrow \delta\mathcal{E}^{q}(k) \rightarrow 0 $$
pour \ $ 2 \leq q \leq n-2 $ \ conduisent aux isomorphismes, puisque les faisceaux \ $ \mathcal{E}^q(k)$ \ sont fins et la famille de supports \ $ \Phi $ \ est paracompactifiante
$$ \Gamma_{\Phi}(Y,  \delta \mathcal{E}^{n-2}(k))\big/ \delta \Gamma_{\Phi}(Y, \mathcal{E}^{n-2}(k)) \simeq H^1_{\Phi}(Y, \delta \mathcal{E}^{n-3}(k)) \simeq \cdots \simeq H^{n-3}_{\Phi}(Y, \delta \mathcal{E}^{1}(k)).$$
Bien sur pour \ $n=3 $ \ on doit remplacer \ $ H^{n-3}_{\Phi}(Y, \delta \mathcal{E}^{1}(k))$ \ par son quotient par \ $ \delta H^0_{\Phi}(Y, \mathcal{E}^1(k)) .$

\noindent La suite exacte
$$ 0 \rightarrow Ker\,\delta^1(k) \rightarrow \mathcal{E}^1(k) \rightarrow \delta \mathcal{E}^{1}(k)) \rightarrow 0 $$
donne alors l'isomorphisme \ $ H^{n-3}_{\Phi}(Y, \delta \mathcal{E}^{1}(k)) \simeq H^{n-2}_{\Phi}(Y, Ker\,\delta^1(k)) .$ 

\noindent Les suites exactes
\begin{align*}
0 \rightarrow \delta \mathcal{E}^0(k) \rightarrow Ker\,\delta^1(k) \rightarrow h^1(k) \rightarrow 0  \\
 0 \rightarrow h^0(k) \rightarrow \mathcal{E}^0(k) \rightarrow \delta \mathcal{E}^0(k) \rightarrow 0
 \end{align*}
 combin\'ees avec le fait que \ $ j_{k,k+1} : h^1(k) \rightarrow h^1(k+1) $ \ est nulle, donnent un isomorphisme 
\begin{align*}
 \underset{k \rightarrow \infty}{lim} \mathbb{H}^{n-1}_{\Phi}(k) \simeq \underset{k \rightarrow \infty}{lim} H^{n-2}_{\Phi}(Y, Ker\,\delta^1(k)) \simeq \underset{k \rightarrow \infty}{lim} H^{n-2}_{\Phi}(Y, \delta\mathcal{E}^0(k)) \\
  \simeq \underset{k \rightarrow \infty}{lim} H^{n-1}_{\Phi}(Y, h^0(k)) \simeq H^{n-1}_{\Phi}(Y, \mathbb{C}).
  \end{align*}
  Pour terminer la preuve de la proposition (8.17), il suffit donc de voir que l'application canonique
  $$ H^{n-1}_{c\,\cap\,S}(Y, \mathbb{C}) \rightarrow H^{n-1}_c(Y, \mathbb{C}) $$
  est un isomorphisme pour \ $ n \geq 3 $, et entre groupes nuls pour \ $n \geq 4 $.
  
  \noindent Ceci r\'esulte du lemme (8.20). \ $ \hfill \blacksquare $
  
  \bigskip

\noindent Ceci ach\`eve donc la preuve de l'injectivit\'e de la variation pour \ $n \geq 3$. On a donc \'egalement \'etabli la non d\'eg\'en\'erescence de la forme hermitienne canonique (voir (7.1)) et l'\'egalit\'e de dimension 
$$ \dim_{\mathbb{C}}\, \big( Im\, \theta\big) = \dim_{\mathbb{C}}\,\big( H^1_{\lbrace 0 \rbrace}(S, H^{n-1}(0))\big).$$
Cette \'egalit\'e implique que \ $ Im\, \theta$ \ co{\"i}ncide avec l'orthogonal dans \ $ H^1(S^*, H^{n-1}(0))$ \ du sous espace \ $ H^0(S,  H^{n-1}(0))$ \ de \ $H^0(S^*,  H^{n-1}(0))$ \ pour l'accouplement non d\'eg\'en\'er\'e d\'eduit de la forme hermitienne canonique sur le syst\`eme local \ $ H^{n-1}(0))$\ (voir (3.1)).

\bigskip

\subsection{Le cas \ $n = 2$.}

Nous allons montrer que les r\'esultats pr\'ec\'edents s'\'etendent "mutadis mutandis" au cas \ $n =2 $ \ quitte \`a remplacer l'espace vectoriel monodromique \ $H^2_{c\,\cap\,S}(0) $ \ par le quotient
$$ \frac{H^2_{c\,\cap\,S}(0)}{H^1(S^*, \mathbb{C})} $$
o\`u l'inclusion \ $j$ \ de \ $H^1(S^*, \mathbb{C})$\footnote{muni de la monodromie triviale, c'est \`a dire \'egale \`a l'identit\'e.}a \'et\'e d\'efinie au Lemme (4.8) (voir aussi (4.12.1)).

\noindent Nous allons seulement pr\'eciser les points de la d\'emonstration du cas \ $n \geq 3$ \ qui sont \`a modifier de fa{\c c}on significative dans ce cas.

\bigskip

\subsubsection{}

Commen{\c c}ons par remarquer que toute la premi\`ere partie de la preuve de l'injectivit\'e de la variation, \`a savoir les sections (8.1) \`a (8.10), reste valable telle quelle\footnote{Grace, entre autres, au fait que le sous-espace \ $j(H^1(S^*, \mathbb{C}))$\ de \ $H^2_{c\,\cap\,S}(0)$ \ est form\'e de vecteurs invariants par la monodromie.}

\noindent Il n'en n'est pas de m\^eme pour la seconde partie de cette preuve.
\begin{itemize}
\item{1)} Les lemmes (8.11) et (8.13) doivent \^etre remplac\'es par les variantes (8.12) et (8.14) qui les suivent.
\item{2)} L'assertion 2) de (8.15) doit \^etre remplac\'ee par l'assertion suivante:

\smallskip

\noindent \quad  L'application \ $Res$ \ est surjective et son noyau a pour image par \ $r$ \ l'image injective par la restriction \`a \ $ (X \setminus \bar{X}' )\setminus Y $ \ du sous espace \ $j(H^1(S^*, \mathbb{C}))$ \ de \ $H^2_{c\,\cap\,S}(0)$. 

\noindent  Ceci est une cons\'equence facile de la remarque qui conclut le (8.14)

\item{3)} On a pour \ $n = 2$, \ $H^1_c(0) \simeq 0 $ \ et \ $ H^1_{c\,\cap\,S}(0) \simeq H^1_{\lbrace 0 \rbrace} (S,\mathbb{C}) $. De plus, un \'el\'ement de \ $\mathcal{K}\big/ j(H^1(S^*, \mathbb{C})$ \ invariant par la monodromie donne une classe nulle par l'application \ $Res$ \ du diagramme de (8.17) si et seulement si cet \'el\'ement est nul.

\end{itemize}

\bigskip


\subsubsection{Preuve de l'assertion 3).}
L'annulation de \ $H^1_c(0)$ \ r\'esulte du fait que la fibre de Milnor \ $F$ \ de \ $f$ \ \`a l'origine est une vari\'et\'e de Stein de dimension  2 ; donc  \ $H^3(F, \mathbb{C}) \simeq 0$ \  et par dualit\'e de Poincar\'e on a \ $H^1_c(F, \mathbb{C}) \simeq 0 .$

\noindent L'isomorphisme 
 $$ \frac{ H^0(S^*, \mathbb{C}) }{ H^0(S, \mathbb{C})}\simeq H^1_{\lbrace 0 \rbrace} (S, \mathbb{C}) \rightarrow H^1_{c\,\cap\,S}(Y, \mathbb{C})$$
 est donn\'e par la fl\`eche
$$ \varepsilon : H^0(S^*, \mathbb{C}) \rightarrow H^1_{c\,\cap\,S}(Y, \mathbb{C}) $$
d\'efinie de la fa{\c c}on suivante : \`a la fonction localement constante  \ $\varphi$ \ au voisinage de \ $S^*$ \ on associe la classe \ $[\varphi.d\chi ]$ \ o\`u \ $\chi \in \mathcal{C}^{\infty}(Y) $ \ vaut identiquement au voisinage de \ $\partial S $ \ et s'annule identiquement d\`es que l'on s'\'eloigne de \ $\partial S$. On constate facilement que ceci ne d\'epend pas du choix de la fonction \ $\chi$ \ et passe au quotient par \ $H^0(S, \mathbb{C})$\footnote{Si \ $\varphi$ \ est constante au voisinage de \ $S$ \ on peut prolonger la constante \`a \ $Y$ \ et alors \ $ \varphi.d\chi = d\big((1 - \chi).\varphi \big) $ \ montre que la classe \ $ [\varphi.d\chi] $ \ est nulle dans \ $H^1_c(0)$.}.

\noindent L'application induite par \ $\varepsilon$ \ est injective :  si \ $ \varphi.d\chi = d\psi $ \ o\`u \ $ \psi \in \mathcal{C}^{\infty}_{c\,\cap\, S}(Y)$, on aura \ $ d(\psi - \chi.\varphi) = 0 $ \ sur \ $Y$, et donc \ $\psi - \chi.\varphi $ \ sera localement constante (donc constante) sur \ $Y$. Comme \ $\psi$ \ est nulle et \ $\chi$ \ vaut identiquement \ $1$ \ pr\`es de \ $\partial S$, ceci montre que \ $\varphi$ \ se prolonge en une constante sur \ $Y$ \ et donc induit  la classe nulle dans \ $H^1_{\lbrace 0 \rbrace} (S, \mathbb{C})$.

\noindent L'application  induite par \ $\varepsilon$ \ est surjective : si \ $\alpha \in \mathcal{C}^{\infty}_{c\,\cap\, S}(Y)^1 $ \ est \ $d-$ferm\'ee, on peut l'\'ecrire, puisque \ $Y$ \ est contractible, \ $ \alpha = d\beta $ \ avec \ $ \beta \in \mathcal{C}^{\infty}(Y)^0$. On en d\'eduit que \ $\beta$ \ est localement constante au voisinage de \ $\partial S$. Soit \ $\varphi$ \ la fonction localement constante au voisinage de \ $S^*$ \ qu'elle d\'efinie. Alors l'\'egalit\'e \ $\alpha - d(\chi.\varphi) = d\big((1 - \chi).\beta\big) $ \ montre, puisque \ $ (1 - \chi).\beta $ \ est \`a support dans \ $ c\,\cap\, S$ \ que l'on a \ $\varepsilon([\varphi]) = [\alpha]$.

\noindent Il nous reste \`a montrer que l'application \'evidente
$$ H^1_{c\,\cap\,S}(Y, \mathbb{C}) \rightarrow H^1_{c\,\cap\,S}(0) $$
est un isomorphisme.

\noindent Elle est injective : si \ $\alpha \in \mathcal{C}^{\infty}_{c\,\cap\, S}(Y)^1 $ \ est \ $d-$ferm\'ee et si \ $ j_{1,k}(\alpha) = \delta \beta $ \ avec \ $ \beta \in \Gamma_{c\,\cap\,S}(Y,\mathcal{E}^0(k))$, on a 
\begin{align*}
\quad \quad & \alpha = d\beta_k - \frac{df}{f} \wedge \beta_{k-1} \\
{\rm et} \quad & d\beta_j - \frac{df}{f} \wedge \beta_{j-1} = 0 \quad \forall j \in [1,k-1]
\end{align*}
avec la convention \ $\beta_0 : = 0$. Alors \ $ d\beta_1 = 0 $ \ et donc \ $\beta_1 = 0 $ \ puisque c'est une fonction localement constante et nulle pr\`es de \ $\partial S$. On en d\'eduit de m\^eme que les \ $\beta_j$ \ sont nulles pour \ $ j \in [1,k-1]$. Il nous reste alors seulement l'\'equation sur \ $Y$ \
$$ \alpha = d\beta_k .$$
Mais une fonction semi-m\'eromorphe dont la diff\'erentielle est \ $\mathcal{C}^{\infty}$ \ est \ $\mathcal{C}^{\infty}$\footnote{On se ram\`ene au cas m\'eromorphe par Dolbeault (local) et on conclut par le lemme de de Rham holomorphe.}. On a donc \ $[\alpha] = 0 $ \ dans \ $ H^1_{c\,\cap\,S}(Y, \mathbb{C})$.

\noindent Elle est surjective : si \ $ \gamma \in \Gamma_{c\,\cap\,S}(Y, Ker\,\delta^1(k))$, comme on a la nullit\'e de \ $H^0_{c\,\cap\,S}(Y, H^1(0))  $ \ ainsi que l'isomorphisme  \ $ \underset{k \rightarrow \infty}{\lim} \ h^1(k) \simeq H^1(0)$, quitte \`a augmenter \ $k$, on peut supposer que \ $ \alpha \in  \Gamma_{c\,\cap\,S}(Y, \delta\mathcal{E}^0(k)).$ Mais la suite exacte
$$ 0 \rightarrow h^0(k) \rightarrow \mathcal{E}^0(k)\rightarrow \delta \mathcal{E}^0(k)\rightarrow 0 $$
montre que l'on a un isomorphisme
$$ \partial : \frac{ \Gamma_{c\,\cap\,S}(Y, \delta\mathcal{E}^0(k))}{\delta  \Gamma_{c\,\cap\,S}(Y, \mathcal{E}^0(k))} \rightarrow H^1_{c\,\cap\,S}(Y, h^0(k)) \simeq  H^1_{c\,\cap\,S}(Y, \mathbb{C}).$$
On v\'erifie alors facilement que l'application consid\'er\'ee co{\"i}ncide avec \ $ \partial^{-1}$.

\bigskip

\noindent V\'erifions enfin la derni\`ere partie de l'assertion 3). On doit donc examiner le cas o\`u la $2-$ forme ferm\'ee \ $\hat{w}$ \ peut \^etre choisie \ $\mathcal{C}^{\infty} $ \ sur \ $X^*$. Mais comme on a, \ $ H^2_{mod\,S}(X\setminus\bar{X}', \mathbb{C}) \simeq H^1(S^*, \mathbb{C})$, cela signifie que l'on peut choisir \ $ \hat{w} = d\chi \wedge \varphi $ \ o\`u \ $\varphi$ \ est une $1-$forme \ $\mathcal{C}^{\infty} $ \ et \ $d-$ferm\'ee au voisinage de \ $S^*$ \ et o\`u la fonction \ $\chi$ est identiquement \'egale \`a  1  pr\`es de \ $S^*$. Mais alors la classe initiale est dans \ $j(H^1(S^*,\mathbb{C})$ \ d'apr\`es la remarque importante du (8.13) et de sa variante dans (8.14)

\bigskip

\noindent Terminons ce paragraphe 8  en redonnant l'\'enonc\'e de  notre r\'esultat :

\subsection{Th\'eor\`eme 2.}
\textit{Sous les hypoth\`eses standards pour \ $n \geq 3 $ \ on a les propri\'et\'es suivantes
\begin{itemize}
\item{1)} $\dim H^n_{c\,\cap\,S}(0) = \dim H^n(0) $.
\item{2)} La variation \ $ var : H^n_{c\,\cap\,S}(0) \rightarrow H^n_c(0) $ \ est un ismorphisme d'espace vectoriels monodromiques.
\item{3)} La forme hermitienne canonique
$$\mathcal{H} :  H^n_{c\,\cap\,S}(0)\times H^n(0) \rightarrow \mathbb{C} $$
est non d\'eg\'en\'er\'ee.
\end{itemize}
Pour \ $n = 2 $ \ tout ceci reste vrai \`a condition de remplacer \ $H^2_{c\,\cap\,S}(0)$ \ par son quotient par \ $j(H^1(S^*, \mathbb{C}))$.}

\section{Applications et un exemple pour conclure.}

\noindent Il n'est pas difficile de d\'eduire les analogues pour la valeur propre \ $1$ \ des th\'eor\`emes 13 et 14 de \ $[B.91]$ \ de l'\'egalit\'e entre \ $ Im(\theta) $ \ et l'orthogonale pour \ $\tilde{h}$ \ de \ $H^0(S, H^{n-1}(0))$ \ dans \ $H^1(S^*, H^{n-1}(0))$. On obtient les th\'eor\`emes suivants dont la ligne de d\'emonstration, maintenant que l'on sait que \ $Im(\theta)$ \ est exactement l'orthogonal de \ $H^0(S, H^{n-1}(0))$ \ dans \ $ H^1(S^*, H^{n-1}(0))$ \ pour \ $\tilde{h}$, suit pas \`a pas celle de [B.91].

\subsection{Th\'eor\`eme 3.}

\textit{Sous les hypoth\`eses standards pour la valeur propre 1, consid\'erons \ $e \in H^n(0) $ \ v\'erifiant \ $\mathcal{N}^k(e) = 0 $ \ avec \ $k \geq k_0$\footnote{Rappelons que, par d\'efinition, \ $k_0$ \ est l'ordre de nilpotence de la monodromie agissant sur le syst\`eme local \ $H^{n-1}(0)$ \ sur \ $S^*$.}. Alors les conditions suivantes sont \'equivalentes :}
\begin{itemize}
\item{1)} \ $\tilde{ob}_k(e) = 0$ .
\item{2)} \textit{Si \ $w \in \Gamma(Y, Ker\,\delta^n(k))$ \ v\'erifie \ $ r^n(k)(w) = e $, pour chaque \ $j \in \mathbb{Z} $, les fonctionnelles analytiques
$$ P_{k+l} \big(\lambda = 0, \int_X \vert f \vert^{2\lambda} \bar{f}^{-j}\frac{df}{f}\wedge \bar{w}_k \wedge \square \big) $$
sont nulles pour \ $l \geq 2$.}
\end{itemize}

\noindent \textit{De plus, pour v\'erifier 2) il suffit de le faire pour \ $l = 2$ \ et pour \ $j \in [0,n] $.}

\bigskip

\subsection{Th\'eor\`eme 4.}

\textit{Sous les hypoth\`eses standards pour la valeur propre 1, supposons que le prolongement m\'eromorphe de \ $ \int_X \vert f \vert^{2\lambda} \square $ \ admette un p\^ole d'ordre \ $ \geq k $ \ en un entier n\'egatif, avec \ $ k \geq \sup{(k_0, k_1)} + 2$ \ o\`u \ $k_0$ \ et \ $k_1$ \ sont respectivement les ordres de nilpotence de la monodromie agissant sur le syst\`eme local \ $ H^{n-1}(0)_{\vert S^*} $ \ et sur l'espace vectoriel \ $H^n(0)$.}

\noindent \textit{ Alors l'ordre des p\^oles aux entiers n\'egatifs assez grand est exactement \'egal \`a \ $k$ \ et l'on a \ $ k_0 \leq k_1 = k - 2 $ \ avec  \ $\tilde{ob}_{k_1} \not\equiv 0 .$}

\bigskip

\subsection{Exemple.}

Consid\'erons la fonction \ $ f(x,y,z) = x^2.(x^2 + y^2) + z^4 .$ Comme cet exemple est tr\`es similaire \`a celui \'etudi\'e en d\'etail \`a la fin de [B.91], nous allons seulement esquisser son \'etude.

\noindent Tout d'abord l'homog\'ene{\"i}t\'e de \ $f$ \ nous assure que la monodromie agit de fa{\c c}on semi-simple sur \ $H^2(0).$ Cependant, nous allons montrer que le prolongement m\'eromorphe de
$$\int_{X} \vert f \vert^{2\lambda}\rho. \vert z \vert^2 dx\wedge d\bar{x}\wedge dy\wedge d\bar{y}\wedge dz\wedge d\bar{z} $$
admet en \ $ \lambda = -1 $ \ un p\^ole triple, si \ $\rho \in \mathcal{C}^{\infty}_c(X)$ \ vaut identiquement 1  pr\`es de l'origine. Pour cela il suffit, par transformation de Mellin complexe, de montrer que le d\'eveloppement asymptotique, quand \ $ \vert s \vert \rightarrow 0$ \ de l'int\'egrale
$$\int_{f = s} \rho. \vert z \vert^2\frac{ dx\wedge dy}{z^2}\wedge \frac{ d\bar{x}\wedge d\bar{y}}{\bar{z}^2}$$
commence par un terme non nul en \ $ (Log\vert s\vert )^2$, puisque l'on a sur \ $ \lbrace f = s \rbrace$
$$ 4.z.dx\wedge dy\wedge dz \big/ df = \frac{ dx\wedge dy}{z^2} .$$
En rempla{\c c}ant la fonction de troncature \ $\rho$ \ par la fonction caract\'eristique du polydisque \ $ (\vert x \vert \leq 1) \cap( \vert y \vert \leq 2) $ (voir [B.91] pour justifier que cela ne modifiera pas le p\^ole triple cherch\'e) on est ramen\'e \`a \'etudier le premier terme du d\'eveloppement asymptotique de la fonction
$$ \varphi(s) = \int_{(\vert x \vert \leq 1 )\cap (\vert y \vert \leq 2)} \frac{dx\wedge d\bar{x}\wedge dy\wedge d\bar{y}}{\vert s - x^2(x^2 + y^2)\vert} .$$
Apr\`es le changement de variable \ $ x = u.s^{\frac{1}{4}} , y = t.u.s^{\frac{1}{4}} $ \ on trouve
$$  \varphi(s) = c.\int_{(\vert u\vert \leq \vert s \vert^{-\frac{1}{4}}) \cap (\vert t.u \vert \leq 2.\vert s \vert^{-\frac{1}{4}})} \frac{\vert u \vert^2 .du\wedge d\bar{u}\wedge dt \wedge d\bar{t}}{\vert 1 - u^4.(1 + t^2)\vert} $$
o\`u \ $c$ \ est une constante non nulle.

\noindent Pour \ $ u$ \ fix\'e, on a
$$\int_{(\vert t.u \vert \leq 2.\vert s \vert^{-\frac{1}{4}})} \frac{dt\wedge d\bar{t}}{\vert 1 - u^4.(1 + t^2)\vert} \simeq c'. \vert u \vert^{-4}.Log(\vert u.s^{\frac{1}{4}}\vert) $$
ce qui donne
$$ \varphi(s) \simeq c''.(Log\vert s \vert )^2 $$
avec \ $c'$ \ et \ $c"$ \ des constantes non nulles.

\noindent Ceci montre que l'on a bien un p\^ole triple.

\noindent Il est facile de donner un \'el\'ement non nul de \ $H^1_{\lbrace 0 \rbrace}(S, H^1(0)) $ \ dans cet exemple: Fixons \ $ 0 < \alpha \ll \varepsilon $. Pr\`es d'un point du cercle \ $ \lbrace x = z = 0\ ; \vert y \vert = \alpha \rbrace $, on peut choisir comme coordonn\'ees locales \ $ \xi = x.\sqrt{y^2 + x^2} , y , z $ \ de sorte que la fonction \ $f$ \ dans ces coordonn\'ees se r\'eduit \`a \ $ \xi^2 + z^4 .$ Posons 
$$ v : = 2x.d\xi - \xi.dx .$$
Alors la forme diff\'erentielle holomorphe \ $z.v$, v\'erifie \ $ d(z.v) = \frac{df}{f}\wedge z.v $ \ et elle induit une section globale, uniforme et non nulle\footnote{puisque le mon\^ome \ $z$ \ n'est pas dans l'id\'eal jacobien de \ $ \xi^2 + z^4 .$}du syst\`eme local \ $H^1(0)$ \ le long de ce cercle. Le lecteur se convaincra facilement qu'elle n'est pas prolongeable \`a l'origine et donne donc un \'el\'ement non nul de \ $H^1_{\lbrace 0 \rbrace}(S, H^1(0)) $.

\bigskip

\section*{Appendice I :  Description topologique de \ $H^n_{c\,\cap\,S}(0).$}

\noindent Soit \ $ \mathbb{H} \overset{exp(2i\pi.\Box)}{\longrightarrow} D^* $ \ le rev\^etement universel du disque point\'e, et posons
$$ \hat{X} : = X \underset{D}{\times} \mathbb{H} .$$
Pour chaque voisinage ouvert \ $\mathcal{U}$ \ de \ $ S\, \cap\, \partial X $ \ notons par \ $\hat{\mathcal{U}}$ \ l'image r\'eciproque sur \ $\hat{X}$ \ de \ $\mathcal{U}$. D\'efinissons la famille paracompactifiante \ $\Phi$ \ de ferm\'es de \ $\hat{X}$\  en posant
$$ G \in \Phi \quad {\rm ssi} \quad \exists \ \mathcal{U} \supset S \cap \partial X \quad / \quad G \cap \hat{\mathcal{U}} = \emptyset .$$
On remarquera que cette famille de supports est invariante par l'action de l'automorphisme de monodromie de \ $\hat{X}$ \ induite par l'automorphisme \\ $\zeta \rightarrow \zeta +1 $ \ du demi-plan \ $\mathbb{H}$.
Nous allons montrer la

\bigskip

\noindent {\bf Proposition.}

\bigskip

\noindent \textit{ On a un isomorphisme "naturel" d'espaces vectoriels monodromiques }
$$ \varepsilon : H^n_{c\,\cap\, S}(0) \rightarrow  H^n_{\Phi}(\hat{X}, \mathbb{C})_{=1}. $$

\bigskip

\noindent \textit{D\'emonstration.} Commen{\c c}ons par remarquer que l'on a un diff\'eomorphisme compatible \`a \ $f$\footnote{O\`u \ $f$ \ est d\'efinie sur \ $\hat{X}$ \ comme la compos\'ee de la projection sur \ $\mathbb{H}$ \ avec la fonction \ $exp(2i\pi.\Box)$.}, grace \`a J. Milnor [Mi.68]
$$ \hat{X} \rightarrow F \times \mathbb{H} $$
o\`u \ $F$ \ d\'enote la fibre de Milnor de \ $f$ \ \`a l'origine. En particulier on a un isomorphisme "naturel" d'espaces vectoriels monodromiques \ $ H^n(\hat{X}, \mathbb{C}) \simeq H^n(F, \mathbb{C}) $ \ qui donne, pour la valeur propre \ $1$ \ l'isomorphisme monodromique  $ H^n(\hat{X}, \mathbb{C})_{=1} \simeq H^n(0).$

\noindent D\'efinissons \ $\varepsilon$. Soit donc \ $ w \in \Gamma_{c\,\cap\,S}(Y, \mathcal{E}^n(k))\,\cap \,Ker\,\delta $ \ et d\'efinissons sur \ $\hat{X}$ \ la \ $n-$forme \ $\mathcal{C}^{\infty} $
$$ W : = \sum_{j=0}^{k-1} \frac{(-\zeta)^j}{j!}.p^*(w_{k-j})$$
o\`u \ $ p : \hat{X} \rightarrow X $ \ est la projection et o\`u \ $\zeta \in \mathbb{H}.$

\noindent On a \ $dW = 0$ \ car on a suppos\'e que \ $\delta w = 0 $, et on a sur \ $\hat{X}$ \ l'\'egalit\'e \ $ p^*(\frac{df}{f}) = d\zeta .$ De plus, le fait que le support de \ $w$ \ rencontre \ $S$ \ suivant un compact montre que \ $W$ \ est \`a support dans \ $\Phi$. On en d\'eduit facilement que l'application \ $\varepsilon$ \ est bien d\'efinie en posant 
$$ \varepsilon ([w]) = [W] .$$
Pour v\'erifier que \ $\varepsilon$ \ commute aux monodromies, il nous suffit alors d'\'etablir la relation
$$ \varepsilon(\mathcal{N}(w)) = - \frac{\partial W}{\partial\zeta} \qquad \qquad (@) $$
puisque d'apr\`es la formule de Taylor, on a pour un polyn\^ome \ $P$
$$exp(\frac{\partial}{\partial\zeta})(P)(\zeta) = P(\zeta + 1).$$
La v\'erification de (@) est imm\'ediate. On en conclut que la classe \ $[W]$ \ est bien dans \ $ H^n_{\Phi}(\hat{X}, \mathbb{C})_{=1}.$

\noindent Montrons maintenant l'injectivit\'e de \ $\varepsilon$.  Il suffit de consid\'erer le cas d'une classe \ $ [w]$ \ invariante par la monodromie qui donne \ $0$ \ dans \ $ H^n_{\Phi}(\hat{X}, \mathbb{C})_{=1}$,  grace \`a l'aspect monodromique de \ $\varepsilon$. Mais dans ce cas, l'image par l'application d'oubli de support de \ $[W]$ \ dans \ $H^n(\hat{X}, \mathbb{C})_{=1} \simeq H^n(0) $ \ montre que \ $[w] $ \ est dans le noyau de \ $can_{c\,\cap\,S}$. On a donc seulement \`a traiter le cas o\`u \ $ [w] = i([\alpha]) = [\tilde{\alpha}\wedge d\chi ]$ \ avec \ $\alpha \in H^1_{\lbrace 0 \rbrace}(S, H^{n-1}(0))$ \ d'apr\`es le th\'eor\`eme 1\footnote{L'adaptation de ceci pour \ $n = 2$ \ ne pr\'esente pas de difficult\'e et elle est laiss\'ee au lecteur.}

\noindent  Soit donc une telle classe invariante \ $[w] = [\tilde{\alpha}\wedge d\chi ]$ \ v\'erifiant \ $\varepsilon([w]) = 0 $. Ceci signifie que l'on peut trouver une \ $(n-1)-$forme \ $v \in \mathcal{C}_{\Phi}^{\infty}(\hat{X}) $ \ v\'erifiant \ $dv = p^*(w) $. Fixons \ $\mathcal{U}$ \ un voisinage ouvert de \ $ S\, \cap\, \partial X $ \ assez petit pour que l'on ait \ $ Supp(v)\, \cap\, \hat{\mathcal{U}} = \emptyset $. 

\noindent Supposons maintenant que \ $[w] \not= 0 $ \ dans \ $H^n_{c\,\cap\,S}(0)$ \ et \ $ n \geq 3$. Comme on sait que \ $ [w] = i[\alpha] $ \ le choix de la fonction \ $\chi $\footnote{voir la d\'efinition de l'application \ $i$ \ au (4.9).} permet de supposer que l'on a \ $ Supp(\hat{w})\,\cap \,S \subset \mathcal{U} $ \ o\`u \ $\hat{w} = \tilde{\alpha}\wedge d\chi $ \ est \'egalement un repr\'esentant de la classe \ $[w]$.

\noindent  La forme sesquilin\'eaire \ $\mathcal{H}$ \ \'etant non d\'eg\'en\'er\'ee, on peut trouver une classe \ $[w'] \in H^n(0) $ \ telle que l'on ait
$$ P_2 \big(\lambda = 0, \int_X \vert f \vert^{2\lambda} \frac{df}{f}\wedge \frac{\bar{df}}{\bar{f}}\wedge w \wedge \bar{w'}_k \big) \not= 0. $$
Soit \ $\sigma \in \mathcal{C}_c^{\infty}(\mathcal{U}) $ \ valant identiquement \ $1$ \ au voisinage de \ $ S\,\cap \, Supp(w ) $. Remarquons que, quitte \`a choisir \ $\mathcal{U}$ \ plus petit, on peut toujours supposer que \ $ Supp(\sigma)\,\cap\,\mathcal{U} = \emptyset $. La formule de Stokes et le prolongement analytique donne alors
$$ P_1 \big(\lambda = 0, \int_X \vert f \vert^{2\lambda} d\sigma \wedge \frac{\bar{df}}{\bar{f}}\wedge w \wedge \bar{w'}_k \big) \not= 0. $$
ce qui signifie que dans le d\'eveloppement asymptotique en \ $ s = 0 $ \ de l'int\'egrale fibre
$$ s\frac{\partial}{\partial s} \int_{f = s} \sigma.w\wedge\bar{w'}_k $$
le terme constant est non nul. Mais ceci est absurde car \ $p^*(w) = dv $ \ et \ $Supp(\sigma) \subset \mathcal{U} $ \ montrent que l' int\'egrale
$$  \int_{f = s} \sigma.w\wedge\bar{w'}_k  = -  \int_{f = s} d\sigma \wedge v \wedge\bar{w'}_k $$
est \ $\mathcal{C}^{\infty}$ \ pr\`es de \ $ s = 0$\footnote{Plus exactement, la partie dans \ $ \sum_{j = 0}^n \mathbb{C}[[s,\bar{s}]].(Log\vert s \vert^2)^j $ \ de son d\'eveloppement asymptotique en \ $s = 0$ \ est en fait dans \ $ \mathbb{C}[[s,\bar{s}]].$}  puisque \ $ Supp(v)\,\cap \hat{\mathcal{U}} = \emptyset . $

\noindent Terminons la preuve de l'injectivit\'e de \ $\varepsilon$ \ pour \ $ n = 2$. Il suffit de voir que \ $j(H^1(S^*, \mathbb{C})) $ \ s'injecte dans \ $ H^2_{\Phi}(\hat{X}, \mathbb{C})_{=1}$. Si \ $ \varphi$ \ est $d-$ferm\'ee au voisinage de \ $S^*$ \ et \ $\chi \equiv 1$ \ pr\`es de \ $S^*$ \ et s'annule d\`es que l'on s'\'eloigne, on a \ $ j[\varphi] = [ d\chi \wedge \varphi] $ \ et \ $ \varepsilon[d\chi \wedge \varphi] = [ p^*(d\chi \wedge \varphi)] .$ Mais si \ $ p^*(d\chi \wedge \varphi) = d\lambda $ \ o\`u \ $ \lambda \in \mathcal{C}^{\infty}_{\Phi}(\hat{X})^{n-1} $ \ on aura \ $ p^*(\chi.\varphi) - \lambda $ \ qui sera $d-$ferm\'ee et localement $d-$exacte pr\`es de \ $ \partial S $. La classe qu'elle d\'efinie dans \ $ H^{n-1}(F, \mathbb{C}) $ \ est donc nulle\footnote{Car une section globale du faisceau  \ $H^{n-1}(0) $ \ qui est nulle pr\`es de \ $\partial S$ \ est nulle.}. On en conclut que \ $\varphi$ \ est globalement $d-$exacte au voisinage de \ $\partial S$ \ et donc induit la classe nulle de \ $H^1(S^*, \mathbb{C}).$

\bigskip

\noindent Montrons maintenant la surjectivit\'e pour \ $ n \geq 2$, en montrant par r\'ecurrence sur \ $ k \geq 1$, que si \ $ [W] \in H^n_{\Phi}(\hat{X}, \mathbb{C})_{=1}$ \ v\'erifie \ $ (T - 1)^k [W] = 0 $, elle est bien dans l'image de \ $\varepsilon$.

\noindent Pour \ $k = 1$, on a \ $ T[W] = [W] $ \ et on peut trouver une $(n-1)-$forme   \ $V \in \mathcal{C}^{\infty}_{\Phi}(\hat{X}) $ \ qui v\'erifie \ $ -\frac{\partial}{\partial \zeta} W = dV .$ On \'ecrit alors\footnote{C'est toujours possible : on utilise ici le  complexe \ $(\mathcal{E}^{\bullet}[Log\, f], d^{\bullet})$ \ o\`u \ $-2i\pi. \zeta = Log\,f $ \
 sur \ $\hat{X}$, pour calculer la partie spectrale pour la valeur propre 1 de la monodromie des cycles \'evanescents de \ $f$.} $ V =  -\frac{\partial}{\partial \zeta} U $ \ avec \ $ U \in \mathcal{C}^{\infty}_{\Phi}(\hat{X})^{n-1}$. Alors \ $ W - dU$ \ est un repr\'esentant $T-$invariant de la classe \ $[W]$ \ et il existe une forme \ $ w \in \mathcal{C}^{\infty}_{c\,\cap\,S}(X \setminus f^{-1}(0))^n $ \ v\'erifiant \ $ dw = 0 $ \ et \ $p^*(w) = W - dU .$
On conclut en utilisant le th\'eor\`eme de A.Grothendieck [Gr.65] pour trouver une $n-$forme semi-m\'eromorphe $d-$ferm\'ee \ $\tilde{w}$ \ \`a p\^oles dans \ $f^{-1}(0)$ \ et \`a support \ $c\,\cap\,S$ \ induisant la m\^eme classe que \ $w$ \ dans \ $ H^n_{c\,\cap\,S}(X \setminus f^{-1}(0), \mathbb{C}).$

\bigskip

\noindent Supposons maintenant le r\'esultat montr\'e pour \ $ k \leq k_0$ \ avec \ $ k_0 \geq 1$ \ et montrons-le pour \ $ k = k_0 +1.$ Alors d'apr\`es l'hypoth\`ese de r\'ecurrence, on peut \'ecrire \ $  -\frac{\partial}{\partial \zeta} W = \tilde{W} + dV $ \ o\`u \ $[ \tilde{W}]$ \ est dans l'image de \ $ \varepsilon$ \ et o\`u \ $ V \in \mathcal{C}^{\infty}_{\Phi}(\hat{X})^{n-1}. $ On \'ecrit \`a nouveau \ $ V = -\frac{\partial}{\partial \zeta}  U $ \ avec \ $ U \in \mathcal{C}^{\infty}_{\Phi}(\hat{X})^{n-1}$ \ et arrive \`a \ $  -\frac{\partial}{\partial \zeta}(W - dU) = \tilde{W} = \varepsilon(\tilde{w}) .$ Mais alors la $n-$forme
$$\int_0^{\zeta} \tilde{W}(z).dz +  (W - dU)$$
 est invariante par la monodromie. Elle est donc image r\'eciproque d'une forme \ $\mathcal{C}^{\infty}$ \ et \ $d-$ferm\'ee sur \ $ X \setminus f^{-1}(0) $ \ \`a support dans \ $ c\,\cap\,S$. On conclut alors facilement gr\^ace au m\^eme argument que  dans le cas \ $ k = 1$, en int\'egrant terme \`a terme la formule donnant \ $ \varepsilon(\tilde{w})$. $\hfill \blacksquare$

\bigskip

\newpage

\section*{Appendice II : Interpr\'etation en termes de complexes de~cycles proches et \'evanescents.\\  {\em  Par Claude Sabbah.}}

\bigskip
Soit $f:(\mathbb{C}^{n+1},0)\to(\mathbb{C},0)$ un germe de fonction holomorphe ($n\geq1$). Soit $\psi_f\QQ$ (resp.~$\phi_f\QQ$) le complexe des cycles proches (resp.~\'evanescents) de $f$ pour le faisceau constant $\QQ$ (voir [D. 73] et  aussi [K.S. 90]). Il est muni d'un op\'erateur de monodromie $T$. \`A~un d\'ecalage de $n$ pr\`es, c'est un faisceau pervers (voir [Br. 86] et aussi [K.S. 90]).
Il admet une d\'ecomposition suivant les valeurs propres
$$\psi_f\QQ=\oplus_{u\in\QQ\cap[0,1[}\psi_{f,\exp(-2i\pi u)}\QQ\quad \text{resp.~}\phi_f\QQ=\oplus_{u\in\QQ\cap[0,1[}\phi_{f,\exp(-2i\pi u)}\QQ
$$
dans la cat\'egorie des faisceaux pervers (voir \emph{loc.\kern2pt cit.}). L'op\'erateur $T-\exp(-2i\pi u)\Id$ est nilpotent sur $\psi_{f,\exp(-2i\pi u)}\QQ$ (resp.~$\phi_{f,\exp(-2i\pi u)}\QQ$).

Pour chaque $u$, on dispose d'un morphisme ``canonique''
$$
\can:\psi_{f,\exp(-2i\pi u)}\QQ\to\phi_{f,\exp(-2i\pi u)}\QQ
$$
et d'un morphisme ``variation''
$$
\var:\phi_{f,\exp(-2i\pi u)}\QQ\to\psi_{f,\exp(-2i\pi u)}\QQ
$$
dont les compos\'es $\var\circ\can$ et $\can\circ\var$ sont \'egaux \`a $T-\Id$ sur les complexes correspondants. En un point~$x$ de~$f^{-1}(0)$, le germe en $x$ du faisceau de cohomologie $H^p(u):=\cH^p\psi_{f,\exp(-2i\pi u)}\QQ$ est l'espace propre g\'en\'eralis\'e $H^p(F_x,\QQ)_{\exp(-2i\pi u)}$ de la monodromie, si $F_x$ est la fibre de Milnor (d'un bon repr\'esentant) de $f$ en~$x$.

Lorsque $u\neq0$, ces deux morphismes sont des isomorphismes, et on identifie $\psi_{f,\exp(-2i\pi u)}\QQ$ et $\phi_{f,\exp(-2i\pi u)}\QQ$ \emph{via} $\can$, de sorte que $\var$ devient \'egal \`a $T-\Id$.

Soit $i_{f^{-1}(0)}:\{f=0\}\hookrightarrow \CC^{n+1}$ l'inclusion. Lorsque $u=0$, on a un triangle de la cat\'egorie d\'eriv\'ee des complexes born\'es \`a cohomologie constructible sur $f^{-1}(0)$:
\[\tag{$*$}
i_{f=0}^{-1}\QQ\to\psi_{f,1}\QQ\to\phi_{f,1}\QQ\To{+1}
\]
qui dit que les faisceaux de cohomologie de $\phi_{f,1}\QQ$ sont les faisceaux de cohomologie r\'eduite de $\psi_{f,1}\QQ$.

\smallskip\noindent
\textit{Hypoth\`ese}.~Fixons une valeur propre $\exp(-2i\pi u)$ ($u\in[0,1[\cap\QQ$). Nous faisons l'hypoth\`ese que le complexe de faisceaux $\phi_{f,\exp(-2i\pi u)}\QQ$ est \`a support sur une \emph{courbe} $(S,0)$. Nous noterons
$$
i:\{0\} \hookrightarrow S,\quad j:S\setminus \{0\}\hookrightarrow S,\quad i_S:S \hookrightarrow \{f=0\},
$$
les inclusions. Pour simplifier, nous noterons $\psi_{f,\exp(-2i\pi u)},\phi_{f,\exp(-2i\pi u)}$ les complexes $\psi_{f,\exp(-2i\pi u)}\QQ,\phi_{f,\exp(-2i\pi u)}\QQ$.

\smallskip\noindent
\textbf{1.}~La perversit\'e de $\phi_{f,\exp(-2i\pi u)}$ implique que $\phi_{f,\exp(-2i\pi u)}$ n'a de cohomologie qu'en degr\'e $n$ et $n-1$. Le faisceau $\cH^n\phi_{f,\exp(-2i\pi u)}$, not\'e $\cH^n$, est \`a support l'origine, et le faisceau $\cH^{n-1}\phi_{f,\exp(-2i\pi u)}$, not\'e $\cH^{n-1}$, est \`a support sur $S$.

On en d\'eduit (sans hypoth\`ese si $u\neq0$ et sous l'hypoth\`ese $n\geq2$ si $u=0$, en utilisant le triangle $(*)$)  que, 
\begin{itemize}
\item
$\cH^n$ est la partie correspondant \`a la valeur propre $\exp(-2i\pi u)$ du $n$-i\`eme groupe de cohomologie de la fibre de Milnor de $f$ en $0$,
\item
$\cH^{n-1}$ est le faisceau analogue pour la cohomologie en degr\'e $n-1$.
\end{itemize}

Comme le faisceau constant $\QQ$ sur $\CC^{n+1}$ est auto-dual (par dualit\'e de Poincar\'e-Verdier) \`a d\'ecalage pr\`es de $2(n+1)$, on d\'eduit, puisque le foncteur des cycles \'evanescents commute \`a la dualit\'e \`a un d\'ecalage de $1$ pr\`es (voir [Br. 86] et [K.S. 90]), que $\phi_{f,\exp(-2i\pi u)}$ est dual de $\phi_{f,\exp(2i\pi u)}$ \`a un d\'ecalage de $2n$ pr\`es.
D'un autre c\^ot\'e, le foncteur $i^{-1}$ est transform\'e en $i^!=R\Gamma_{0}$ \`a un d\'ecalage de $2$ pr\`es.

Le complexe $i^{-1}\phi_{f,\exp(-2i\pi u)}$  est un complexe d'espaces vectoriels ayant deux groupes de cohomologie $\cH^n$ et $\cH^{n-1}_0:=i^{-1}\cH^{n-1}$ (dans les degr\'es indiqu\'es).

Par dualit\'e de Poincar\'e-Verdier, on en d\'eduit que $i^!\phi_{f,\exp(2i\pi u)}$ n'a que deux groupes de cohomologie, qui sont
$$
\cH^n(i^!\phi_{f,\exp(2i\pi u)})\simeq{\cH^n}^\vee\quad\text{et}\quad\cH^{n+1}(i^!\phi_{f,\exp(2i\pi u)})\simeq{\cH^{n-1}_0}^\vee
$$
(o\`u ${}^\vee$ d\'esigne le dual des espaces vectoriels). En particulier, on a un accouplement non d\'eg\'en\'er\'e $\cH^n(i^!\phi_{f,\exp(2i\pi u)}) \otimes\cH^n\to\QQ$.

\smallskip\noindent
\textbf{2.}~Puisque $\phi_{f,\exp(-2i\pi u)}$ n'a que deux faisceaux de cohomologie, on a un triangle
\[
\cH^{n-1}\phi_{f,\exp(-2i\pi u)}[-(n-1)]\to\phi_{f,\exp(-2i\pi u)}\to\cH^n[-n]\To{+1}
\]
et donc, en appliquant $i^!$, on a un triangle correspondant, qui donne une suite exacte longue
\begin{multline*}
0\rightarrow H_{0}^{1}(S,\cH^{n-1})\rightarrow\cH^n(i^!\phi_{f,\exp(-2i\pi u)})\rightarrow\cH^n\rightarrow H_{0}^{2}(S,\cH^{n-1})\\
\rightarrow\cH^{n+1}(i^!\phi_{f,\exp(-2i\pi u)})\rightarrow0
\end{multline*}
avec $0$ aux deux bouts puisque $\cH^n$ est \`a support l'origine. Pour la m\^eme raison, le morphisme du milieu $\can:\cH^n(i^!\phi_{f,\exp(-2i\pi u)})\to\cH^n$ est le $H^n$ du morphisme de complexes $i^!\phi_{f,\exp(-2i\pi u)}\to i^{-1}\phi_{f,\exp(-2i\pi u)}$ (oubli des supports).

\smallskip\noindent
\textbf{3.~Interpr\'etation de $\cH^n(i^!\phi_{f,\exp(-2i\pi u)})$.}
J'interpr\`ete le groupe $H^n_{c\cap S}(u)$ comme le groupe de cohomologie de degr\'e $n$ du complexe $i^!i_S^{-1}\psi_{f,\exp(-2i\pi u)}$.

Si $u\neq0$, on a donc $\cH^n(i^!\phi_{f,\exp(-2i\pi u)})=H^n_{c\cap S}(u)$.

Supposons maintenant que $u=0$. Apr\`es application de $i_S^{-1}$, le triangle $(*)$ devient, puisque $\phi_{f,1}$ est \`a support dans $S$ (donc $i_S^{-1}\phi_{f,1}=\phi_{f,1}$)
\[
\QQ_S\to i_S^{-1}\psi_{f,1}\to\phi_{f,1}\To{+1}
\]
et, apr\`es application de $i^!$,
\[
i^!\QQ_S\to i^!i_S^{-1}\psi_{f,1}\to i^!\phi_{f,1}\To{+1}.
\]
La cohomologie du complexe $i^!\QQ_S=R\Gamma_{0}\QQ_S$ n'existe qu'en degr\'e $2$. Donc
\[
\cH^n(i^!\phi_{f,1})=
\begin{cases}
H^n_{c\cap S}(0)&\text{si }n\geq3,\\
H^n_{c\cap S}(0)/\text{image}\big(H^2_0(\QQ_S)\big)&\text{si }n=2,
\end{cases}
\]
o\`u $H^2_0(\QQ_S)\simeq H^1(S^*,\QQ)$.

\newpage

\section{Bibliographie.}

\begin{itemize}

\item{[A.V.G]} Arnold, V., Varchenko, A. et Goussein-Zad\'e, S. \textit{Singularit\'es des applications diff\'erentiables}, vol.2 Monodromie et comportement asymptotique des int\'egrales, Moscou 1986
(traduction fran{\c c}aise, Ed. Mir).

\item{[B.84.a]} Barlet, D.  \textit{Contribution effective de la monodromie aux d\'eveloppements asymptotiques}, Ann. Scient. ENS 4-i\`eme s\'erie, 17 (1984) p.239-315.

\item{[B.84.b]} Barlet, D.  \textit{Contribution du cup-produit de la fibre de Milnor aux p\^oles de \ $\vert f \vert^{2\lambda}$ }, Ann. Inst. Fourier (Grenoble), 34 (1984) p.75-107.

\item{[B.85]} Barlet, D.  \textit{La forme hermitienne canonique sur la cohomologie de la fibre de Milnor d'une hypersurface \`a singularit\'e isol\'ee}, Invent. Math. 81 (1985) p.115-153 .

\item{[B.90]} Barlet, D.   \textit{La forme hermitienne canonique pour une singularit\'e presque isol\'ee}, in Complex Analysis (K.Diederich eds) Vieweg,Wuppertal (1990) p.20-28  .

\item{[B.91]} Barlet, D.  \textit{Emm\^elements de strates cons\'ecutives pour les cycles \'evanescents}, Ann. Scient. ENS 4-i\`eme s\'erie  24 (1991) p.401-506.

\item{[B.97]} Barlet, D.   \textit{La variation pour une hypersurface \`a singularit\'e isol\'ee relativement \`a la valeur propre  1 }, Revue de l'Inst. E. Cartan (Nancy) 15 (1997) p.1-29  .

\item{[B.02]} Barlet, D.   \textit{Singularit\'es r\'eelles isol\'ees et d\'eveloppements asymptotiques d'int\'egrales oscillantes}, in S\'eminaire et Congr\`es 9 (Actes des journ\'ees math\'ematiques \`a la m\'emoire de J. Leray) Guillop\'e, L. et Robert, D. \'editeurs, Soci\'et\'e Math\'ematique de France (2004), p. 25-50.
\item{[B.03]}  Barlet, D.   \textit{Interaction de strates cons\'ecutives  II }, Publ. du RIMS Kyoto Univ., vol. 41 (2005), p. 139-173.
\item{[B.04]}  Barlet, D.   \textit{Sur certaines singularit\'es non isol\'ees d'hypersurfaces I},
preprint Institut E.Cartan (Nancy) 2004/$n^{\circ}$ 03.

\item{[B.M.89]} Barlet, D. et Maire, H.-M.   \textit{Asymptotic expansion of complex integrals via Mellin transform}, Journ. Funct. Anal. 83 (1989) p.233-257.

\item{[Br. 86]} Brylinski, J.-L. \textit{Transformations canoniques, dualit\'e projective, th\'eorie de {Lefschetz}, transformation de {Fourier} et sommes trigonom\'etriques}, in {\em G\'eom\'etrie et ana\-lyse microlocales}, Ast{\'e}risque, vol.~140-141, Soci{\'e}t{\'e} Math{\'e}matique de  France, 1986, p.3-134.

\item{[D. 73]} Deligne, P. \textit{Le formalisme des cycles \'evanescents (expos\'es 13 et 14)}, in {\em SGA~7~II}, Lect. Notes in Math., vol.340, Springer-Verlag, 1973, p.82-173.

\item{[Go.58]} Godement, R. \textit{Th\'eorie des faisceaux}, Hermann (1958) Nancago.

\item{[Gr.65]} Grothendieck, A. \textit{On the De Rham cohomology of algebraic varieties}, Publ. Math. IHES, 29 (1966) p.96-103.

\item{[K.76]} Kashiwara, M. \textit{b-Function and Holonomic Systems, Rationality of Roots of b-Functions}, Invent. Math. 38 (1976) p.33-53.

\item{[K.84]} Kashiwara, M. \textit{The Riemann-Hilbert Problem for Holonomic Systems}, Publ. RIMS Kyoto Univ. 20 (1984) p.319-365.

\item{[K. S. 90]} Kashiwara, M. and Schapira, P. \textit{Sheaves on Manifolds}, {\em Grundlehren der mathematischen Wissenschaften}, vol.292, Springer-Verlag, 1990.

\item{[Le.59]} Leray, J.  \textit{Le probl\`eme de Cauchy III }, Bull.Soc. Math. France 87 (1959)  p.81-180  .

\item{[Lo.86]} Loeser, F. \textit{A propos de la forme hermitienne canonique d'une singularit\'e isol\'ee d'hypersurface}, Bull. Soc. Math. France, 114 (1986) p.385-392.

\item{[Ma.74]} Malgrange, B.  \textit{Int\'egrales asymptotiques et monodromie}, Ann. Scient. ENS 4-i\`eme serie   7 (1974)  p.405-430  .

\item{[Mi.68] } Milnor, J.  \textit{Singular Points of Complex Hypersurfaces}, Ann. of Math. Studies  61 (1968) Princeton  .

\end{itemize}

  \end{document}